\documentclass[12pt]{article}     

\usepackage{graphicx}
\usepackage{mathptmx}
\usepackage{bm}
\usepackage{amsmath}    
\usepackage[percent]{overpic}
\usepackage{multirow}
\usepackage[caption=false]{subfig}
\usepackage{lscape}
\usepackage{epsfig}
\usepackage{latexsym}
\newcommand{\dss}{\displaystyle}
\newcommand{\tss}{\textstyle}
\newtheorem{remark}{Remark}

\begin{document}

\title{Backward bifurcation underlies rich dynamics \\ 
in simple disease models
\thanks{This work was supported by the Natural Sciences and Engineering
Research Council of Canada (NSERC).}}

\author{
Wenjing Zhang
\thanks{Wenjing Zhang, email: zhang.wenjing14@gmail.com} 
\and Pei Yu, \thanks{Pei Yu, email: pyu@uwo.ca} 
\thanks{Applied Mathematics, Western University, 
London, Ontario, Canada N6A 5B7}
\and Lindi M. Wahl \thanks{Lindi M. Wahl, email: lwahl@uwo.ca}
\thanks{Applied Mathematics, Western University, 
London, Ontario, Canada N6A 5B7}}


\maketitle

\begin{abstract}
In this paper, dynamical systems theory and bifurcation theory 
are applied to investigate the rich dynamical behaviours observed
in three simple disease models.  The 2- and 3-dimensional models
we investigate have arisen in previous investigations of epidemiology,
in-host disease, and autoimmunity.  These closely related models
display interesting dynamical behaviors 
including bistability, recurrence, and regular oscillations, each
of which has possible clinical or public health implications.
In this contribution we elucidate the key role of backward bifurcation
in the parameter regimes leading to the behaviors of interest. 
We demonstrate that backward bifurcation facilitates the
appearance of Hopf bifurcations, and the varied dynamical behaviors
are then determined by the properties of the
Hopf bifurcation(s), including their location and direction.
A Maple program developed earlier is implemented to determine the 
stability of limit cycles bifurcating from the Hopf bifurcation.
Numerical simulations are presented to illustrate phenomena of interest
such as bistability, recurrence and oscillation.  We also
discuss the physical motivations for the models
and the clinical implications of the resulting dynamics.
\end{abstract}

\section{Keywords}
Convex incidence rate, backward bifurcation,
supercritical and subcritical Hopf bifurcations, recurrence, bistability,
disease model

\section{Mathematics Subject Classification (2000)}
$$34C \bullet 34D \bullet 92B \bullet 92D$$

\section{Introduction}
In the mathematical modelling of epidemic diseases, 
the fate of the disease can be predicted through the 
uninfected and infected equilibria and their stability.
The basic reproduction number, $R_0$, represents the average number 
of new infectives introduced into an otherwise disease-free system by 
a single infective, 
and is usually chosen as the bifurcation parameter.
If the model involves a forward bifurcation, the uninfected equilibrium is
in general globally asymptotically stable~\cite{Korobeinikov2005}, 
characterized by $R_0<1$, and infection fails to invade in this parameter 
regime. 
The threshold $R_0=1$ defines a bifurcation (or critical) point, 
and when $R_0>1$, a stable infected equilibrium emerges. 
This simple exchange of stability implies that complex dynamics 
will not typically occur in forward bifurcation.  

In contrast, backward bifurcation describes a scenario
in which a turning point of the infected equilibrium exists in a region
where all state variables are positive, and $R_0 < 1$.
This induces multiple infected equilibria, disrupting
the global stability of the uninfected equilibrium.
Multiple stable states (e.g., bistability) may likewise appear in
\cite{dushoff1998backwards,Blayneh2010,arino2003global},
and Yu et al. (submitted for publication). 
Instead of converging globally to the uninfected equilibrium when $R_0<1$, 
the solution may approach an infected equilibrium,
depending on initial conditions.

In practice, the phenomenon of backward bifurcation
gives rise to new challenges in disease control,
since reducing $R_0$ such that 
$R_0<1$ is not sufficient to eliminate the disease
\cite{hadeler1997backward,Brauer2004418}.
Instead, $R_0$ needs to be reduced past the critical value given
by the turning point~\cite{hadeler1997backward}, 
since the result in Yu et al. (submitted for publication) 
shows that the
uninfected equilibrium in backward bifurcation is globally stable 
if $R_0$ is smaller than the turning point.
Furthermore, an infective outbreak or catastrophe may occur
if $R_0$ increases and crosses unity, while
the upper branch of the infected equilibrium remains stable
\cite{dushoff1998backwards,GómezAcevedo2005101,ZWY2013,ZWY2014}.
In addition, oscillation or even recurrent phenomena may occur
if uninfected and infected
equilibria coexist in a parameter range, and both are unstable~\cite{ZWY2013,ZWY2014}. 
\cite{hadeler1997backward} predicted oscillations arising
from backward bifurcation,
and \cite{Brauer2004418} pointed out that the unstable infected 
equilibrium ``commonly arises from Hopf bifurcation'', 
but did not demonstrate oscillations.

Several mechanisms leading to backward bifurcation have
been proposed, such as partially effective vaccination programs
\cite{Brauer2004418,arino2003global}, educational influence on
infectives' behavior~\cite{hadeler1997backward},
the interaction among multi-group models
\cite{CastilloChavez1989,CastilloChavez1989327,huang1992stability}
and multiple stages of infection~\cite{simon1992reproduction}.
In this study, we will investigate the emergence of backward 
bifurcation in three simple disease models which have arisen
in the study of epidemiology, in-host disease and autoimmunity.
In each case, we find that backward bifurcation facilitates the
emergence of Hopf bifurcation(s), and Hopf bifurcation in turn underlies
a range of complex and clinically relevant dynamical behaviors. 

A central theme in our investigation is the role of the incidence
rate in the epidemiological and in-host disease models.
The incidence rate describes the speed at which an infection spreads;
it denotes the rate at which susceptibles become infectives.
Under the assumptions of mass action, 
incidence is written as the product of the infection force
and the number of susceptibles. 
For example, if $S$ and $I$ denote the 
susceptible and infective population size respectively,
a bilinear incidence rate, $f(S, I) = \beta S I$ 
(where $\beta$ is a positive constant),
is linear in each of the state variables: $S$ and $I$.

The possibility of saturation effects
\cite{capasso1978generalization,Brown199510} has motivated the modification 
of the incidence rate from bilinear to nonlinear.
Saturation occurs when the number of susceptible contacts per infective
drops off as the proportion of infectives increases. 
A nonlinear incidence rate, therefore, typically increases sublinearly
with respect to the growth of the infective population,
and may finally reach an upper bound.
The development of nonlinear incidence 
was first investigated in the form $\beta I^p S^q$,
where $\beta$, $p$, and $q$ are positive constants
\cite{Liu1986,Liu1987,hethcote1989epidemiological,
hethcote1991some,DerrickDriessche1993,li1995global}.
Other forms of nonlinear incidence have also been
analysed, such as $k I^p S/(1+\alpha I^l)$~\cite{Liu1986}, and 
$k S \ln(1+vP/k)$~\cite{briggs1995dynamics}.

Since the nonlinear incidence functions described above
were often developed to incorporate saturation effects, 
these functions are typically concave at realistic parameter values.
\cite{Korobeinikov2005} used this feature to 
derive general results for disease models with concave incidence.
They proved that standard epidemiological models with 
concave incidence functions will have globally asymptotically stable 
uninfected and infected equilibria for $R_0<1$ and $R_0>1$, respectively.

More specifically, denoting the incidence rate function 
as $f(S,\,I,\,N)$, where $N$ is the population size, the classical SIRS model 
considered in~\cite{Korobeinikov2005} 
takes the form
\begin{equation} 
\frac{{\mathrm d}S}{{\mathrm d} t} =  \mu N - f(S,\,I,\,N) - \mu S + \alpha R,  \quad  
\frac{{\mathrm d}I}{{\mathrm d} t} =  f(S,\,I,\,N) - (\delta + \mu) I, \quad  
\frac{{\mathrm d}R}{{\mathrm d} t}=  \delta I -\alpha R - \mu R, 
\label{Paper3_Eq1} 
\end{equation}
where $\mu$, $\delta$, and $\alpha$ represent the birth/death rate, 
the recovery rate and the loss of immunity rate, respectively.
When $\alpha=0$, system~(\ref{Paper3_Eq1}) becomes an SIR model.
Assuming that the total population size is constant, that is, $N=S+I+R$,
the above system can be reduced to a 2-dimensional model: 
\begin{equation}
\frac{{\mathrm d}S}{{\mathrm d} t}=(\alpha+\mu)N - f(S,\,I,\,N) - \alpha I - (\alpha+\mu)S, \quad
\frac{{\mathrm d}I}{{\mathrm d} t}= f(S,\,I,\,N) - (\delta + \mu) I.
\label{Paper3_Eq2}
\end{equation}
Moreover, it is assumed in~\cite{Korobeinikov2005} that the
function $f(S,\,I,\,N)$, 
denoting the incidence rate,
satisfies the following three conditions: 
\begin{subequations}
\begin{align}
& f(S,\,0,\,N)=f(0,\,I,\,N)=0, \label{Paper3_Eq3a}\\[0.5ex]
&\frac{\partial f(S,\,I,\,N)}{\partial I}>0, \quad 
\frac{\partial f(S,\,I,\,N)}{\partial S}>0, \quad
\forall \; S,\,I >0  \label{Paper3_Eq3b}\\[0.5ex]
&\frac{\partial^2 f(S,\,I,\,N)}{\partial I^2}\leq 0, \quad \forall \; S,\,I >0
\label{Paper3_Eq3c}.
\end{align}
\label{Paper3_Eq3}
\end{subequations}
The first two conditions (\ref{Paper3_Eq3a}) and (\ref{Paper3_Eq3b}) are 
necessary to ensure that the model is biologically meaningful. 
The third condition~(\ref{Paper3_Eq3c}) implies that 
the incidence rate $f(S,\,I,\,N)$,
is concave with respect to the number of infectives.
It is also assumed that 
$\frac{\partial f(S,\,I,\,N)}{\partial I}$ evaluated at the uninfected
equilibrium is proportional to the basic reproduction number $R_0$
\cite{vandenDriessche200229}, and thus should be a positive finite number
\cite{Korobeinikov2005}.
Korobeinikov and Maini
first considered $\dot{I}=0$, or $f(S,\,I,\,N) - (\delta + \mu) I=0$, 
and showed that forward bifurcation occurs in model~(\ref{Paper3_Eq2})
with a concave incidence function. They further proved that
the uninfected equilibrium $Q_0=(S_0,\,I_0)=(N,\,0)$ and
the infected equilibrium $\bar{Q}=(\bar{S},\,\bar{I})$ are 
globally asymptotically stable, when
$R_0=\frac{1}{\delta+\mu}\frac{\partial f(S_0,\,I_0,\,N)}{\partial I}<1$
and $R_0>1$, respectively. 

In the sections to follow, for an incidence rate function $f(S,\,I)$, satisfying
 (\ref{Paper3_Eq3a}) and (\ref{Paper3_Eq3b}), 
we define $f(S,\,I)$ 
as concave, if it satisfies~(\ref{Paper3_Eq3c});
as convex, if $\frac{\partial^2 f(S,\,I)}{\partial I^2} > 0$,
$\forall \; I>0$; 
and as convex-concave, if there exist $0<I_1<I_2\leq +\infty$, 
such that $\frac{\partial f(S,\,I)}{\partial I}>0$, 
$\forall \; I \in (0,I_2)$, and $\frac{\partial^2 f(S,\,I)}{\partial I^2}> 0$,
$\forall \; I \in (0,I_1)$, $\frac{\partial^2 f(S,\,I)}{\partial I^2} = 0$,
for $I=I_1$, $\frac{\partial^2 f(S,\,I)}{\partial I^2} < 0$,
$\forall \; I \in (I_1,I_2)$.

Several models closely related to (\ref{Paper3_Eq2}) have been previously studied.
For example, by adding a saturating treatment term to 
model~(\ref{Paper3_Eq2}) with a concave incidence rate, 
\cite{Zhou2012} showed that this model 
may yield backward bifurcation and Hopf bifurcation.
With an even more sophisticated nonlinear incidence rate function:
$k I^p S/(1+\alpha I^l)$, where $p=l=2$,
\cite{Ruan2003135} proved that a reduced 2-dimensional SIRS
model could exhibit backward bifurcation, Hopf bifurcation, and
even Bogdanov-Takens bifurcation and homoclinic bifurcation.
Although the choice of $p=l=2$ was not motivated by a specific
physical process,
this important result demonstrates that a nonlinear incidence rate can
induce backward bifurcation, and further generate complex 
dynamics in a simple disease model.
 
One of the focal points of our study will be a convex incidence function
which arose in a 4-dimensional HIV antioxidant therapy model~\cite{GW2009}.
In this model, the infectivity of infected cells was proposed to be
an increasing function of the density 
of reactive oxygen species, which themselves increase as the infection
progresses.  In~\cite{GW2009}, meaningful parameter values were carefully
chosen by data fitting to both experimental and clinical results.
In this parameter regime, the model was observed to capture the phenomenon
of viral blips, that is, long periods of undetectable viral load punctuated by
brief episodes of high viral load.
Viral blips have been observed clinically in HIV patients 
under highly active antiretroviral therapy~\cite{CCS1996,DZV1999,PWM2003,PMW2008},
and have received much attention in the research literature,
both by experimentalists~\cite{Fung2012,Garretta2012,Grennan15042012}
and mathematicians~\cite{Fraser2001June,JonesPerelson2005,CC2011,RP2009a,RP2009b}. Nonetheless,
the mechanisms underlying this phenomenon 
are still not thoroughly understood~\cite{Grennan15042012,RP2009b}.

We recently re-examined the model developed in~\cite{GW2009}, 
with the aim of providing new insight into the mechanism of
HIV viral blips~\cite{ZWY2013,ZWY2014}.
Focusing on the dynamics of the slow manifold of this model, 
we reduced the dimension of 
the 4-dimensional model by using quasi-steady state assumptions.
After a further generalization and parameter rescaling process,
a 2-dimensional in-host HIV model~\cite{ZWY2013,ZWY2014} was obtained, 
given by 
\begin{equation}
\dss \frac{{\mathrm d} X}{{\mathrm d} \tau} = 1-DX-(B+\frac{AY}{Y+C})XY, \qquad
\dss \frac{{\mathrm d} Y}{{\mathrm d} \tau} = (B+\frac{AY}{Y+C})XY-Y,
\label{Paper3_Eq4}
\end{equation}
where $X$ and $Y$ denote the concentrations
of the uninfected and infected cells respectively.
The constant influx rate and the death rate of $Y$ have been scaled to $1$.
The death rate of $X$ is $D$.
The 2-dimensional infection model above~(\ref{Paper3_Eq4}),
reduced from the 4-dimensional HIV model~\cite{GW2009},
preserves the viral blips observed in the HIV model. 

Importantly, system~(\ref{Paper3_Eq4}) is equivalent to the SIR model~(\ref{Paper3_Eq2}), except that the incidence function is convex, as we will 
show in section~\ref{HopfConvex}. 
This equivalence can be demonstrated
if we 
set $S=e_1 x$, 
$I=e_2 y$, and $t=e_3 \tau$ with $e_1=e_2=\frac{\mu N}{\delta +\mu}$
and $e_3=\frac{1}{\delta+\mu}$.  In this case, system~(\ref{Paper3_Eq2}) is rescaled to 
\begin{equation*}
\dss \frac{\mathrm d x}{\mathrm d \tau}
=1-\frac{\mu}{\delta+\mu}x -\frac{1}{\mu N} f(x,\,y), \qquad
\dss \frac{\mathrm d y}{\mathrm d \tau}
=\frac{1}{\mu N} f(x,\,y) - y,
\end{equation*}
which takes the same form as system~(\ref{Paper3_Eq4}).
Therefore, although system~(\ref{Paper3_Eq2}) arises in epidemiology and system~(\ref{Paper3_Eq4}) was derived as an in-host model, they are mathematically
equivalent in this sense.  We will refer to both systems~(\ref{Paper3_Eq2}) and (\ref{Paper3_Eq4}) as infection models.

In previous work~\cite{ZWY2013,ZWY2014}, we analyze the recurrent
behavior which emerges in system~(\ref{Paper3_Eq4}) in some detail.  
Recurrence is a particular form of oscillatory behavior
characterized by long periods of time
close to the uninfected equilibrium, punctuated by brief episodes
of high infection~\cite{YHW2006}.  Thus HIV viral blips are an example of 
recurrent behavior, but recurrence is a more general feature of many
diseases~\cite{YHW2006,ZWY2014}. We have demonstrated that
the increasing and saturating infectivity function of system~(\ref{Paper3_Eq4})
is critical to the emergence of recurrent behaviour.  
This form of an infectivity function corresponds to 
a convex incidence rate function in the associated 2-dimensional infection 
model~(\ref{Paper3_Eq4}), and can likewise induce recurrence in this model.
Convex incidence has been previously suggested to model
`cooperation effects' in epidemiology~\cite{Korobeinikov2005},
or cooperative phenomena in reactions between enzyme and
substrate, as proposed by~\cite{murray2002mathematical}.

The rest of this paper is organized as follows.
In Section 2, we study two 2-dimensional infection models, both 
closely related to system~(\ref{Paper3_Eq2}).  We show that 
system~(\ref{Paper3_Eq2}) 
with either (a) a concave incidence rate and saturating treatment term 
or (b) a convex incidence rate as shown in system~(\ref{Paper3_Eq4}), 
can exhibit backward bifurcation;
we then identify the necessary terms 
in the system equations which cause this phenomenon.
In Section 3, we demonstrate that in both models, 
backward bifurcation increases the likelihood of a
Hopf bifurcation on the upper branch of the infected equilibrium.
Studying system~(\ref{Paper3_Eq4}) in greater detail, 
we illustrate how the location of the Hopf bifurcations
and their directions (supercritical or subcritical),
determine the possible dynamical behaviors, concluding
that backward bifurcation facilitates Hopf bifurcation(s), which
then underly the rich behaviours observed in these models.
In Section 4, we explore backward bifurcation further,
presenting an autoimmune disease model which exhibits 
negative backward bifurcation, that is, a bifurcation
for which the turning point when $R_0 < 1$
is located in a region where one or more state variables is negative.
Although this bifurcation introduces two branches of the infected
equilibrium,
we demonstrate that, in the biologically feasible area, only forward bifurcation
exists in this model.
We then present a modification to this autoimmune model, motivated
by the recent discovery of a new cell type, which generates
a negative backward bifurcation and
Hopf bifurcation, and allows recurrent behavior to emerge.
A conclusion is drawn in Section 5.

\section{Backward bifurcation}

In this section, we study backward bifurcation in 
two 2-dimensional infection models. 
In particular, we explore the essential terms and parameter relations 
which are needed to generate backward bifurcation. 
Furthermore, we examine the convex incidence rate, and reveal
its underlying role in determining the emergence of backward bifurcation.

\subsection{Backward bifurcation in the infection model 
with concave incidence}

First, we consider the SIR model with concave incidence, 
described by the following equations~\cite{Zhou2012}: 
\begin{equation}
\dss \frac{{\mathrm d} S}{{\mathrm d} t} = \Lambda - \frac{\beta S I}{1 + k I} - d S, \quad 
\dss \frac{{\mathrm d} I}{{\mathrm d} t} =  \frac{\beta S I}{1 + k I} - (d+\gamma+\epsilon)I, \quad 
\dss \frac{{\mathrm d} R}{{\mathrm d} t} = \gamma I - d R,
\label{Paper3_Eq5}
\end{equation}
where $S$, $I$ and $R$ denote the number of susceptible, infective,
and recovered individuals, 
respectively;
$\Lambda$ is the constant recruitment rate of susceptibles;
$d$, $\gamma$, and $\epsilon$ represent the rates of natural death,
recovery, and the disease-induced mortality, respectively.
Note that the function $\frac{\beta S I}{1 + k I}$ is an incidence
rate of the form $\frac{kI^l S}{1+\alpha I^h}$~\cite{Liu1986}, when $l=h=1$. 
Here, $\beta$ is the infection rate, and $k$ measures 
the inhibition effect.
Since the variable $R$ is not involved in the first 
two equations, system~(\ref{Paper3_Eq5}) can be reduced to 
a 2-dimensional model as
\begin{equation}
\dss \frac{{\mathrm d} S}{{\mathrm d} t} = \Lambda - \frac{\beta S I}{1 + k I} - d S, \quad 
\dss \frac{{\mathrm d} I}{{\mathrm d} t} =  \frac{\beta S I}{1 + k I} - (d+\gamma+\epsilon)I.
\label{Paper3_Eq6}
\end{equation}
In~\cite{Zhou2012}, an additional assumption regarding limited medical treatment 
resources is introduced to the above model, leading to a model 
with a saturating treatment term, given by
\begin{equation}
\dss \frac{{\mathrm d} S}{{\mathrm d} t}=f_1(S,\,I)= \Lambda - \frac{\beta S I}{1 + k I} - d S, \quad 
\dss \frac{{\mathrm d} I}{{\mathrm d} t}=f_2(S,\,I)= \frac{\beta S I}{1 + k I} - (d+\gamma+\epsilon)I - 
 \frac{\alpha I}{\omega + I},
\label{Paper3_Eq7}
\end{equation}
where the real, positive parameter $\alpha$ represents the maximal 
medical resources per unit time, and the real, positive parameter $\omega$ 
is the half-saturation constant.
For simplicity, let the functions on the right-hand side of the equations in 
(\ref{Paper3_Eq7}) be $f_1$ and $f_2$, respectively. 
Then, the equilibrium solutions of system (\ref{Paper3_Eq7}) 
are obtained by solving the following algebraic  equations:
$f_1(S,\,I)=0$ and $f_2(S,\,I)=0$, from which the disease-free equilibrium 
can be easily obtained as
$\bar{\mathrm E}_0 = (\Lambda/d,\,0)$.  
For the infected equilibrium $\bar{\mathrm E}=(\bar{S},\,\bar{I})$, 
$\bar{S}$ is solved from $f_1=0$ as 
$\dss \bar{S}(I)=\frac{\Lambda(1+kI)}{(d k+\beta) I+d}$. 
Then, substituting $S=\bar{S}(I)$ into $f_2=0$ 
yields a quadratic equation of the form 
\begin{equation}
{\mathcal F}(I) = {\mathcal A}I^2 + {\mathcal B} I + {\mathcal C}  = 0,
\label{Paper3_Eq8}
\end{equation}
which in turn gives two roots: 
$\bar{I}_{1,\,2} = \frac{-{\mathcal B}\pm 
\sqrt{{\mathcal B}^2-4{\mathcal A}\,{\mathcal C}}}{2{\mathcal A}}$,
where,
${\mathcal A} = (d+\gamma+\epsilon)(dk+\beta)$, 
${\mathcal B} = [(dk +\beta)\omega+d](d+\gamma+\epsilon)
+(d k+\beta)\alpha-\beta\Lambda$,
${\mathcal C} = [(d+\gamma+\epsilon)\omega+\alpha]d - \beta\Lambda\omega$
for system~(\ref{Paper3_Eq7}).
Since all parameters take positive values, we have ${\mathcal A}>0$.
To get the two positive roots essential for backward bifurcation, 
it is required that 
${\mathcal B}<0$ and ${\mathcal C}>0$. 
Noticing that $\beta, \, \Lambda,\, \omega >0$, we can see that 
the infection force, $\beta$,
the constant influx of the susceptibles, $\Lambda$,
and the effect of medical treatment $\frac{\alpha I}{\omega + I}$ 
are indispensible terms for backward bifurcation. 
The number of positive infected equilibrium solutions changes from 
two to one when the value of $C$ passes from negative to positive, 
which gives a critical point at $C=0$, that is, 
$[(d+\gamma+\epsilon)\omega+\alpha]d = \beta\Lambda\omega$, 
which is equivalent to 
$R_0 = \frac{\beta\Lambda}{(d+\gamma+\epsilon+\alpha/\omega)d}=1$.

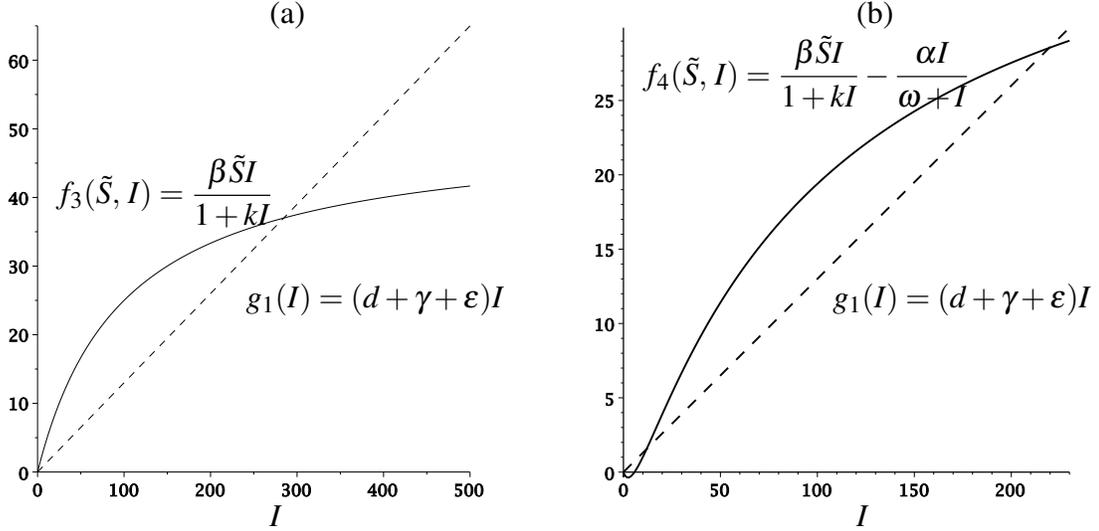
\begin{figure} 
\vspace{0.15in} 
\centering
\begin{overpic}[width=.40\textwidth]{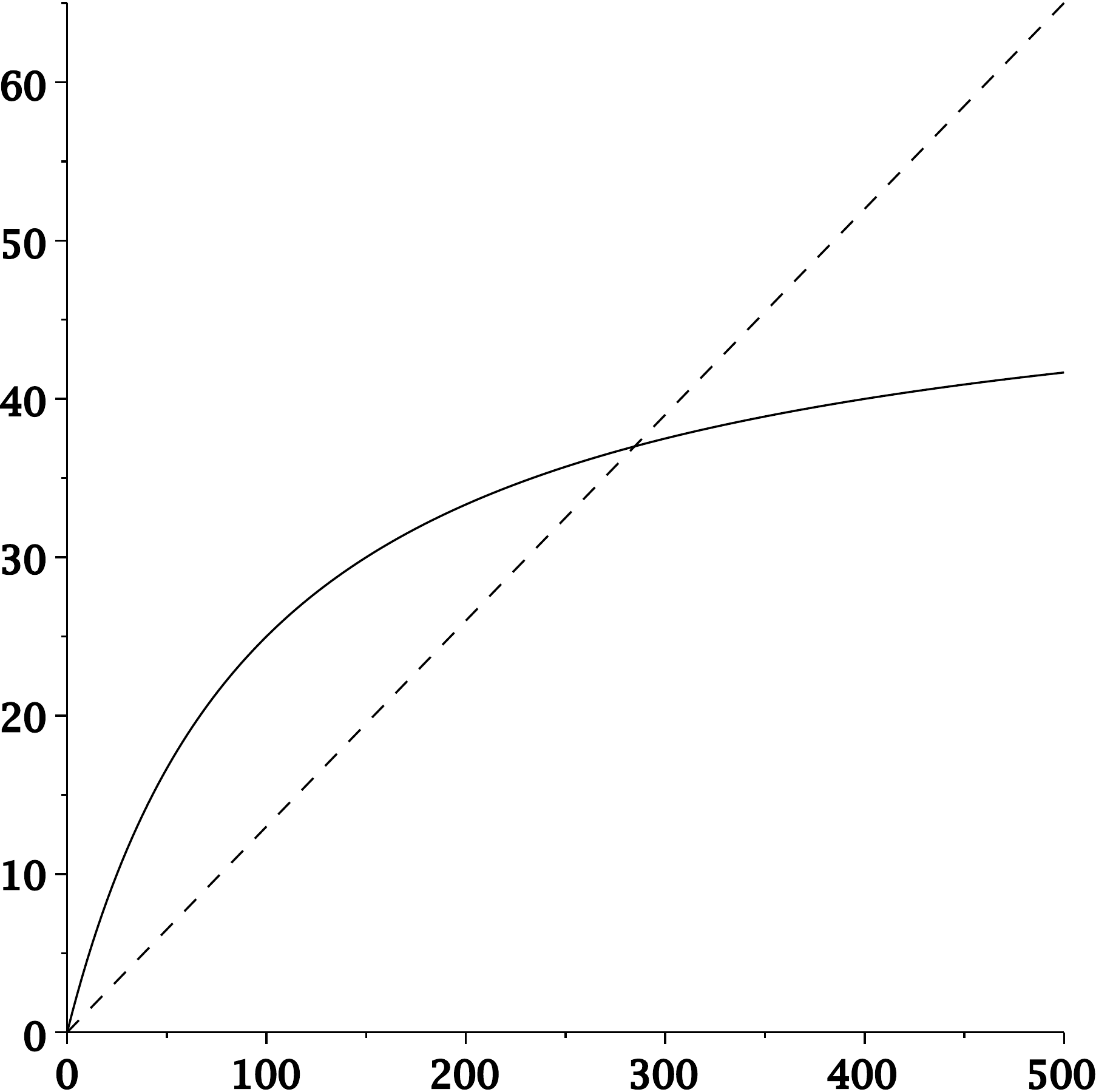}
\put(55,100){(a)}
\put(55,-6){$I$}
\put(10,62){$f_3(\tilde{S},\,I)
=\dss\frac{\beta \tilde{S} I}{1 + k I}$}
\put(50,40){$g_1(I)=(d+\gamma+\epsilon)I$}
\end{overpic}
\hspace{.5in}
\begin{overpic}[width=.40\textwidth]{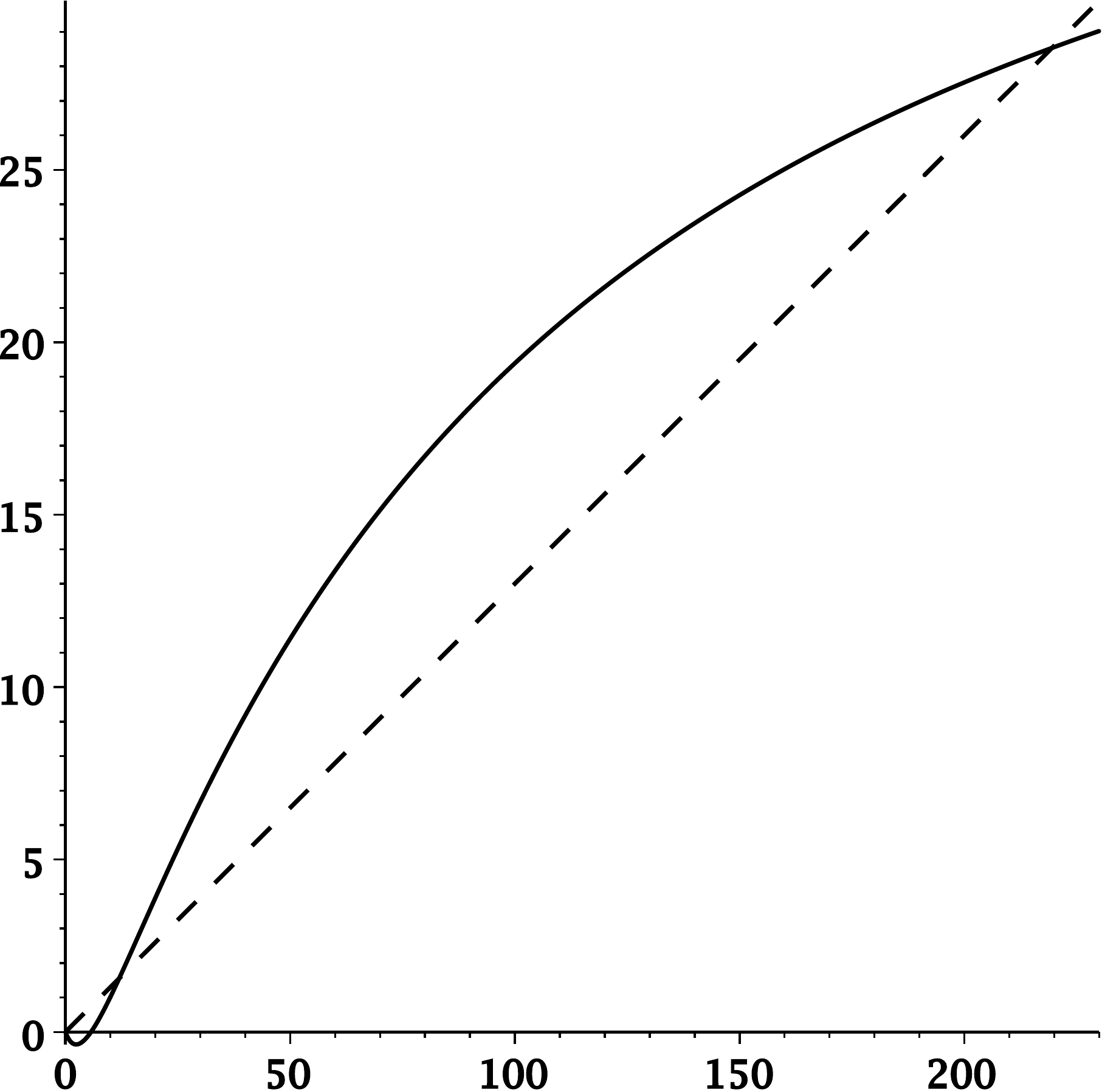}
\put(55,100){(b)}
\put(55,-6){$I$}
\put(10,87){$f_4(\tilde{S},\,I)=
\dss\frac{\beta \tilde{S} I}{1 + k I}-\frac{\alpha I}{\omega +I}$}
\put(50,40){$g_1(I)=(d+\gamma+\epsilon)I$}
\end{overpic}
\hspace{1.0in}
\caption{ 
Graphs of the incidence function $f_3$ in 
system~(\ref{Paper3_Eq5}), (\ref{Paper3_Eq6}) 
and function $f_4$ in system~(\ref{Paper3_Eq7})
with respect to $I$, for which $\tilde{S}=50$ has been used. 
The parameter values are chosen as $\beta=0.01$, $k=0.01$, $\alpha=6$, 
$\omega=7$, $d=0.1$, $\gamma=0.01$, $\epsilon=0.02$, 
according to~\cite{Zhou2012}.
The solid lines denote $f_3$ in (a) and $f_4$ in (b),
while the dashed ray lines in both graphs represent 
$g_1(I)=(d+\gamma+\varepsilon) I$. 
(a) the incidence function $f_3(S,I) = \frac{\beta SI}{1 +kI}$, showing 
one intersection point with $g_1$; and 
(b) the function $f_4(S,I) = \frac{\beta SI}{1 +kI}
-\frac{\alpha I}{\omega +I}$, showing two intersection points 
with $g_1$.} 
\label{Paper3_Fig1}
\end{figure}

On the other hand, we may infer the emergence of backward bifurcation
without solving the equilibrium conditions. 
If we do not consider the medical treatment term 
$\frac{\alpha I}{\omega + I}$ and remove it from system 
(\ref{Paper3_Eq7}), that leads to system~(\ref{Paper3_Eq6}), 
which is a typical example of an SIR model studied by 
(\ref{Paper3_Eq2}). By setting the incidence function as 
$f_3(S,\,I)=\frac{\beta S I}{1 + k I}$, 
we have $f_3(0,\,I)=f_3(S,\,0)=0$;
$\frac{\partial f_3(S,\,I)}{\partial S}= \frac{\beta I}{1+k I}>0$ and
$\frac{\partial f_3(S,\,I)}{\partial I}= \frac{\beta S}{(1+k I)^2}>0$ for all
$S,\,I>0$; and 
$\frac{\partial^2 f_3(S,\,I)}{\partial I^2}=-2\beta k S(1+k I)^{-3}<0$ for all
$S,\,I>0$.
Therefore, the incidence function $f_3(S,\,I)$, satisfies the
conditions given in~(\ref{Paper3_Eq3}). 
In particular, the function is concave, and can only
have one intersection point with the line $(d+\gamma+\epsilon)I$ in the 
$I$-$S$ plane, as shown in Figure~\ref{Paper3_Fig1}(a).
Thus, the uniqueness of the positive infected equilibrium implies that
backward bifurcation cannot occur in this case.
Moreover, according to the result in~\cite{Korobeinikov2005},
the uninfected and infected equilibria are globally asymptotically stable for 
$R_0=\frac{\beta\Lambda}{d(d+\gamma+\epsilon)}<1$ and $R_0>1$, respectively.
No complex dynamical behavior happens in system~(\ref{Paper3_Eq6}). 

In contrast, when we introduce the loss of the infectives 
due to medical treatment,
the dynamics of system (\ref{Paper3_Eq7}) differ greatly from
system (\ref{Paper3_Eq6}). In particular, 
backward bifurcation emerges and complex dynamical behaviors may occur.
To clarify this effect, we denote the function induced by $\dot{I}=0$ from~(\ref{Paper3_Eq7}) as
$f_4(S,\,I)=\frac{\beta S I}{1 + k I}-\frac{\alpha I}{\omega + I}$.
Note that $f_4(S,\,I)$ is not an incidence rate. 
But, if we fix $S=\tilde{S}>0$, there exist $0<I_1<I_2<+\infty$, such that
$\frac{\partial f_4(\tilde{S},\,I)}{\partial I}
=\frac{1}{(1+kI)^2(\omega+I)^2}
[\beta \tilde{S}(\omega+I)^2-\alpha \omega (1+kI)^2]>0$, 
$\forall \; I \in (0,\,I_2)$; and
$\frac{\partial^2 f_4(\tilde{S},\,I)}{\partial I^2} 
= -2k\beta \tilde{S} (1+kI)^{-3}
+ 2\alpha \omega (\omega+I)^{-3}>0$, 
$\forall \; I \in (0,\,I_1)$,
$\frac{\partial^2 f_4(\tilde{S},\,I)}{\partial I^2}=0$,
for $I=I_1$,
$\frac{\partial^2 f_4(\tilde{S},\,I)}{\partial I^2}<0$,
$\forall \; I \in (I_1,\,I_2)$.
Thus, $f_4(\tilde{S},\,I)$ actually has a convex-concave `$S$' shape, and may have
two positive intersection points with the ray line, 
$g_1(I)=(d+\gamma+\epsilon)I$, in the first quadrant; 
see Figure~\ref{Paper3_Fig1}(b).
These intersections contribute the two positive equilibrium solutions 
that are a necessary feature of backward bifurcation.

In summary we may conclude that the necessary terms which should be 
contained in system (\ref{Paper3_Eq7}) in order to have
backward bifurcation are the constant influx $\Lambda$,
the infection force $\beta$,
and the saturating medical treatment $\frac{\alpha I}{\omega +I}$.

\subsection{Backward bifurcation in the infection model 
with convex incidence}\label{HopfConvex}

Now we consider the 2-dimensional infection model~(\ref{Paper3_Eq4})
which exhibits viral blips, studied in~\cite{ZWY2013,ZWY2014}.
The motivation for this model was a series of clinical discoveries
indicating that viral infection can increase 
the density of a harmful chemical substance
\cite{gil2003contribution,LK1999,Schwarz1996,IG1997}, 
thereby amplifying an associated biochemical reaction~\cite{SMDW2005}, 
and thus accelerating the infection rate~\cite{gil2003contribution}.
This cooperative phenomenon in viral infection is expressed by
an increasing, saturating infectivity function:
$\tss(B+\frac{AY}{Y+C})$. 
According to the principle of mass action, the incidence function
is then denoted as $\tss (B+\frac{AY}{Y+C})XY$, which is a convex
function with respect to the infectives' density $Y$.

To analyze the occurrence of possible backward bifurcation, we first examine 
the two equilibrium solutions from the following equations: 
\begin{equation}
f_5(X,\,Y) = 1-DX-(B+\frac{AY}{Y+C})XY= 0, \quad
f_6(X,\,Y) = (B+\frac{AY}{Y+C})XY-Y=0,
\label{Paper3_Eq9}
\end{equation}
where all parameters $A$, $B$, $C$ and $D$ are positive constants. 
It is easy to find the uninfected equilibrium 
$\bar{\mathrm E}_0=(\bar{X}_0,\,\bar{Y}_0)=(\frac{1}{D},\,0)$, 
whose characteristic polynomial has two roots: $\lambda_1=-D<0$,
and $\lambda_2 = \frac{B}{D}-1$, which gives $R_0=\frac{B}{D}$.
Consequently, $\bar{\mathrm E}_0$ is stable (unstable) for $R_0<1\,(>1)$.
To find the infected equilibrium solution, setting $f_6(X,\,Y) =0$ 
yields $\bar{X}_1(Y)=\frac{Y+C}{(A+B)Y+BC}$, which is then substituted 
into $f_5(X,\,Y)=0$ to give the following quadratic equation: 
\begin{equation}
{\mathcal F}_5(Y)=(A+B)Y^2+(BC+D-A-B)Y+C(D-B)=0.
\label{Paper3_Eq10}
\end{equation}
In order to have two real, positive roots, 
two conditions must be satisfied, that is, 
$BC+D-A-B<0$ and $D-B>0$, or in compact form, $0<D-B<A-BC$. 
The condition $D-B>0$ is equivalent to $0<R_0=\frac{B}{D}<1$, 
which is a necessary condition for backward bifurcation.
Moreover, the positive influx constant, having been scaled to $1$, 
is a necessary term for the positive equilibrium of $Y$.
Therefore, the positive influx rate term and 
the increasing and saturating infectivity function 
are necessary for backward bifurcation.

\begin{figure}
\vspace{0.15in} 
\centering
\begin{overpic}[width=.40\textwidth]{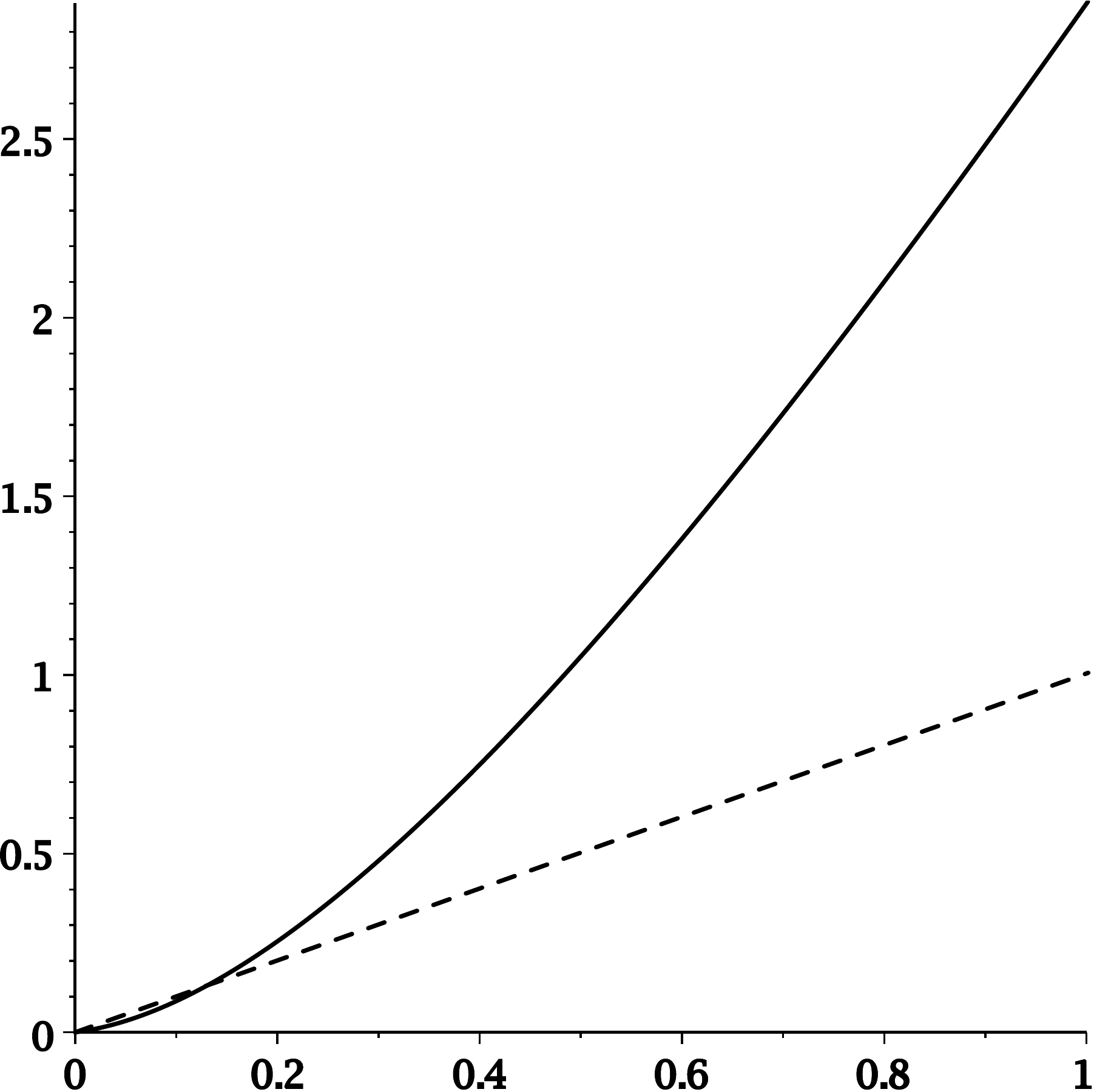}
\put(50,-5){$Y$}
\put(50,100){(a)}
\put(70,55){$f_7(\tilde{X},\,Y)$ in (\ref{Paper3_Eq11})}
\put(70,20){$g_2(Y)=Y$}
\put(10,55){\includegraphics[width=.15\textwidth]{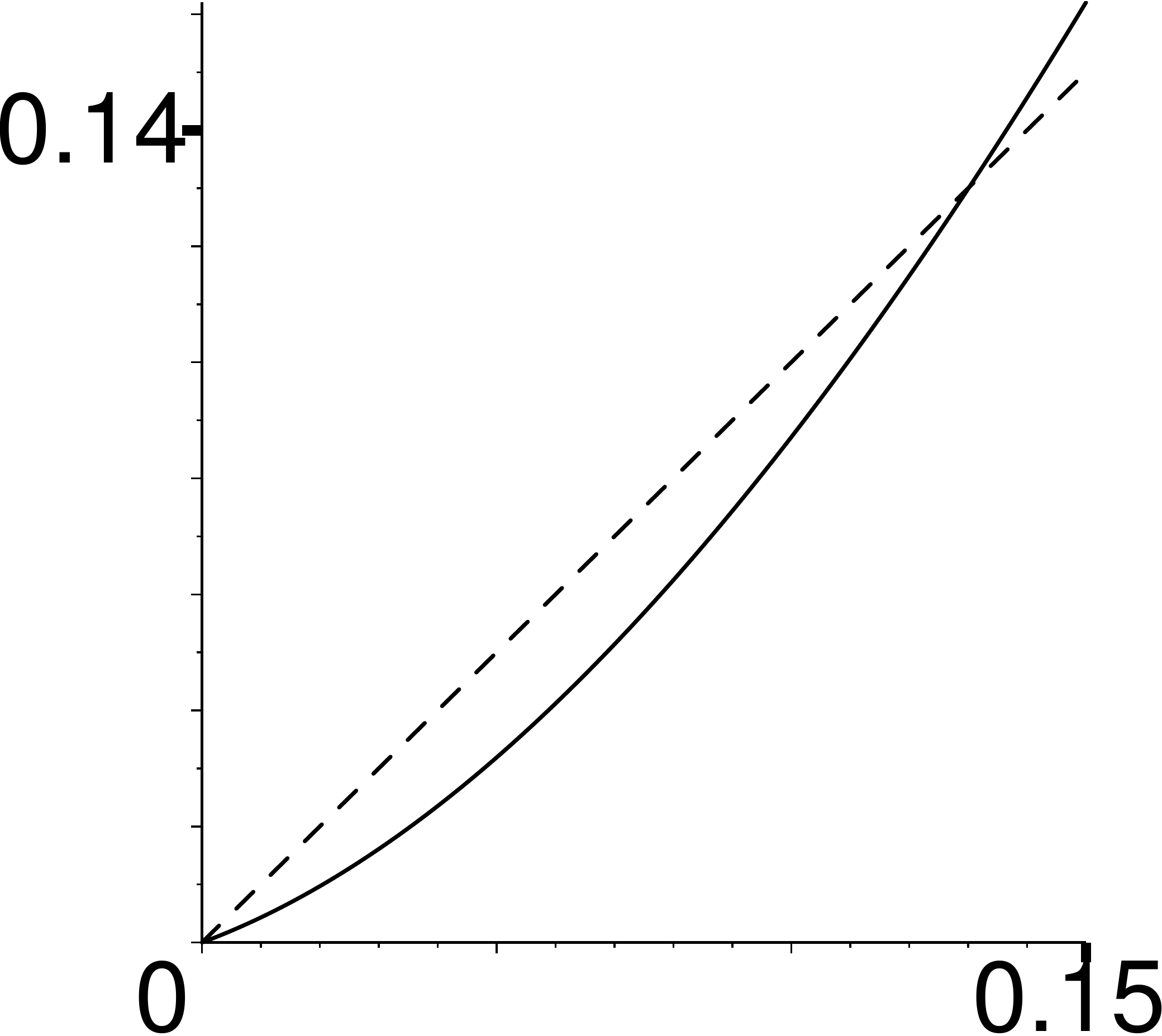}}
\put(27,53){\tiny $Y$}
\put(7,20){\line(1,0){23}}
\put(30,20){\line(0,-1){15}}
\end{overpic}
\hspace{0.5in}
\begin{overpic}[width=.40\textwidth]{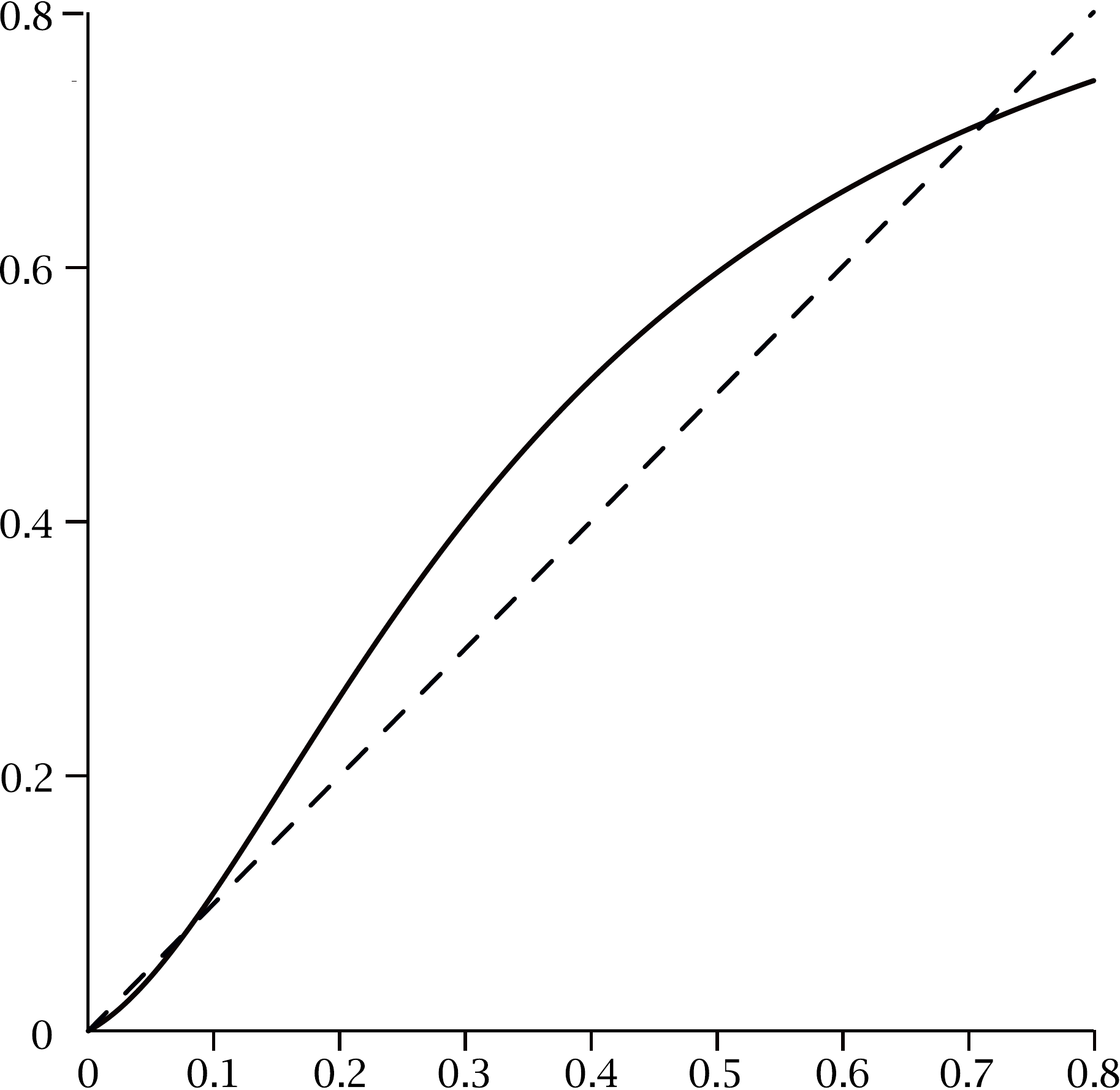}
\put(50,-5){$Y$}
\put(50,100){(b)}
\put(15,65){$f_7(Y)$ in (\ref{Paper3_Eq12})}
\put(70,95){$g_2(Y)=Y$}
\put(60,10){\includegraphics[width=.15\textwidth]{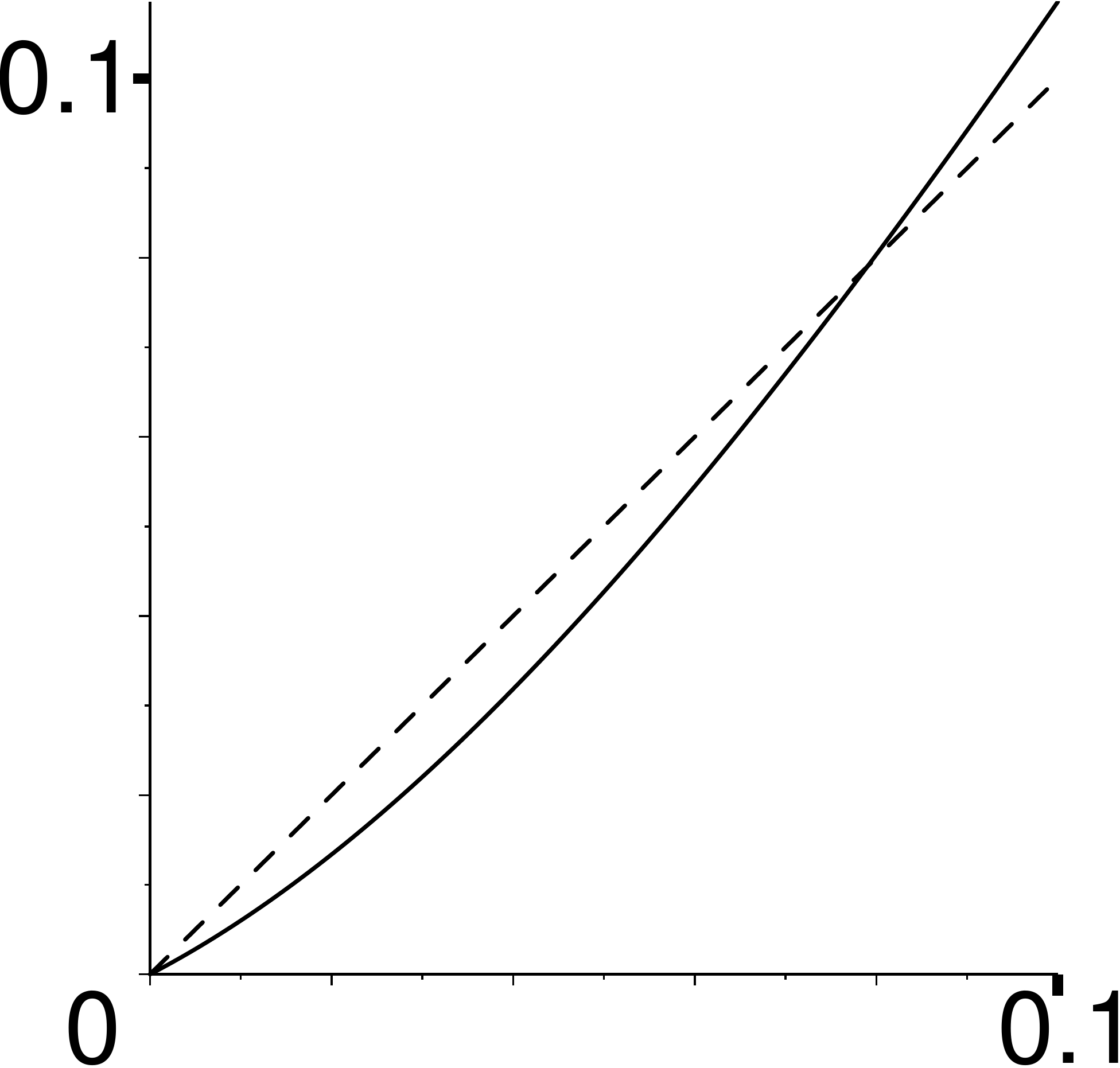}}
\put(78,8){\tiny $Y$}
\put(7,20){\line(1,0){18}}
\put(25,20){\line(0,-1){15}}
\end{overpic}
\vspace{0.1in}
\caption{Graphs of the incidence functions 
$f_7(\tilde{X},Y)$ and $f_7(Y)$
for the parameter values $A=0.364$,
$B=0.03$, $C=0.823$, and $D=0.057$. 
The incidence functions are denoted by the solid lines, while 
the ray lines, determined by $g_2(Y)=Y$, are denoted by dotted lines:
(a) the incidence function $f_7(\tilde{X},Y)$, showing 
one intersection point with $g_2$ with an inset, 
with a fixed value $\tilde{X}=12.54$; and 
(b) the incidence function $f_7(Y)$, showing two intersection points 
with an inset.} 
\label{Paper3_Fig2}
\end{figure}
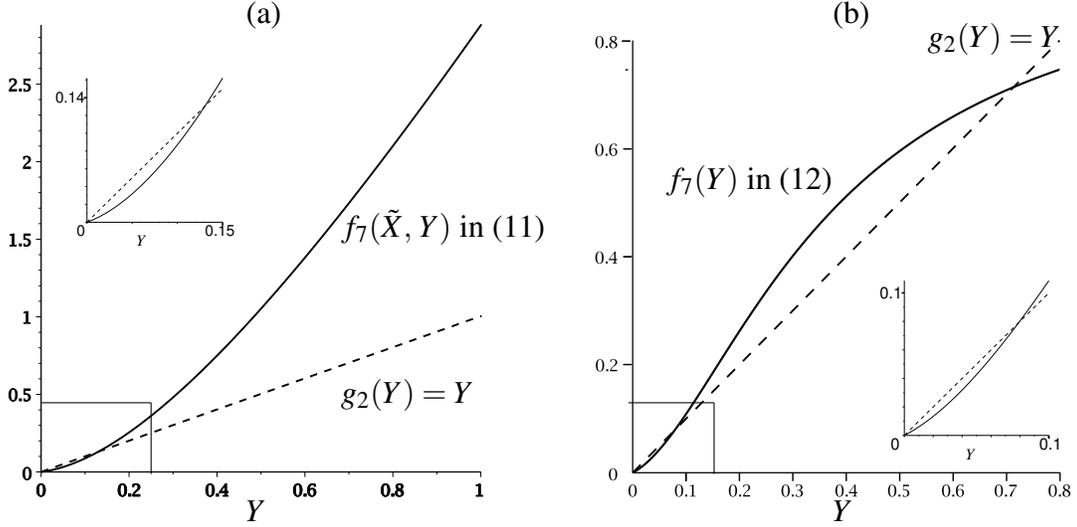

In the rest of the subsection, we further examine 
the incidence function, 
\begin{equation}
f_7(X,\,Y)=(B+\frac{AY}{Y+C})XY,
\label{Paper3_Eq11}
\end{equation}
without solving the equilibrium solutions.
The incidence function $f_7$ obviously satisfies 
the condition (\ref{Paper3_Eq3a}), as well as 
the condition (\ref{Paper3_Eq3b}) since
$\frac{\partial}{\partial X}f_7(X,\,Y)=[B+AY(Y+C)^{-1}]Y >0$ and
$\frac{\partial}{\partial Y}f_7(X,\,Y)= ACXY(Y+C)^{-2}+[B+AY(Y+C)^{-1}]X>0$
for all $X,\,Y>0$. However, the second partial derivative of 
$f_7(X,\,Y)$ with respect to $Y$, 
$\frac{\partial^2}{\partial Y^2}f_7(X,\,Y)= 2AC^2X(X+C)^{-3}>0$ for all
$X,\,Y>0$, showing that $f_7(X,\,Y)$ is a convex function with respect to 
the variable $Y$.
Consequently, $f_7(X,\,Y)$ can only have one intersection with
$g_2(Y)=Y$, implying that only one equilibrium solution would exist 
if we only consider the second equation in (\ref{Paper3_Eq9}),
as shown Figure~\ref{Paper3_Fig2} (a).
However, when considering both conditions given in (\ref{Paper3_Eq9}) for 
equilibrium solutions, we 
will have two intersection points between $f_7$ and $g_2$.
According to the first equation in (\ref{Paper3_Eq9}), that is $f_5(X,\,Y)=0$, 
we can use $Y$ to express $X$ in the equilibrium 
state as $\bar{X}(Y)=(Y+C)[(A+B)Y^2+(BC+D)Y+DC]^{-1}$.
Substituting $\bar{X}(Y)$ into $f_7(X,\,Y)$ in~(\ref{Paper3_Eq11}), we obtain
\begin{equation}
f_7(Y)=Y[(A+B)Y+BC][(A+B)Y^2+(BC+D)Y+CD]^{-1},
\label{Paper3_Eq12}
\end{equation}
and $\frac{\partial}{\partial Y}f_7(Y)=
D[(A+B)Y^2+2(A+B)CY+BC^2][(A+B)Y^2+(BC+D)Y+CD]^{-2}>0$ for all
$X,\,Y>0$.
However, the sign of $\frac{\partial^2}{\partial Y^2}f_7(Y)= 
-2D[(A+B)^2Y^3+3C(A+B)^2Y^2+3(A+B)BC^2Y+(B^2C-AD)C^2]
[(A+B)Y^2+(BC+D)Y+CD]^{-3}$, 
could alter at the inflection point
from positive to negative as $Y$ increases.
Therefore, with appropriate parameter values, $f_7(Y)$ can have 
a convex-concave `$S$' shape, yielding two intersection points with 
the ray line, $g_2(y)$, in the first quadrant of the $X$-$Y$ plane,
as shown in Figure~\ref{Paper3_Fig2} (b). 
The above discussion, as illustrated in Figure~\ref{Paper3_Fig2},
implies that system (\ref{Paper3_Eq4}) 
can have two positive equilibrium solutions when $R_0 < 1$,
and thus backward bifurcation may occur.

\begin{remark}
Summarizing the discussions and results given in this section 
indicates that a disease model with a convex-concave incidence function
may lead to backward bifurcation, which in turn implies: 
(a) the system has at least 
two equilibrium solutions, and the two equilibrium solutions intersect at 
a transcritical bifurcation point; and (b) at least one of the equilibrium 
solutions is determined by a nonlinear equation.
\label{Paper3_Remark1}
\end{remark}

\section{Hopf bifurcation}

In the previous section, we studied backward bifurcation 
and established the necessary conditions for the occurrence of 
backward bifurcation in two models. In this section, we turn to 
Hopf bifurcation, since it typically underlies the change of
stability in the 
upper branch of the infected equilibrium,
the key condition in determining whether a model can exhibit 
oscillation or even recurrence. 
Again, we will present detailed studies for the two models.  

\subsection{Hopf bifurcation in the infection model 
with concave incidence}

In this subsection, we study two cases of an infection model with
concave incidence: system~(\ref{Paper3_Eq6}) and (\ref{Paper3_Eq7}).
First, we discuss the equilibrium solutions and 
their stability by using the Jacobian matrix, denoted by $J$,
and examining the corresponding characteristic polynomial,
\begin{equation}
P|_{J}(L) = L^2 + \mathrm{Tr}(J) L + \mathrm{Det}(J).
\label{Paper3_Eq13}
\end{equation}
Bifurcation analysis is conducted by choosing $\Lambda$ as the 
bifurcation parameter.

First, we consider the case without saturating medical treatment, 
system (\ref{Paper3_Eq6}). This system satisfies 
the three conditions in (\ref{Paper3_Eq3}), and consequently,
its uninfected equilibrium 
$\bar{\mathrm{E}}_0= (\frac{\Lambda}{d},\,0)$ 
is globally asymptotically stable if
$R_0=\frac{\beta\Lambda}{(d+\gamma+\epsilon)d}\le 1$, while the infected
equilibrium $\bar{\mathrm E}_1 = (\frac{k\Lambda+d+\gamma+\epsilon}{dk+\beta},
\,\frac{\beta\Lambda-(d+\gamma+\epsilon)d}{(dk+\beta)(d+\gamma+\epsilon)})$
emerges and is globally asymptotically stable if $R_0 > 1$.
Therefore, for this case the system has only one 
transcritical bifurcation point at $R_0=1$ and no complex dynamics 
can occur.

Next, with the saturating treatment term, system (\ref{Paper3_Eq7}) violates 
the conditions established for model (\ref{Paper3_Eq3}), 
but leads to the possibility of complex dynamical behaviors. 
In fact, evaluating the Jacobian matrix $J_1=J|_{(\ref{Paper3_Eq7})}
(\bar{\mathrm E}_0)$ at the uninfected equilibrium, 
$\bar{\mathrm E}_0= (\frac{\Lambda}{d},\,0)$, yields the characteristic 
polynomial in the form of (\ref{Paper3_Eq13}), denoted by $P|_{J_1}(L)$,
with $\mathrm{Tr}(J_1)= \left(-\frac{\beta\Lambda}{d}+\epsilon \right.$ 
$\left. +\frac{\alpha}{\omega}+2d \right)$, and $\mathrm{Det}(J_1) = 
\left(-\beta \Lambda+d^2+d\epsilon+\frac{\alpha d}{\omega}\right) = 
\mathrm{Tr}(J_1) \, d - d^2$. This indicates that \\$\mathrm{Det}(J_1)<0$ when
$\mathrm{Tr}(J_1)=0$, and thus Hopf bifurcation cannot occur from 
$\bar{\mathrm E}_0$.
On the other hand, a static bifurcation can occur when 
$\mathrm{Det}(J_1) =0 $, that is, 
$\Lambda_S = \frac{1}{\beta}(d^2+d \epsilon + \frac{\alpha d}{\omega})$,
where the subscript `$S$' refers to \emph{static bifurcation}. 
Therefore, $\bar{\mathrm E}_0$ is 
stable (unstable) for $\Lambda< \Lambda_S \ (> \Lambda_S)$, 
or $R_0< 1 (>1)$, with 
$R_0=\beta\Lambda d^{-1}(d+\gamma+\epsilon+\frac{\alpha}{\omega})^{-1}$
\cite{Zhou2012}.

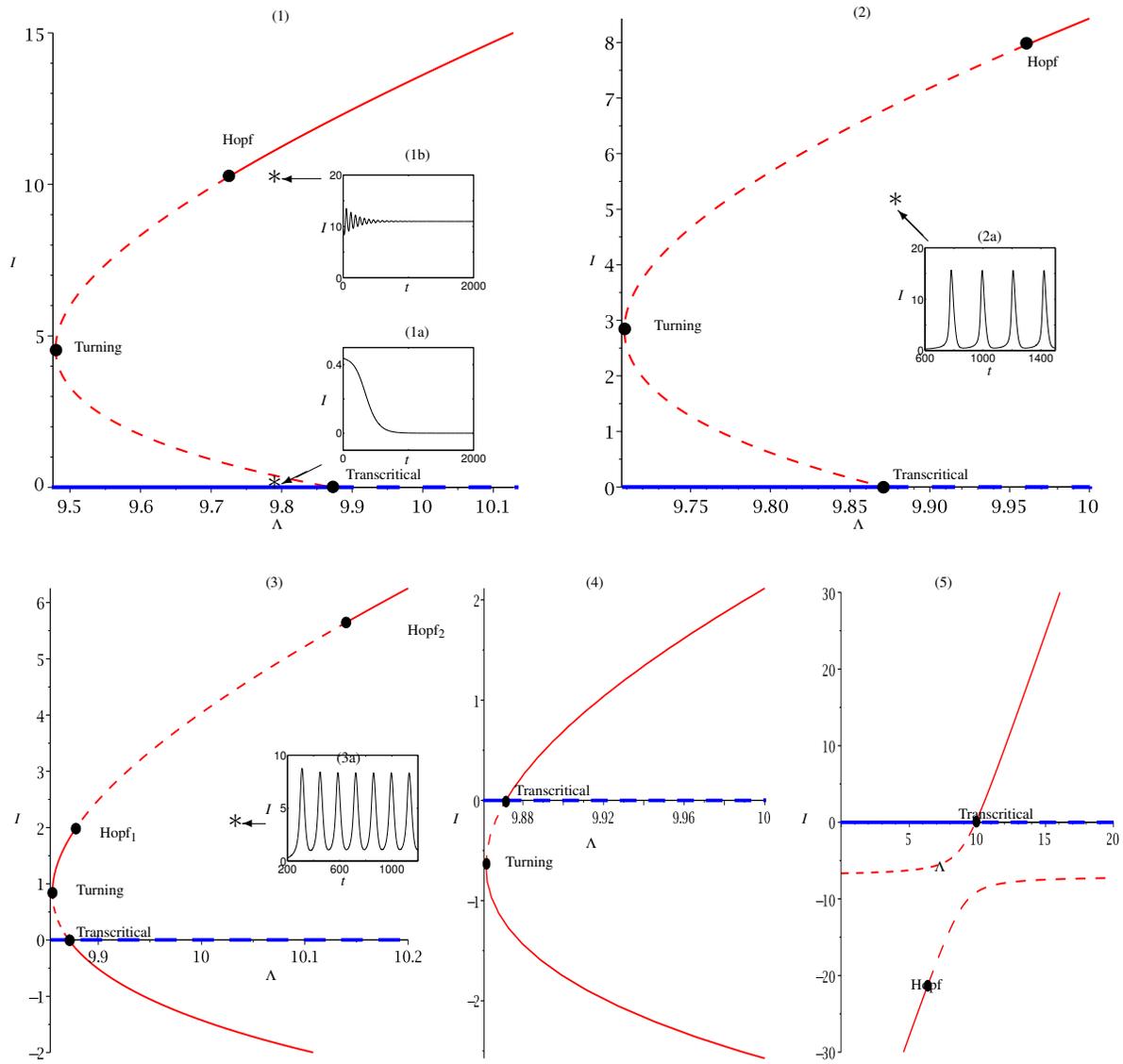
\begin{figure}
\vspace{0.15in} 
\centering
\hspace{-0.5cm}
\begin{overpic}[width=.44\textwidth]
{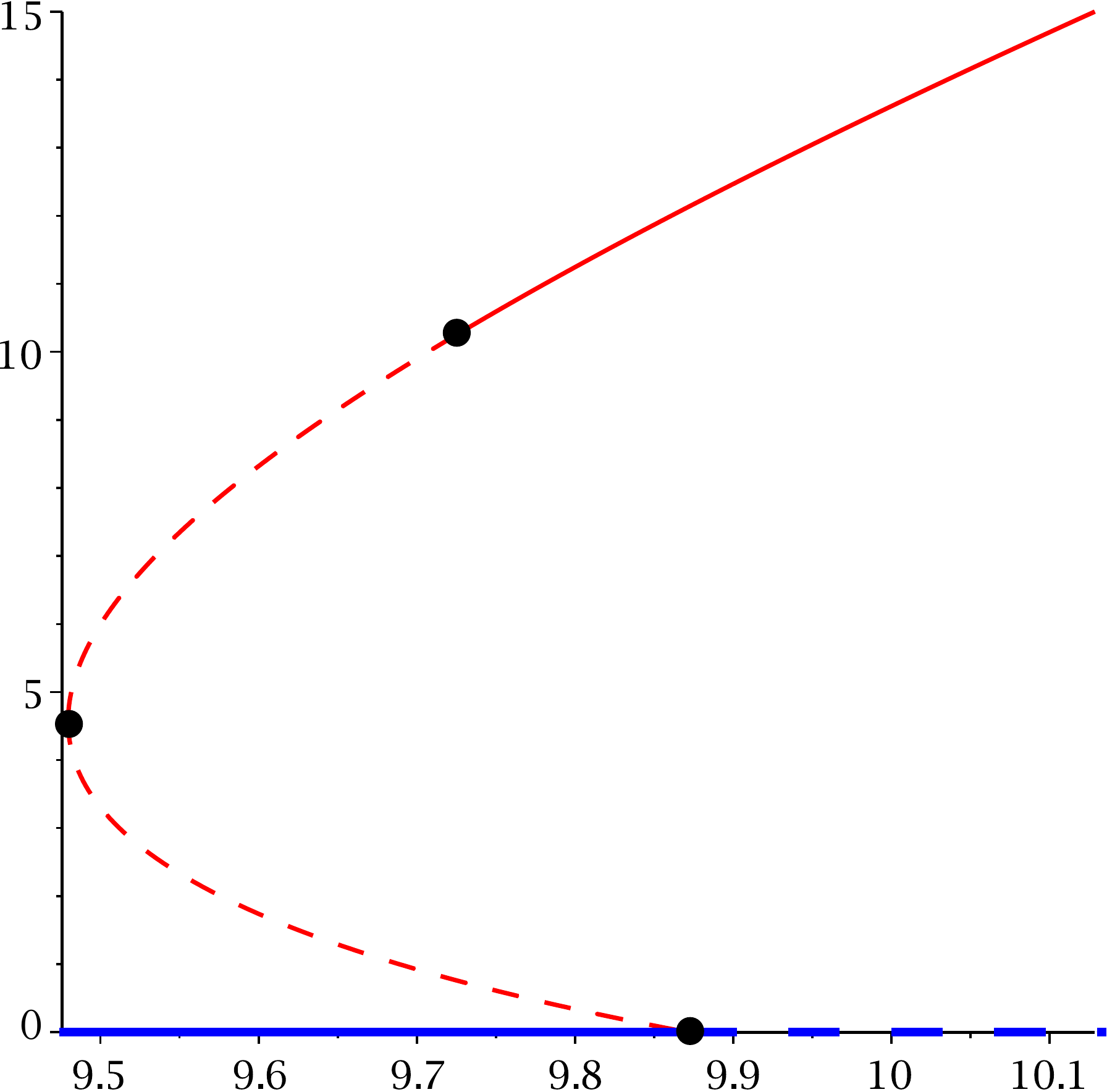}
\put(50,100){\tiny (1)}
\put(-3,50){\tiny $I$}
\put(50,-3){\tiny$\Lambda$}
\put(65,7){\tiny Transcritical}
\put(10,33){\tiny Turning}
\put(40,75){\tiny Hopf}
\put(49,5){$\ast$}
\put(49,67){$\ast$}
\put(60,45){\includegraphics[width=.15\textwidth]{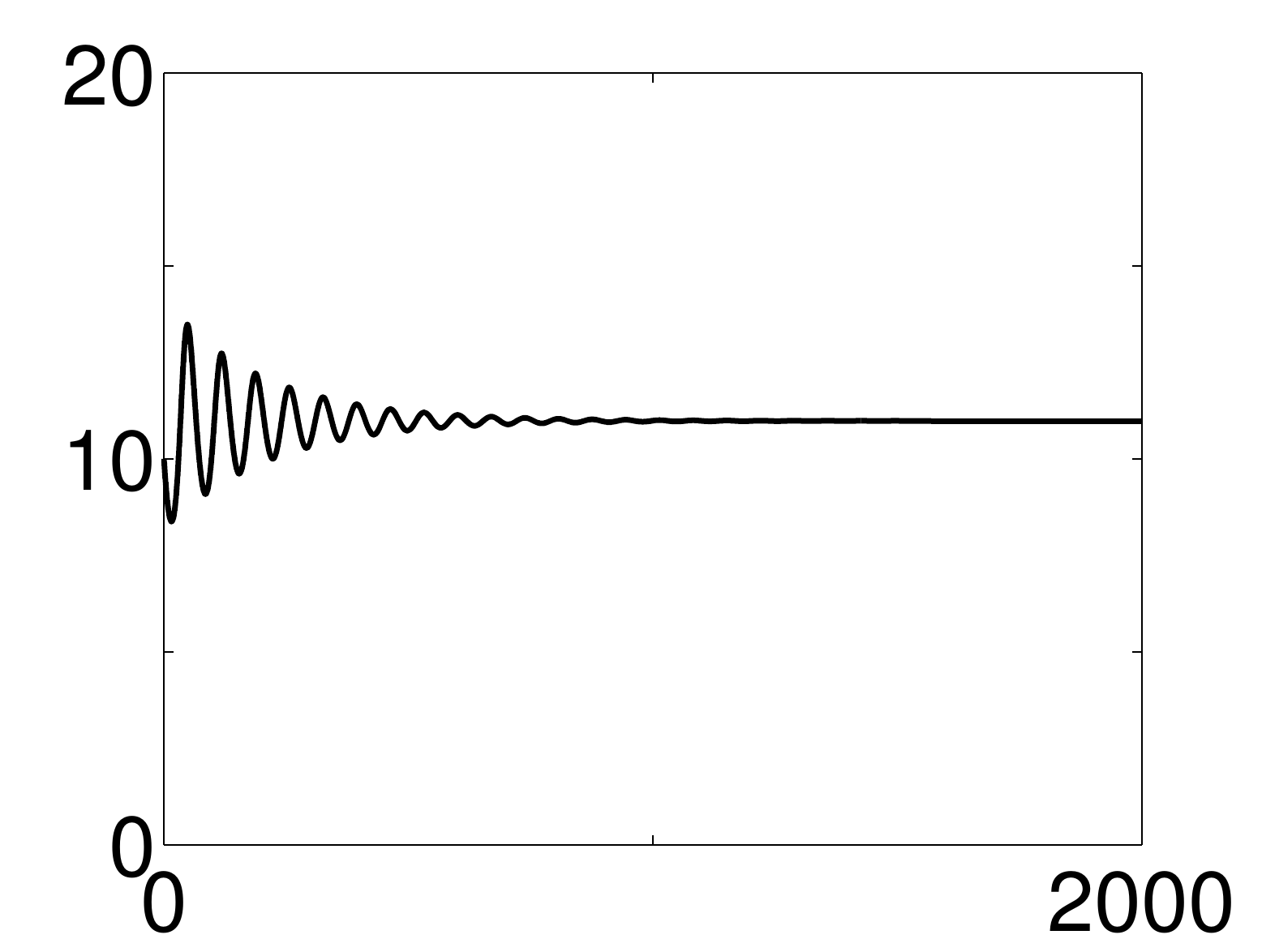}}
\put(60,68){\vector(-1,0){8}}
\put(77,72){\tiny (1b)}
\put(60,57){\tiny $I$}
\put(77,45){\tiny $t$}
\put(60,10){\includegraphics[width=.15\textwidth]{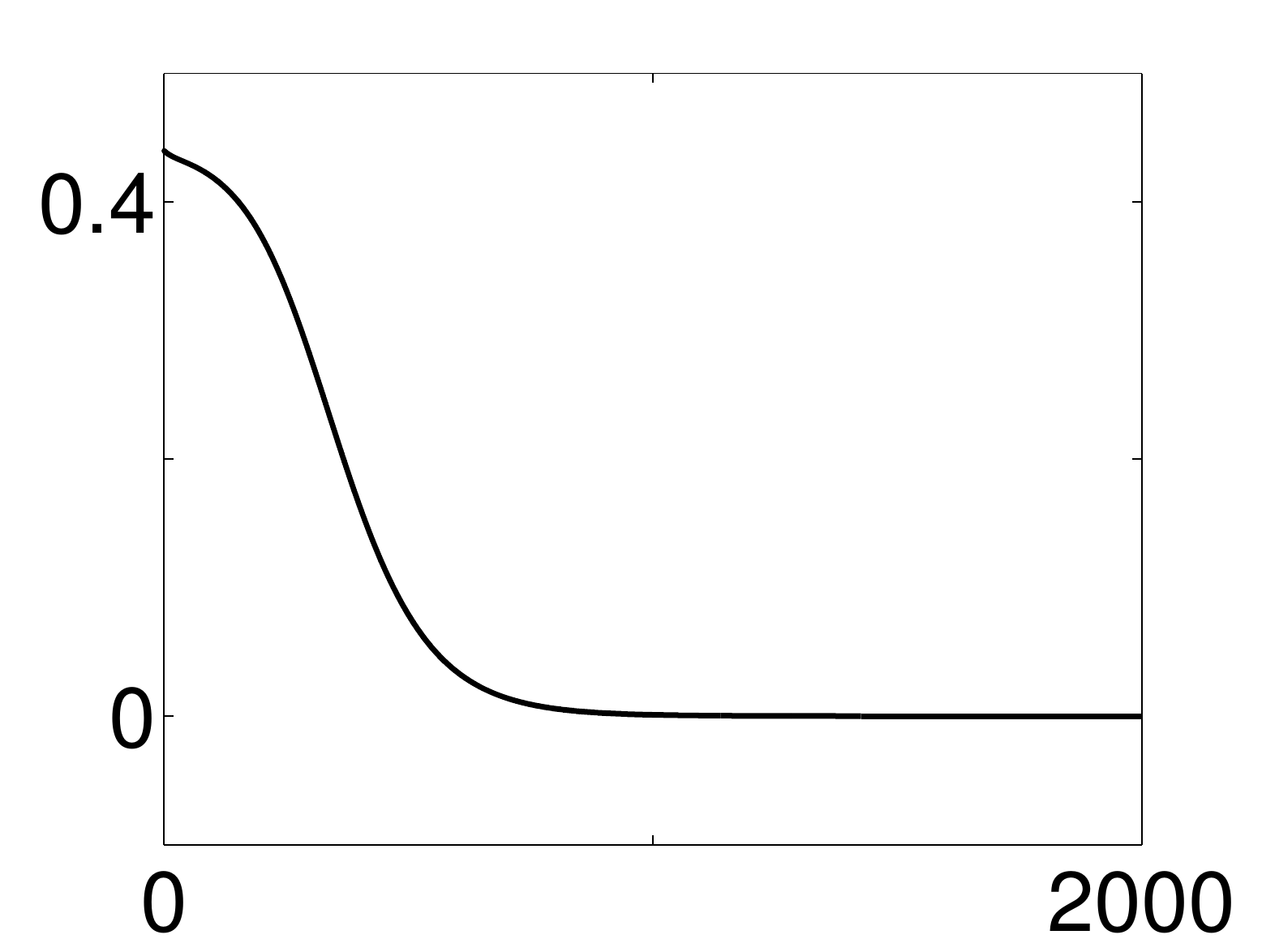}}
\put(60,10){\vector(-2,-1){8}}
\put(77,36){\tiny (1a)}
\put(60,22){\tiny $I$}
\put(77,10){\tiny $t$}
\end{overpic}
\hspace{1.0cm}
\begin{overpic}[width=.44\textwidth]
{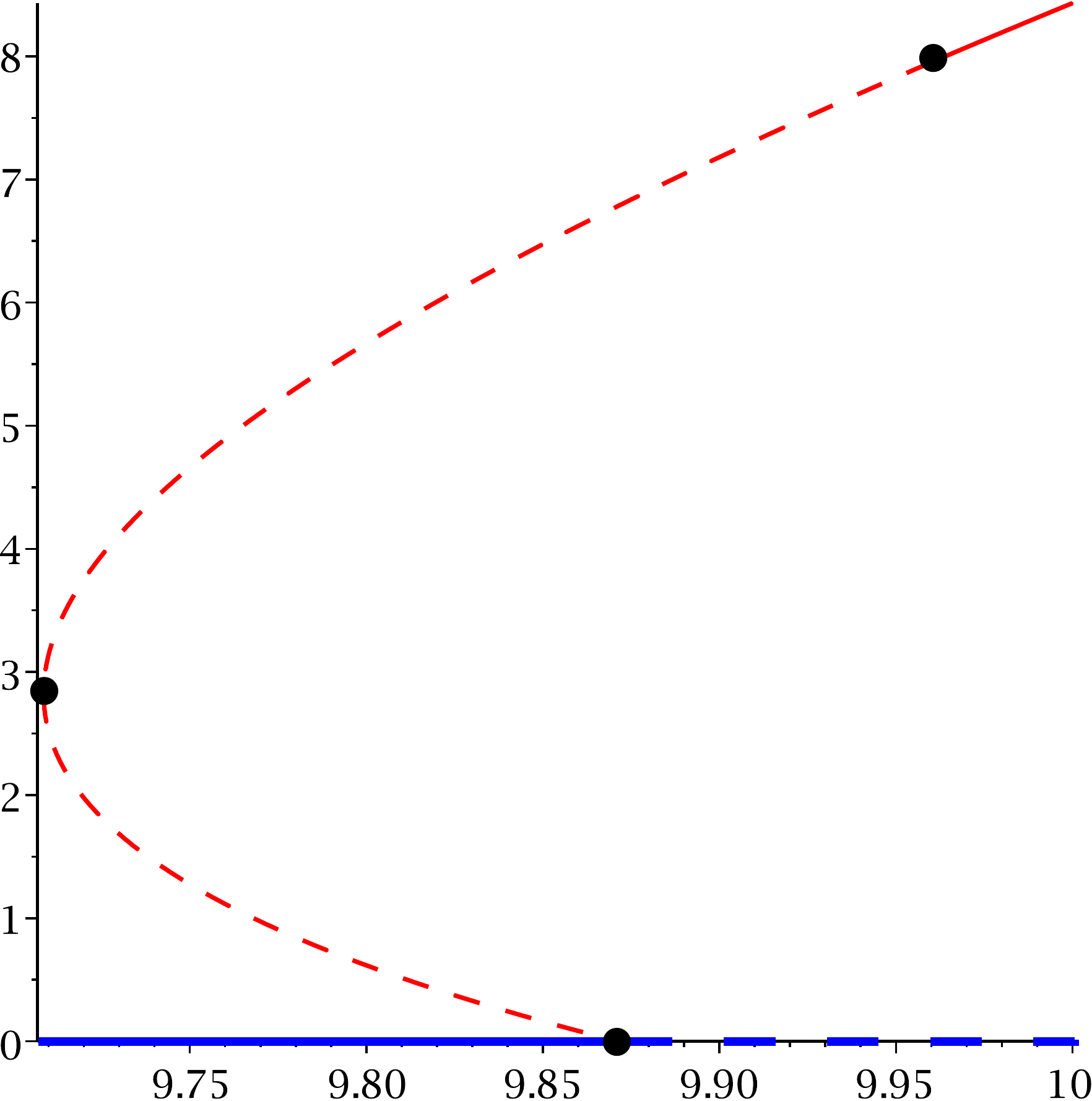}
\put(50,100){\tiny (2)}
\put(-3,50){\tiny $I$}
\put(50,-3){\tiny$\Lambda$}
\put(58,7){\tiny Transcritical}
\put(10,37){\tiny Turning}
\put(85,90){\tiny Hopf}
\put(57,62){$\ast$}
\put(60,30){\includegraphics[width=.15\textwidth]{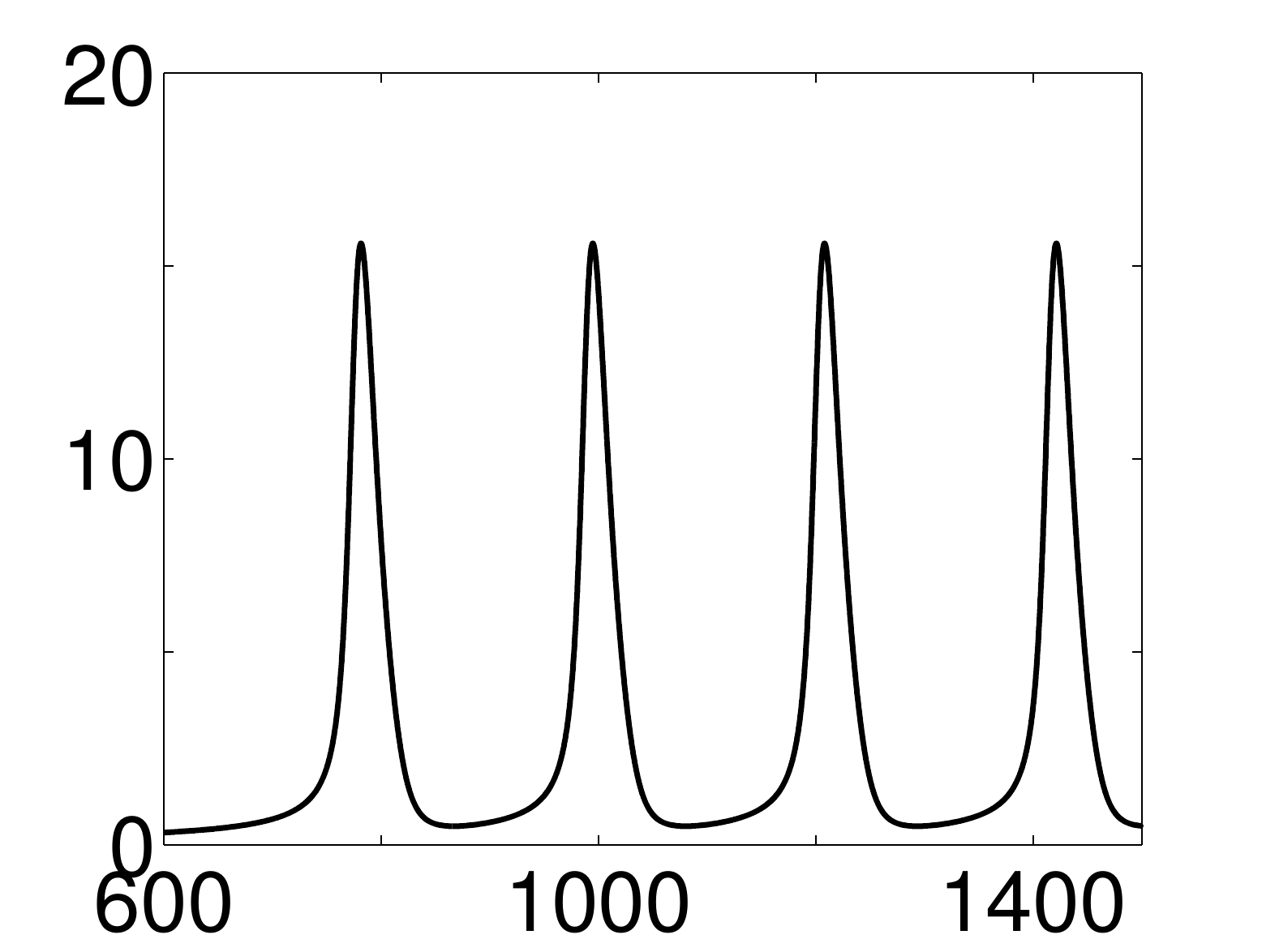}}
\put(65,55){\vector(-1,1){6}}
\put(75,55){\tiny (2a)}
\put(59,43){\tiny $I$}
\put(77,28){\tiny $t$}
\end{overpic}
\\ \vspace{1cm}
\hspace{-0.2cm}
\begin{overpic}[width=.35\textwidth,height=.29\textheight]{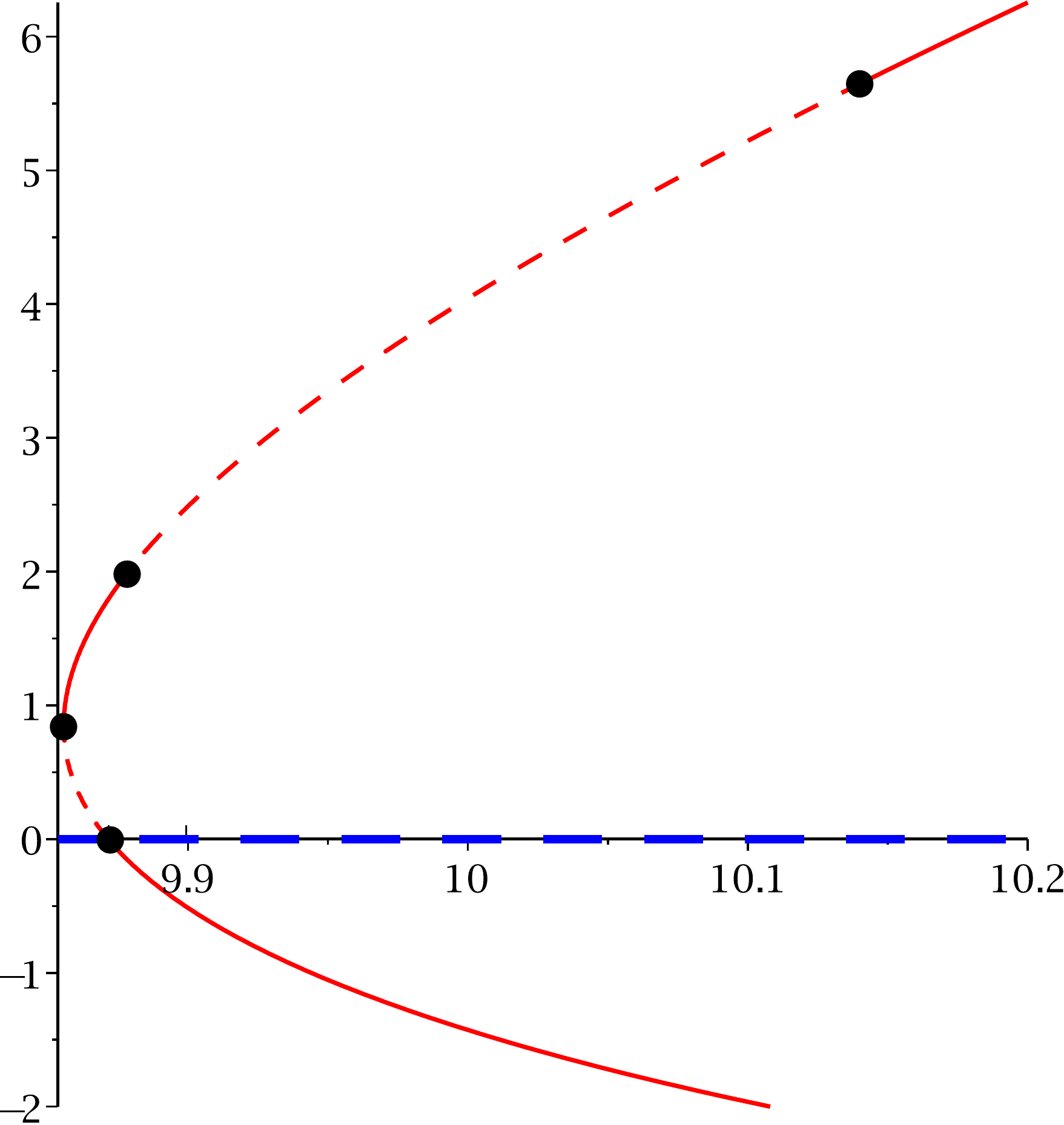}
\put(50,100){\tiny (3)}
\put(-3,50){\tiny $I$}
\put(50,17){\tiny$\Lambda$}
\put(10,26){\tiny Transcritical}
\put(10,35){\tiny Turning}
\put(80,90){\tiny Hopf$_2$}
\put(15,47){\tiny Hopf$_1$}
\put(42,49){$\ast$}
\put(50,39){\includegraphics[width=.15\textwidth,height=.08\textheight]{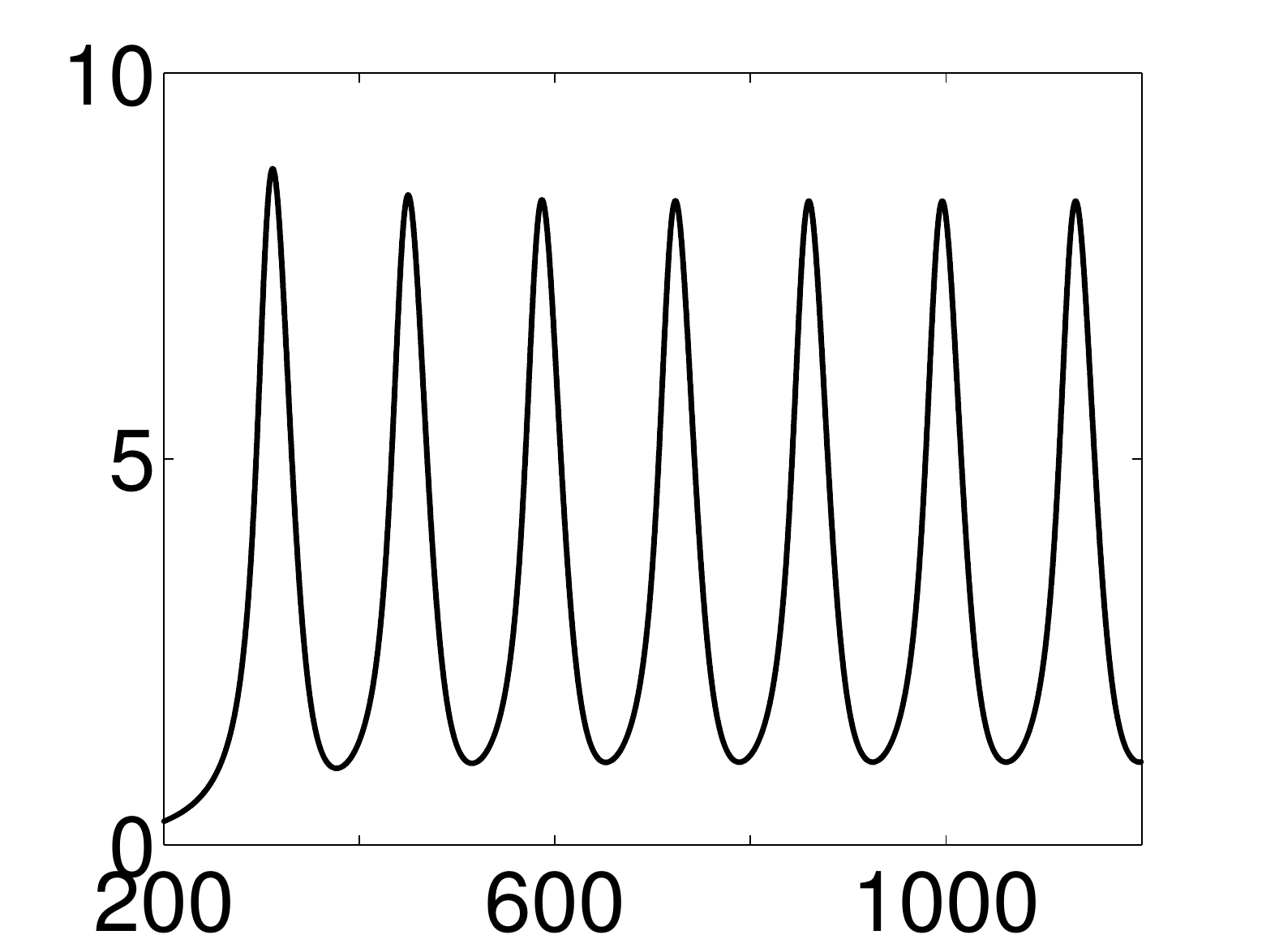}}
\put(50,50){\vector(-1,0){5}}
\put(65,63){\tiny (3a)}
\put(50,52){\tiny $I$}
\put(66,37){\tiny $t$}
\end{overpic}
\hspace{.5cm}
\begin{overpic}[width=.27\textwidth,height=.29\textheight]{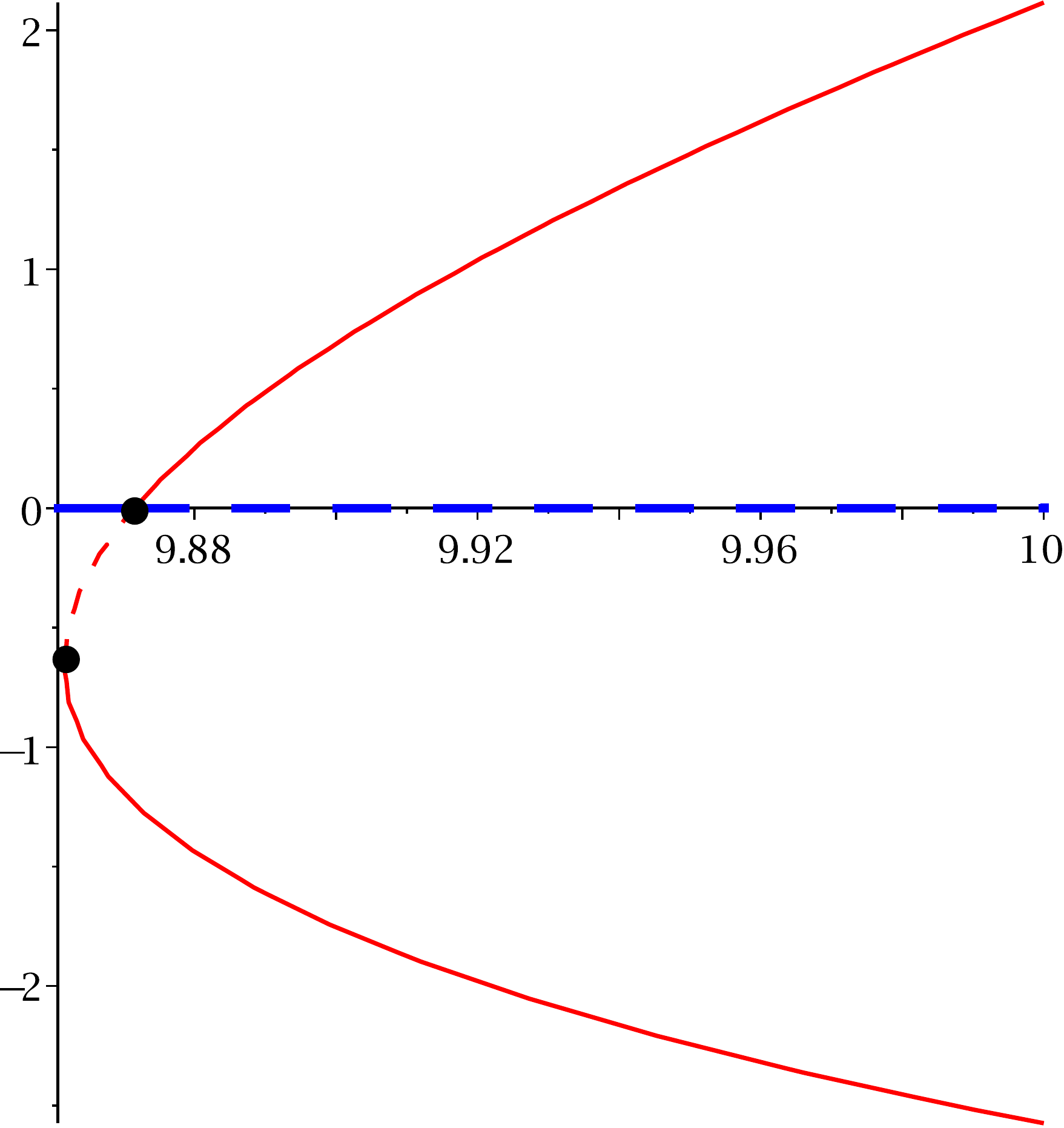}
\put(25,100){\tiny (4)}
\put(-3,50){\tiny $I$}
\put(25,45){\tiny$\Lambda$}
\put(10,56){\tiny Transcritical}
\put(8,41){\tiny Turning}
\end{overpic}
\hspace{.5cm}
\begin{overpic}[width=.27\textwidth,height=.29\textheight]{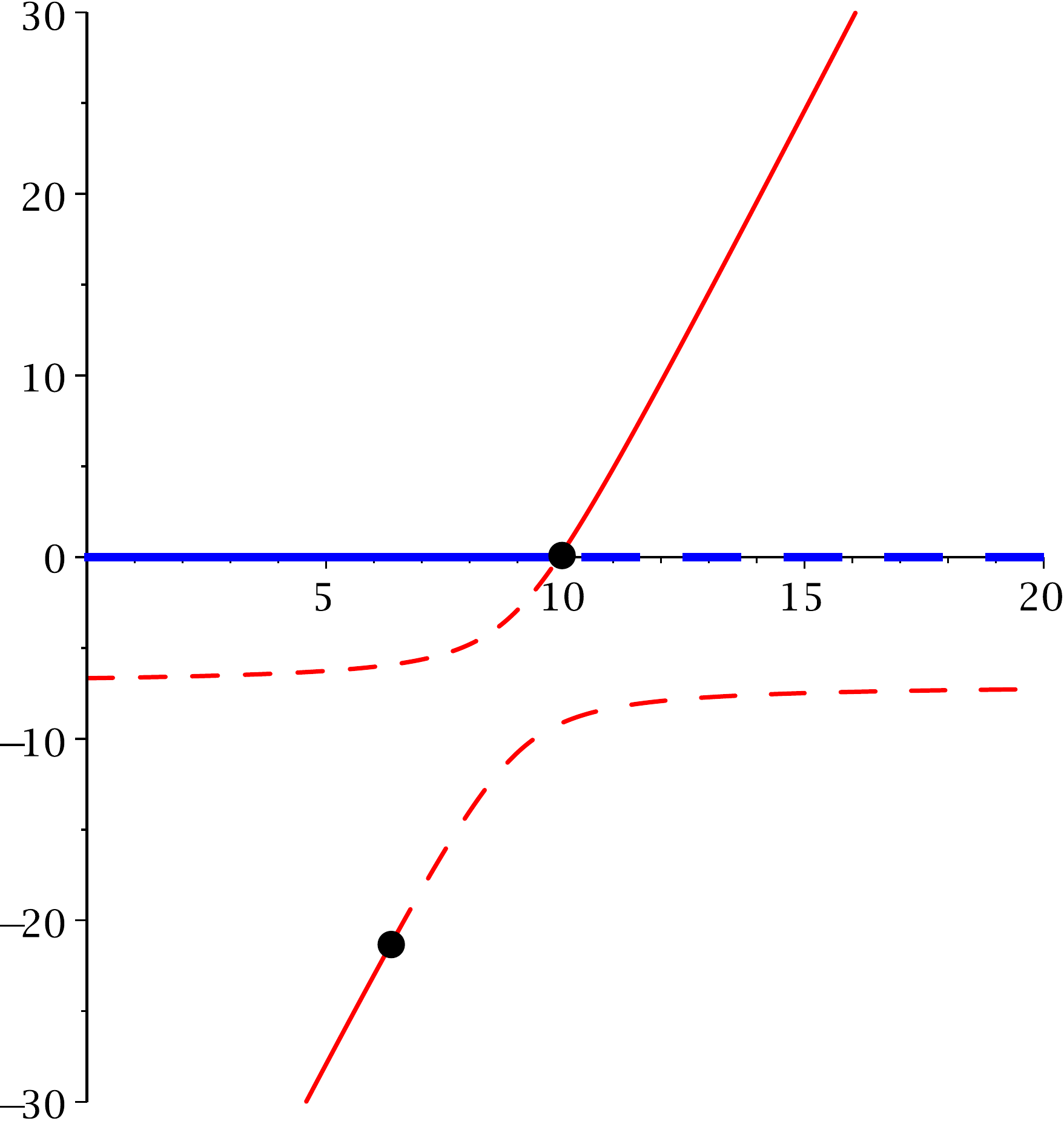}
\put(25,100){\tiny (5)}
\put(-3,50){\tiny $I$}
\put(25,40){\tiny$\Lambda$}
\put(30,51){\tiny Transcritical}
\put(20,15){\tiny Hopf}
\end{overpic}
\\ \vspace{1cm}
\caption{Bifurcation diagrams and simulations associated with 
the five cases given in Table~\ref{Paper3_Table1}, 
demonstrating various dynamical behaviors.} 
\label{Paper3_Fig3}
\end{figure}

We will show that complex dynamical behaviors can emerge in 
system (\ref{Paper3_Eq7}) 
from the infected equilibrium $\bar{\mathrm E}_1 = (\bar{S},\,\bar{I})$, 
where $\bar{I}$ is 
determined from the equation ${\mathcal F}(I) =0$ in (\ref{Paper3_Eq8}).
In the $\Lambda$-$I$ plane, the bifurcation diagram as shown in 
Figure~\ref{Paper3_Fig3} (1)-(4), indicates a turning point on the curve
with appropriate parameter values,
determined by both the quadratic equation~(\ref{Paper3_Eq8}) and the relation
$\frac{\mathrm{d}\Lambda}{\mathrm{d}I}=-\frac
{\partial {\mathcal F}}{\partial I}/
\frac{\partial {\mathcal F}}{\partial \Lambda}=0$, which is equivalent to
$\frac{\partial {\mathcal F}}{\partial I}=0$. Solving 
$\frac{\partial {\mathcal F}}{\partial I}=0$ yields the turning point of $I$, 
denoted by $I_T$ (`$T$' means \emph{turning}), taking the form
\begin{equation*}
\dss I_T = \frac{1}{2}\left[ \frac{\beta \Lambda_T}{(dk+\beta)(d+\epsilon)} 
-\omega - \frac{d}{dk+\beta} - \frac{\alpha}{d+\epsilon} \right],
\end{equation*}
where $\Lambda_T$ is obtained from $\mathcal{F}(I_T)=0$, see~(\ref{Paper3_Eq8}).
Thus, when $I_T>0\ (<0)$, the turning point of the quadratic curve 
appears above (below) the $I$-axis, meaning that 
backward bifurcation occurs for $I>0 \ (<0)$. 
Evaluating the Jacobian matrix at the infected equilibrium 
$\bar{\mathrm E}_1$, 
and further denoting it as $J_2=J|_{(\ref{Paper3_Eq7})}(\bar{\mathrm E}_1)$, 
we obtain the characteristic polynomial in the form of~(\ref{Paper3_Eq13}),
with $\mathrm{Tr}(J_2)=a_{11}/[(\omega+I)^2(kI+1)(dkI+\beta I+d)]$ and
$\mathrm{Det}(J_2)=a_{21}/[(\omega+I)^2(kI+1)(dkI+\beta I+d)]$,
where $a_{11}=a_{1a}-a_{1b}$ and $a_{21}=a_{2a}-a_{2b}$,
with $a_{1b}= \beta \Lambda (\omega+I)^2$ and $a_{2b}=d a_{1b}$, 
and $a_{1a}$ and $a_{2a}$ only contain positive terms 
(their expressions are omitted here for brevity).
Therefore, we can rewrite $\mathrm{Det}(J_2)= \frac{a_{21}}{d} d/
[(\omega+I)^2(kI+1)(dkI+\beta I+d)]$. 
Determining whether a Hopf bifurcation can occur from $\bar{E}$ is 
equivalent to finding whether $\mathrm{Det}(J_2)$ remains positive when
$\mathrm{Tr}(J_2)=0$.
Ignoring the positive factors in the 
following subtraction yields 
\begin{equation*}
\dss
h_1(I) = \frac{\mathrm{Tr}(J_2)-\mathrm{Det}(J_2)/d}
{(\omega+I)^2(kI+1)(dkI+\beta I+d)}=
a_{11}-\frac{1}{d}a_{21}=a_{1a}-\frac{1}{d}a_{2a},
\end{equation*}
where $h_1(I)=\frac{1}{d}(dkI+\beta I+d)[(kI+1)d^2(\omega+I)^2-
\beta \epsilon I (\omega+I)^2-\alpha \beta \omega I]$.
Thus, when $a_{1a}=0$, $\frac{1}{d}a_{2a}$ and $h_1(I)$
have opposite signs, implying that when $\mathrm{Tr}(J_2)=0$,
$\mathrm{Det}(J_2)$ could be positive only if $h_1(I)$ is negative.
Therefore, the necessary condition for system (\ref{Paper3_Eq7}) to 
have a Hopf bifurcation from the infected equilibrium $\bar{\mathrm E}_1$
is that $h_1(I)$ is negative.

In the remaining part of this subsection, we demonstrate
various dynamics which may happen in system~(\ref{Paper3_Eq7}) with different 
parameter values of $k$, as shown in Table~\ref{Paper3_Table1}. 
Taking other parameter values as $\alpha=6$, $\omega=7$, 
$\epsilon= 0.02$, $\gamma = 0.01$, $\beta = 0.01$, and $d=0.1$,
and solving the two equations $\mathrm{Tr}(J_2)=0$ 
and ${\mathcal F}(I) =0$ in~(\ref{Paper3_Eq8}) 
gives the Hopf bifurcation point candidates, $(\Lambda_H,\,I_H)$, 
for which $h_1(I_H)<0$.
Since the formula for the transcritical bifurcation point $\Lambda_S$ 
has no relation with $k$, $(\Lambda_S,\,I_S)=(9.87,\,0)$ 
is a fixed value pair in Table~\ref{Paper3_Table1}.
Bifurcation diagrams and associated numerical simulations
are shown in Figure~\ref{Paper3_Fig3} corresponding to 
the five cases given in Table~\ref{Paper3_Table1}. 
The blue lines and red curves represent the uninfected equilibrium
$\bar{E}_0$ and infected equilibrium $\bar{E}_1$, respectively.
The stable and unstable equilibrium solutions are shown by solid 
and dashed lines/curves, respectively. 
Backward bifurcation occurs in Cases $1$, $2$, and $3$ (see 
Table~\ref{Paper3_Table1}), which are illustrated by the corresponding 
bifurcation diagrams in 
Figures~\ref{Paper3_Fig3}(1), (2), and (3), respectively. 
For Cases $1$ and $2$, only one Hopf bifurcation occurs on the 
upper branch of the infected equilibrium $\bar{E}_1$, and this bifurcation
point exists at the critical point $\Lambda_H < \Lambda_S$ for Case $1$ and 
$\Lambda_H > \Lambda_S$, for Case $2$.
For Case $1$ with $\Lambda = 9.78$, 
the simulated time history converges to $\bar{\mathrm E}_0$ 
with initial condition IC$=[93.6,\, 0.44]$, 
shown in Figure~\ref{Paper3_Fig3}(1a), 
but converges to $\bar{\mathrm E}_1$ with initial condition
IC$=[46.8,\, 10]$, shown in Figure~\ref{Paper3_Fig3}(1b).
This clearly indicates the bistable behavior when $\Lambda_H < \Lambda_S$, and an overlapping stable region for 
both $\bar{\mathrm E}_0$ and $\bar{\mathrm E}_1$ exists
(see Figure~\ref{Paper3_Fig3}(1)).
The recurrent behavior for Case 2 is simulated at $\Lambda=9.87$ 
with IC$=[50,\,5]$, shown in Figure~\ref{Paper3_Fig3}(2a).
For Case 2, $\Lambda_H > \Lambda_S$,
and an overlapping unstable parameter region for both $\bar{\mathrm E}_0$ 
and $\bar{\mathrm E}_1$ occurs between $\Lambda_S$ and $\Lambda_H$
(see Figure~\ref{Paper3_Fig3}(2)).
For Case $3$, two Hopf bifurcations occur on the left side
of $\Lambda_S$, and a stable part in the upper branch of 
$\bar{\mathrm E}_1$ exists when $\Lambda$ passes through the critical value
$\Lambda=\Lambda_S$.
In this case, although backward bifurcation still exists and
the turning point is also located above the $\Lambda$-axis, giving
two branches of biologically feasible $\bar{\mathrm E}_1$, 
only regular oscillating behavior is observed.
The simulated time history is conducted at $\Lambda=10$,
with initial condition IC$=[50,\,2]$, shown in
Figure~\ref{Paper3_Fig3}(3a).
For Case $4$, only forward bifurcation occurs 
in the biologically feasible region, and the turning point for
backward bifurcation moves down to the fourth quadrant, that is,
negative backward bifurcation occurs in this case. 
The whole upper branch of $\bar{\mathrm E}_1$ in the first quadrant 
is stable, therefore, no oscillations (or recurrence) can happen.
Finally, further increases to the value of $k$ change the shape of the red 
curves, as shown in Figure~\ref{Paper3_Fig3}(5), which again 
indicates that no biologically meaningful backward bifurcation 
or oscillations can occur. 
Note that in Figure~\ref{Paper3_Fig3}(5) a Hopf bifurcation point 
exists on the lower branch of the equilibrium solution, which is 
biologically unfeasible since it is entirely below the horizontal axis. 
In conclusion, interesting dynamical behaviors can emerge in system
(\ref{Paper3_Eq7}) if backward bifurcation occurs.

\begin{table}
\caption{\label{Paper3_Table1}
Dynamics of system (\ref{Paper3_Eq7}) for different values of $k$, 
with $\alpha=6$, $\omega=7$, $\epsilon= 0.02$, $\gamma = 0.01$, 
$\beta = 0.01$, $d=0.1$, and a fixed transcritical bifurcation point
$(\Lambda_S,\,I_S)=(9.87,\,0)$.}
\begin{tabular}{|c|p{0.8cm}|p{2.3cm}|p{2.7cm}|p{2.5cm}|p{1.8cm}|p{2.3cm}|}
\hline
\!\!\! Case \!\!\! & \quad $k$ & \quad $(\Lambda_T,\,I_T)$ &  
\quad $h_1(I)<0$ &  $(\Lambda_H,\,I_H)$ 
& Dynamics &  Notes\\ \hline
$1$ & $0.001 $ 
& $(9.48,\,4.57)$ & $I\in[1.72,\,\infty]$ & $(9.73,\,10.28)$ 
& Bistability & $\Lambda_H<\Lambda_S$\\ \hline
$2$ & $0.01$ & $(9.71,\,2.82)$ & $I\in[1.76,\,\infty]$ & $(9.96,\,8.00)$ 
& Recurrence & $\Lambda_H > \Lambda_S$\\ \hline 
$3$ & $0.02$ & $(9.85,\,0.84)$ & $I\in[1.82,\,\infty]$ & $(9.88,\,2.09)$,
$(10.14,\,5.62)$ & Oscillation & 
$\begin{array}{ll} \\[-0.5ex]  
\! {\rm Two \ Hopf}\\ 
\! {\rm critical \ points} \end{array}$  \\ \hline
$4$ & $0.027 $ & $(9.86,\,-0.65) \hspace{-0.20in}$ 
&$I\in[1.85,30.65]$ & No Hopf &
$\begin{array}{ll} \\ [-0.5ex] 
\! {\rm No} \\ [-0.5ex] 
\! {\rm oscillation} \end{array}$ 
& $\begin{array}{ll}
\!{\rm Negative} \\[-0.5ex] 
\!{\rm backward} \\[-0.5ex] 
\!{\rm  bifurcation} \\[-0.5ex]  
\end{array}$  
\\ \hline
$5$ & $0.05$ & No Turning & $I\in[2.01,15.03]$ &
$(6.18,-22.15)$ & 
$\begin{array}{ll} \\ [-0.5ex] 
\! {\rm No} \\ [-0.5ex] 
\! {\rm oscillation} \end{array}$ 
& $\begin{array}{ll} 
\\[-0.5ex] 
\! {\rm No \ backward} \\[-0.5ex] 
\! {\rm bifurcation} \end{array} $ \\
\hline
\end{tabular}
\end{table}

\subsection{Hopf bifurcation in the infection model 
with convex incidence}

In this subsection, we return to system~(\ref{Paper3_Eq4}), that is,
the 2-dimensional HIV model 
with convex incidence derived in~\cite{ZWY2013,ZWY2014}, and
analyze the various dynamical phenomena which
system~(\ref{Paper3_Eq4}) could possibly exhibit. 
To achieve this, we set $B$ as the bifurcation parameter, and 
$A$ as a control parameter; the bifurcation analysis will be carried out 
for various values of $A$. 
Also, simulated time histories are provided to illustrate 
the dynamical behavior predicted in the analysis.

We first consider the uninfected equilibrium $\bar{E}_0=(\frac{1}{D},\,0)$,
which has two eigenvalues. One of them, given by 
$\lambda_1|_{\bar{E}_0}=-D$, is always negative. 
The other one is $\lambda_2|_{\bar{E}_0}=\frac{B}{D}-1$. 
Thus, depending upon the relation between $B$ and $D$,
$\lambda_2|_{\bar{E}_0}=0$ gives a static bifurcation at $B_S=D$
(or $R_0=\frac{B}{D}=1$),
which is further proved to be a transcritical bifurcation.
Here the `S' in subscript stands for \emph{static bifurcation}.
Therefore, $\bar{E}_0$ is stable when $B<D$ (or $R_0<1$), 
loses its stability and becomes unstable when $B$ increases 
to pass through $B_S=D$, that is $B>D$ (or $R_0>1$), 
and no other bifurcations can happen.

Next, we examine the infected equilibrium $\bar{E}_1=(\bar{X},\,\bar{Y})$.
Since $\bar{X}(Y)=\frac{Y+C}{(A+B)Y+BC}$, $\bar{Y}$ is determined by
the quadratic equation~(\ref{Paper3_Eq10}), which gives the
turning point $(B_T,\,Y_T)$ as
\begin{equation*}
B_T = \frac{-A+D+2\sqrt{ACD}}{C+1},\quad 
Y_T = \frac{A+B-BC-D}{A+B},
\end{equation*}
where `$T$' in the subscript stands for \emph{turning bifurcation}.
We perform a further bifurcation analysis on its corresponding characteristic 
polynomial~(\ref{Paper3_Eq13}), which takes the form
\begin{equation}
\begin{array}{ll}
\dss P|_{\bar{E}_1}(\lambda,Y) = \lambda^2+\frac{a_{1a}}{[(A+B)Y+BC](Y+C)}
\lambda+\frac{a_{2a}}{[(A+B)Y+BC](Y+C)},\quad \text{where}\\[2.0ex]
a_{1a} = (A+B)^2Y^3+(2BC+D)(A+B)Y^2+(B^2C^2+ACD+2BCD-AC)Y+BC^2D, \\[.5ex]
a_{2a} = (A+B)^2Y^3+2(A+B)BCY^2+(B^2C-AD)CY.
\end{array}
\label{Paper3_Eq14}
\end{equation}
Therefore, the sign of the subtraction between the trace and determinant is
determined by
$h_2(Y)=a_{1a}-a_{2a}=D(A+B)Y^2+[2CD(A+B)-AC]Y+BC^2D$.
Here the equilibrium solution of $Y$ and other parameters satisfy 
the quadratic equation~(\ref{Paper3_Eq10}),
which leads to an explicit expression, given by
$\bar{B}= -\frac{AY^2+(D-A)Y+CD}{Y^2+(C-1)Y-C}$.
Substituting $B=\bar{B}$ into $h_2(Y)$, we obtain 
\begin{equation*}
\dss h_2(Y)|_{B=\bar{B}}=  a_{1a} - a_{2a}
= \frac{[AC(D-1)-D^2]Y^2-[AC(D-1)+2CD^2]Y-C^2D^2}{Y-1}.
\end{equation*}
Hopf bifurcation may occur when the trace is zero, 
while the determinant is still positive.
This implies $h_2(Y)<0$, which is possible with appropriately 
chosen parameter values.
Hence, by solving $a_{1a}=0$ in~(\ref{Paper3_Eq14}) together with 
the quadratic equation~(\ref{Paper3_Eq10}), we get two pairs of points
denoted by $(B_{h1},\,Y_{h1})$ and $(B_{h2},\,Y_{h2})$, which are candidates
for Hopf bifurcation. 
Then validating the above two points by substituting them back into 
the characteristic polynomial~(\ref{Paper3_Eq14}), respectively, 
we denote the Hopf bifurcation point as $(B_H,\,Y_H)$ if 
this validation confirms their existence. 
According to~(Yu et al., submitted for publication), 
Hopf bifurcation can 
happen only from the upper branch of the infected equilibrium $\bar{E}_1$.

The various dynamical behaviors which may appear in system~(\ref{Paper3_Eq4}) 
have been classified in Table~\ref{Paper3_Table2} for different values of 
the parameter $A$, with fixed values of $C=0.823$ and $D=0.057$. 
Thus, the transcritical bifurcation point is fixed 
for all cases: $B_S=D=0.057$ and $Y_S=0$. 
The two solutions $B_{h1}$ and $B_{h2}$ are 
solved from the two equations (\ref{Paper3_Eq14}) $P|_{\bar{E}_1}(\lambda,Y) =0$ 
and (\ref{Paper3_Eq10}) ${\mathcal F}_5(Y)=0$, respectively. 
They become a Hopf bifurcation point only if 
their corresponding $Y$ values
($Y_{h1}$ and $Y_{h2}$, respectively) are in the range such that
$h_2(Y)<0$. Otherwise, system (\ref{Paper3_Eq4}) has a pair of real eigenvalues
with opposite signs at $(B_{h1},\,Y_{h1})$ or $(B_{h2},\,Y_{h2})$,
which is denoted by the superscript `${\ast}$' (which is actually 
a saddle point) in Table~\ref{Paper3_Table2}, while the Hopf bifurcation
point is denoted by the superscript `$H$' in Table~\ref{Paper3_Table2}.

\begin{table}
\caption{\label{Paper3_Table2} Parameter values taken to illustrate various dynamics of system
(\ref{Paper3_Eq4}). The fixed transcritical bifurcation point: $(B_S,\,Y_S)=(0.057,\,0)$}
\begin{tabular}{|c|p{0.8cm}|p{3.5cm}|p{3.5cm}|p{3.5cm}|}
\hline
Case&\ $A$ & \hspace{0.22in} $(B_T,\,Y_T)$ & $h_2(Y)<0,\,Y \! \in$ & 
\hspace{0.16in} $(B_{h1},\,Y_{h1})$ \\ \hline
$1$& $0.80$ & $(-0.1950,\,0.5850)$ & $(0.0036,\,0.9830)$ 
& $(0.0355,\,0.8725)^H$ \\ \hline
$2$ & $0.71$ & $(-0.1580,\,0.5660)$ & $(0.0040,\,0.9800)$ & 
$(0.0539,\,0.0038)^{\ast}$ \\ \hline
$3$ & $0.60$ & $(-0.1140,\,0.5380)$ & $(0.0048,\,0.9769)$ & 
$(0.0540,\,0.0045)^{\ast}$ \\ \hline 
$4$ & $0.07$ & $(0.0557,\,0.0909)$ & $(0.0476,\,0.8030)$ &
$(0.0560,\,0.0470)^{\ast}$ \\ \hline
$5$ & $0.06$ & $(0.056558,\,0.05581)$ & $(0.0574,\,0.7700)$ &
$(0.056559,\,0.0574)^H$ \\ \hline
$6$ & $0.05$ & $(0.05697,\,0.01442)$ & $(0.0724,\,0.7232)$ &
$(0.0574,\,0.0741)^H$ \\ \hline
$7$ & $0.04$ & $(0.0569,\,-0.0358)$ & $(0.0986,\,0.6507)$ &
$(0.0592,\,0.1071)^H$ \\ \hline
$8$ & $0.03$ & $(0.0559,\,-0.0994)$ & $(0.1611,\,0.5149)$ &
\hspace{0.40in} ---  \\ \hline\hline
Case&\ $A$ & \hspace{0.16in} $(B_{h2},\,Y_{h2})$ 
&  Dynamics & Notes\\ \hline
$1$& $0.80$ & 
$(0.054,\,0.0034)^{\ast}$ & Unstable limit cycle, Bistable 
& $B_{h1}<B_S$ \\ \hline
$2$ & $0.71$ & $(0.0574,\,0.8650)^H$ & Recurrence &
$B_{h2}>B_S$ \\ \hline
$3$ & $0.60$ & $(0.0819,\,0.8530)^H$ & Recurrence & 
$B_{h2}>B_S$ \\ \hline
$4$ & $0.07$ & $(0.1015,\,0.5612)^H$ & Recurrence &
$B_{h2}>B_S$ \\ \hline 
$5$ & $0.06$ & $(0.0961,\,0.5225)^H$ & Recurrence &
$ B_{h1}< B_S <B_{h2}$ \\ \hline
$6$ & $0.05$ & $(0.0894,\,0.4701)^H$ & Recurrence &
$ B_{h1}< B_S <B_{h2} $ \\ \hline
$7$ & $0.04$ & $(0.0806,\,0.3897)^H$ & Oscillation &
$ B_{h1}< B_S <B_{h2}$, $ Y_T < 0  $ \\ \hline
$8$ & $0.03$ & \hspace{0.40in} --- & $\bar{E}_1$ stable & $Y_T<0$ 
\\ \hline
\end{tabular}
\end{table}

Next, we further examine the direction of the Hopf bifurcation, 
that is, check whether it is a supercritical or subcritical 
Hopf bifurcation. Since the Jacobian matrix of the system 
evaluated at the Hopf bifurcation point 
has a pair of purely imaginary eigenvalues, 
the linearized system (\ref{Paper3_Eq4}) does not determine the 
nonlinear behavior of the system.  
Therefore, we take advantage of normal form theory 
to study the existence of the limit cycles bifurcating from 
the Hopf bifurcation point as well as their stability.
As mentioned earlier, Hopf bifurcation can only occur from the
upper branch of the infected equilibrium $\bar{E}_1$, therefore we first
transform the fixed point $\bar{E}_1$ to the origin by 
a shifting transformation, and, 
in addition, make the parameter transformation 
$B=B_H+\mu$; the Hopf bifurcation point is thus defined as 
$\mu=\mu_H=0$. Then the normal form of system 
(\ref{Paper3_Eq4}) near the critical point, $\mu=\mu_H=0$, takes the form 
up to third-order approximation: 
\begin{equation}
\dot{r}=d\,\mu \,r+a\,r^3 +\mathcal{O}(r^5),
\qquad
\dot{\theta} = \omega_c + c\,\mu + b\,r^2+\mathcal{O}(r^4), 
\label{Paper3_Eq15}
\end{equation}
where $r$ and $\theta$ represent the amplitude and phase of the motion, 
respectively. The first equation of (\ref{Paper3_Eq15}) can be used for 
bifurcation and stability analysis, while the second equation of 
(\ref{Paper3_Eq15}) can be used to determine the frequency of the 
bifurcating  periodic motions.  
The positive $\omega_c$ in the second equation of (\ref{Paper3_Eq15}) 
is the imaginary part of the eigenvalues at the Hopf bifurcation point. 
The parameters $d$ and $c$ can be easily obtained from 
a linear analysis, while $a$ and $b$ must be derived using a nonlinear
analysis, with the Maple program available in, say, \cite{Yu1998}.

Note that the infected equilibrium $\bar{E}_1$ is represented by 
the fixed point $\bar{r}=0$ of system (\ref{Paper3_Eq15}),
while the nonzero fixed point $\bar{r} > 0$ 
(satisfying $\bar{r}^2=\frac{-d \mu}{a}$) 
is an approximate solution for a limit cycle or periodic orbit. 
The periodic orbit is asymptotically stable (unstable) if $a<0$ ($a>0$), and 
the corresponding Hopf bifurcation is called supercritical (subcritical). 
According to the Poincare-Andronov Hopf Bifurcation theorem~\cite{Wiggins1990},
for $\mu$ sufficiently small, there are four possibilities for the 
existence of periodic orbits and their stability, which are classified in
Table~\ref{Paper3_Table3}, based on the four
sets of the parameter values in the normal form~(\ref{Paper3_Eq15}).
Then we use the results presented in Table~\ref{Paper3_Table3} with 
a nonlinear analysis based on normal form theory 
to classify the Hopf bifurcations appearing in Table~\ref{Paper3_Table2}, 
and the results are shown in Table~\ref{Paper3_Table4}.

\begin{table}
\caption{\label{Paper3_Table3} Classification of 
Hopf bifurcations based on the normal form (\ref{Paper3_Eq15}).} 
\begin{tabular}[\textwidth]{|c|c|c|c|c|c|}
\hline
Class & \multicolumn{2}{|c|}{Stability of $\bar{r}=0$} & 
\multicolumn{2}{|c|}{Stability of $\bar{r}^2 \!=\! -\frac{d\mu}{a}$} &
Hopf bifurcation \\ \cline{2-5}  
 & $\mu < 0$ & $\mu > 0$ & $\mu < 0$ & $\mu > 0$ &  \\ \hline
(a): $d>0$, $a>0$ & stable & unstable & unstable & --& subcritical \\
\hline
(b): $d>0$, $a<0$ & stable & unstable & -- & stable & supercritical \\
\hline
(c): $d<0$, $a>0$ & unstable & stable & --& unstable& subcritical \\
\hline
(d): $d<0$, $a<0$ & unstable& stable&stable& -- & supercritical \\
\hline
\end{tabular}
\end{table}

\begin{table} 
\caption{\label{Paper3_Table4} Classification of 
Hopf bifurcations appearing in Table~\ref{Paper3_Table2}.} 
\vspace{0.05in} 
\begin{tabular}{|c|p{1cm}|p{3.5cm}|p{1.5cm}|c|p{2cm}|p{1.0cm}|}
\hline
Case & $A$ & Hopf bifurcation point $(B_H,\,Y_H)$ & $d$ & $a$ 
& Stability of limit cycles &
Table~\ref{Paper3_Table3} class\\
\hline
$1$ & $0.8$ & $(0.0355,\,0.8725)$ & $-1.0722$ &\hspace{0.10in}$0.2114\times 10^{-2}$
& Unstable & (c)\\
\hline
$2$ & $0.71$ & $(0.0574,\,0.8650)$ & $-1.0726$ &\hspace{0.10in}$0.1424\times 10^{-2}$ 
& Unstable & (c) \\
\hline
$3$ & $0.6$ & $(0.0819,\,0.8530)$ & $-1.0733$ &\hspace{0.10in}$0.6755\times 10^{-3}$
& Unstable & (c) \\
\hline
$4$ & $0.07$ & $(0.1015,\,0.5612)$ & $-1.0307$ & $-0.8791\times 10^{-3}$
& stable & (d) \\
\hline
$5$ & $0.06$  \multirow{2}{*}& $(0.056559,\,0.0574)$ &\hspace{0.10in}$884.27$ 
&$-0.1019$ \hspace{0.41in} & Stable & (b)  \\[-0.0ex] \hline
& & $(0.0961,\,0.5225)$ &$-1.0079$ & $-0.8613\times 10^{-3}$ & 
Stable& (d) \\ 
\hline
$6$ & $0.05$ \multirow{2}{*}& $(0.0574,\,0.0741)$ &\hspace{0.10in}$18.232$
& $-0.3145\times 10^{-2}$ & Stable & (b)  \\[-0.0ex] \hline
& &  $(0.0894,\,0.4701)$ & $-0.9629$ & $-0.8457\times 10^{-3}$ 
& Stable & (d) \\
\hline
$7$ & $0.04$ \multirow{2}{*}& $(0.0592,\,0.1071)$ &\hspace{0.10in}$4.7242$ & 
$ -0.1577\times 10^{-2}$ & Stable & (b) \\[-0.0ex] \hline
& & $(0.0805,\,0.3897)$  & $-0.8437$ & $ -0.8438\times 10^{-3}$ &
Stable & (d) \\
\hline
\end{tabular}
\end{table}

To illustrate the analytical results given in Tables~\ref{Paper3_Table2} 
and~\ref{Paper3_Table4}, we provide the bifurcation diagrams in 
Figures~\ref{Paper3_Fig4} (1)-(8).
These figures depict the uninfected equilibrium 
$\bar{E}_0$ and the infected equilibrium $\bar{E}_1$ 
in blue and red, respectively.
The solid and dashed lines
differentiate stable and unstable states of the equilibrium solutions.
The bifurcation points on the equilibrium solutions are highlighted 
by solid black dots. 
Moreover, `Transcritical', `Turning',
`Hopf$_{\text{sub}}$', and `Hopf$_{\text{super}}$', are used to 
denote
\emph{Transcritical bifurcation}, \emph{Turning point},
\emph{subcritical Hopf bifurcation}, and \emph{supercritical Hopf bifurcation},
respectively.
Simulated time histories are used to validate the analytical results,
and to show different dynamical behaviors in each case listed in 
Tables~\ref{Paper3_Table2} and~\ref{Paper3_Table4}.
Subcritical Hopf bifurcation occurs in Cases 1-3, shown in 
Figures~\ref{Paper3_Fig4} (1)-(3).
$A=0.8$ is used in Figure~\ref{Paper3_Fig4} (1) for Case 1. 
Choosing $B=0.036$, we have $E_0=[17.1282566,\,0.023689]$
and $E_1=[2.233533,\,0.8726886]$.
The simulated solution converges to ${\mathrm E}_0$ or ${\mathrm E}_1$, 
with initial condition taken as 
IC$_d=[17.13,\,0.024]$ or IC$_c=[2.233,\,0.873]$, 
shown in Figures~\ref{Paper3_Fig4} (1d) and (1c), respectively.
Figures~\ref{Paper3_Fig4} (1a) and (1b), on the other hand, 
show the unstable limit cycle
bifurcating from the subcritical Hopf bifurcation with IC$_c=[2.233,\,0.873]$.

Figure~\ref{Paper3_Fig4} (2) corresponds to Case 2 with $A=0.71$.
Choosing $B=0.0572\in[B_S,\,B_H]$ yields recurrence, independent of 
the initial conditions, see, for example, 
the result given in Figure~\ref{Paper3_Fig4} (2b) with IC$_b= [2.4,\,0.5]$.
However, for $B=0.06>B_H$, the simulated time history
converges to $E_1$, with an initial condition close to $E_1$, such as
IC$_a=[2.4,\,0.6]$ as shown in Figure~\ref{Paper3_Fig4} (2a);
or shows recurrence with an initial condition far away from $E_1$, such as
IC$_c= [2.4,\,0.4]$, as shown in Figure~\ref{Paper3_Fig4} (2c).

Figure~\ref{Paper3_Fig4} (3) plots the result for Case 3 with $A=0.6$, and 
shows a broader region between the transcritical and Hopf bifurcation points,
associated with a larger recurrent region.
Recurrence occurs independent of the initial conditions for
$B=0.083\in[B_S,\,B_H]$, giving 
$E_0=[12.048,\,0]$ and $E_1=[2.576,\,0.852]$,
as shown in Figures~\ref{Paper3_Fig4} (3a) and (3b), with
IC$_a=[2.7,\,0.84]$ and IC$_b=[14,\,0.1]$, respectively.
But if we choose $B=0.07>B_H$, we have $E_0=[14.286,0]$ and 
$E_1=[2.67,\,0.8478]$. The time history converges to $E_1$
with IC$_c=[2.6,\,0.8]$, or shows recurrence with IC$_d=[2.6,\,0.1]$,
as shown in Figure~\ref{Paper3_Fig4} (3c) and (3d), respectively.

Supercritical Hopf bifurcations occur in Cases 4-7, as shown in
Figures~\ref{Paper3_Fig4} (4)-(7). 
Figure~\ref{Paper3_Fig4} (4) depicts the result for Case $4$ with $A=0.07$.
Only one supercritical Hopf bifurcation happens in this case, and
gives a large recurrent parameter region between the transcritical and 
Hopf bifurcation points.
Although the simulated recurrent behavior does not depend on initial 
conditions, the recurrent pattern will fade out with the growth of the value 
of $B$ from the transcritical point to the Hopf bifurcation point, 
see Figures~\ref{Paper3_Fig4} (4a) and (4b) with the same 
IC$_{a,\,b}=[8,\,0.1]$, but different values of $B$:
$B=0.06$ and $B=0.09$, respectively.

Figure~\ref{Paper3_Fig4} (5) shows the result for Case 5 with $A=0.06$.
A transcritical bifurcation happens between two supercritical Hopf bifurcations.
The recurrent region still starts from the transcritical point and independent 
of the initial conditions, but is narrower than that shown in 
Figure~\ref{Paper3_Fig4} (4).
The simulated recurrent behavior for this case 
is conducted at IC$=[12,\,0.1]$ and $B=0.06$.
Figure~\ref{Paper3_Fig4} (6) corresponds to Case 6 with $A=0.05$, and 
two supercritical Hopf bifurcations occur on the right side of 
the transcritical bifurcation point, which makes the recurrent 
region even narrower and the recurrent pattern less obvious, 
as shown in the simulated time
history with IC$=[10,\,0.1]$ and $B=0.06$.
Negative backward bifurcations occur in Cases 7 and 8, as shown in 
Figure~\ref{Paper3_Fig4} (7) and (8).
Although two Hopf bifurcations are still present in Case 7, see
Figure~\ref{Paper3_Fig4} (7), only a regular oscillating pattern exists.
For Case 8, no Hopf bifurcation happens in the biologically feasible
part of $E_1$, and therefore no more interesting dynamics occur.

In general, backward bifurcation, which occurs above the horizontal axis,
is much more likely to induce Hopf bifurcation.
A Hopf bifurcation can only occur along the upper branch of $\bar{E}_1$,
since $\bar{E}_0$ only changes its stability at
a transcritical bifurcation point, and any point on  
the lower branch of $\bar{E}_1$ is a saddle 
node~(Yu et al., submitted for publication). 
Moreover, Hopf bifurcation can lead to a change in the stability of the upper 
branch of the infected equilibrium $\bar{E}_1$.
Thus the system further develops 
bistable, recurrent, or regular oscillating behavior, corresponding to 
Cases $1-7$ in Tables~\ref{Paper3_Table2} and~\ref{Paper3_Table4}, 
and in Figures~\ref{Paper3_Fig4} (1)-(7).
In particular, bistability happens
when both equilibria $\bar{E}_0$ and $\bar{E}_1$ share a 
stable parameter region, see Case $1$ in Table~\ref{Paper3_Table2} and
Figure~\ref{Paper3_Fig4} (1).

As for recurrent behavior, we observe that recurrence is
more likely to happen if the following
two conditions are satisfied for the upper branch of $\bar{E}_1$:
(1) the equilibrium remains unstable as the bifurcation parameter increases
and crosses 
the trancritical point, where $\bar{E}_0$ and $\bar{E}_1$ intersect, 
such that the two equilibria share an unstable parameter range; and 
(2) at least one Hopf bifurcation occurs from $\bar{E}_1$. 
As shown in Cases $2$-$6$ in Table~\ref{Paper3_Table2}, 
and the corresponding Figures~\ref{Paper3_Fig4} (2)-(5),
the common recurrent parameter region for both subcritical and
supercritical Hopf bifurcations starts beside the transcritical point,
and is located entirely in the unstable parameter region of 
$\bar{E}_0$ and $\bar{E}_1$. 
The simulated recurrent pattern becomes more pronounced
if the value of the bifurcation parameter is close to the transcritical point,
but approaches an oscillatory pattern as the parameter diverges
from the transcritical point,
as shown in Figure~\ref{Paper3_Fig4} (4a) and (4b).
In this common recurrent parameter region, 
recurrence occurs independent of initial conditions;
see Figures~\ref{Paper3_Fig4} (3a) and (3b).
In addition to the common recurrent region, for subcritical bifurcation,
seen in Table~\ref{Paper3_Table2} for Cases (2) and (3) 
and Figures~\ref{Paper3_Fig4} (2) and (3),
recurrence may also appear on the stable side of the subcritical 
Hopf bifurcation point with an initial condition close to $\bar{E}_1$.
Moreover, the subcritical Hopf bifurcation and the transcritical point
should be close to each other for a clear recurrent pattern.
When this is not the case, the periodic solutions 
show a more regular oscillating pattern,
as compared in Figures~\ref{Paper3_Fig4} (2c) and (3d).
Although two Hopf bifurcation points occur in Table~\ref{Paper3_Table2} 
for Case $5$, see Figure~\ref{Paper3_Fig4} (5),
the transcritical point is located inside the unstable range of the upper
branch of $\bar{E}_1$, between the two Hopf bifurcation points.
A recurrent pattern still characterizes the dynamical behavior 
in this case.
However, if the unstable range of $\bar{E}_1$, between
the two Hopf bifurcation points, is located entirely in the unstable range of 
$\bar{E}_0$, and moves further away from the transcritical point, 
the recurrent motion gradually becomes a regular oscillation,
as shown in Figures~\ref{Paper3_Fig4} (6) and (7).

Summarizing the results and discussions presented in the previous 
two sections, we have the following observations. 
\begin{enumerate} 
\item 
Due to the fact that $\bar{E}_0$ only changes its stability at
the transcritical bifurcation point, and the fact that any point on the 
lower branch of $\bar{E}_1$ is a saddle node, 
Hopf bifurcation can only occur from the upper
branch of $\bar{E}_1$.  A Hopf bifurcation may result in convergent, recurrent, bistable, or 
regular oscillating behaviors. 

\item 
Backward bifurcation gives rise to two branches in the infected equilibrium
$\bar{E}_1$. Hopf bifurcation is more likely to happen when the turning 
point of the backward bifurcation is located on the positive part of the 
equilibrium solution in the bifurcation diagram,
as shown in Figures~\ref{Paper3_Fig4} (2)-(6).
This means that 
we have two biologically feasible infected equilibria, 
which is essential to observe bistability, as shown in 
Figure~\ref{Paper3_Fig4} (1).

\item 
However, if the turning point on the infected equilibrium $\bar{E}_1$, or
the backward bifurcation moves down to the negative part of a state 
variable in the bifurcation diagram, that is, negative backward bifurcation
occurs, then Hopf bifurcation is very unlikely to happen.
Although Figure~\ref{Paper3_Fig4} (7) shows an exceptional case, 
the parameter range for such a Hopf bifurcation is very narrow.

\item 
The bifurcation diagram for system (\ref{Paper3_Eq4}) with $A=0.03$, 
shown in Figure~\ref{Paper3_Fig4} (8),
is a typical model with negative backward bifurcation.
Such negative backward bifurcation may occur in higher-dimensional 
systems. However, by considering more state variables, 
which make the system more complicated, Hopf bifurcation can happen 
in the upper branch of the negative backward bifurcation. 
We will discuss this possibility in more detail in the next section 
by examining an autoimmune disease model.
\end{enumerate}

The results obtained in this section suggest the following summary. 

\begin{remark}
If a disease model contains a backward bifurcation 
on an equilibrium solution, 
then as the system parameters are varied, there may exist none, 
one or two Hopf bifurcations from the equilibrium solution, 
which may be supercritical or subcritical.
If further this equilibrium has a transcritical bifurcation point at 
which it exchanges its stability with another equilibrium, then 
recurrence can occur between the transcritical 
and Hopf bifurcation points and near the transcritical point, 
where both equilibrium solutions are unstable, 
and bistability happens when Hopf bifurcation makes a shared
stable parameter region for both equilibria.
\label{Paper3_Remark2} 
\end{remark}

\begin{figure}
\centering
{\tiny
\hspace{-0.3in}
\subfloat{
\begin{overpic}[width=0.4\textwidth]{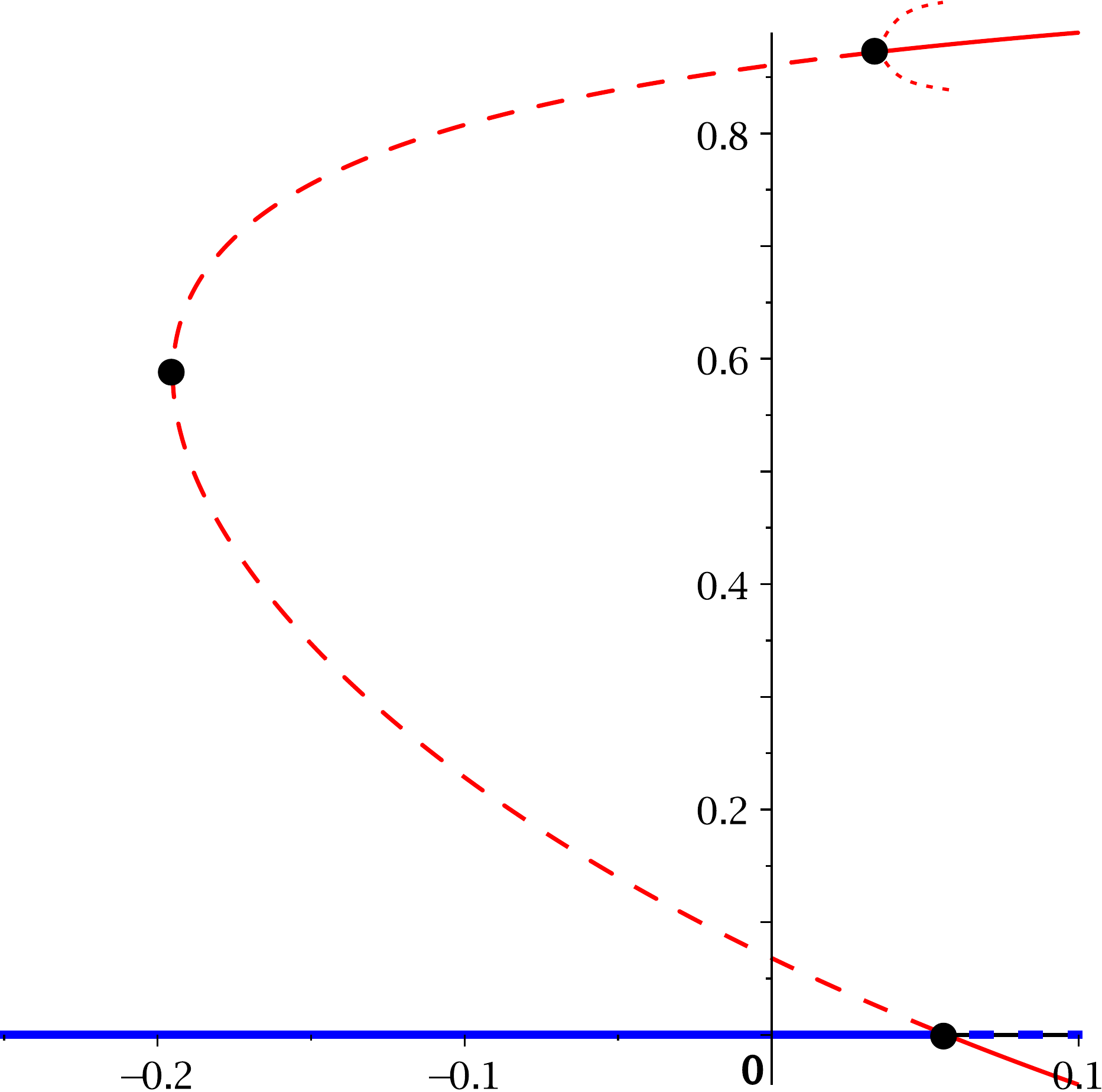}
\put(50,100){(1)}
\put(50,-2){$B$}
\put(65,55){$Y$}
\put(0,65){\tiny Turning}
\put(65,95){\tiny Hopf$_{\text{sub}}$}
\put(88,6){\tiny Transcritical}
\put(81,91){$\ast$}
\put(38,60){\includegraphics[width=0.1\textwidth,height=0.07\textheight]
{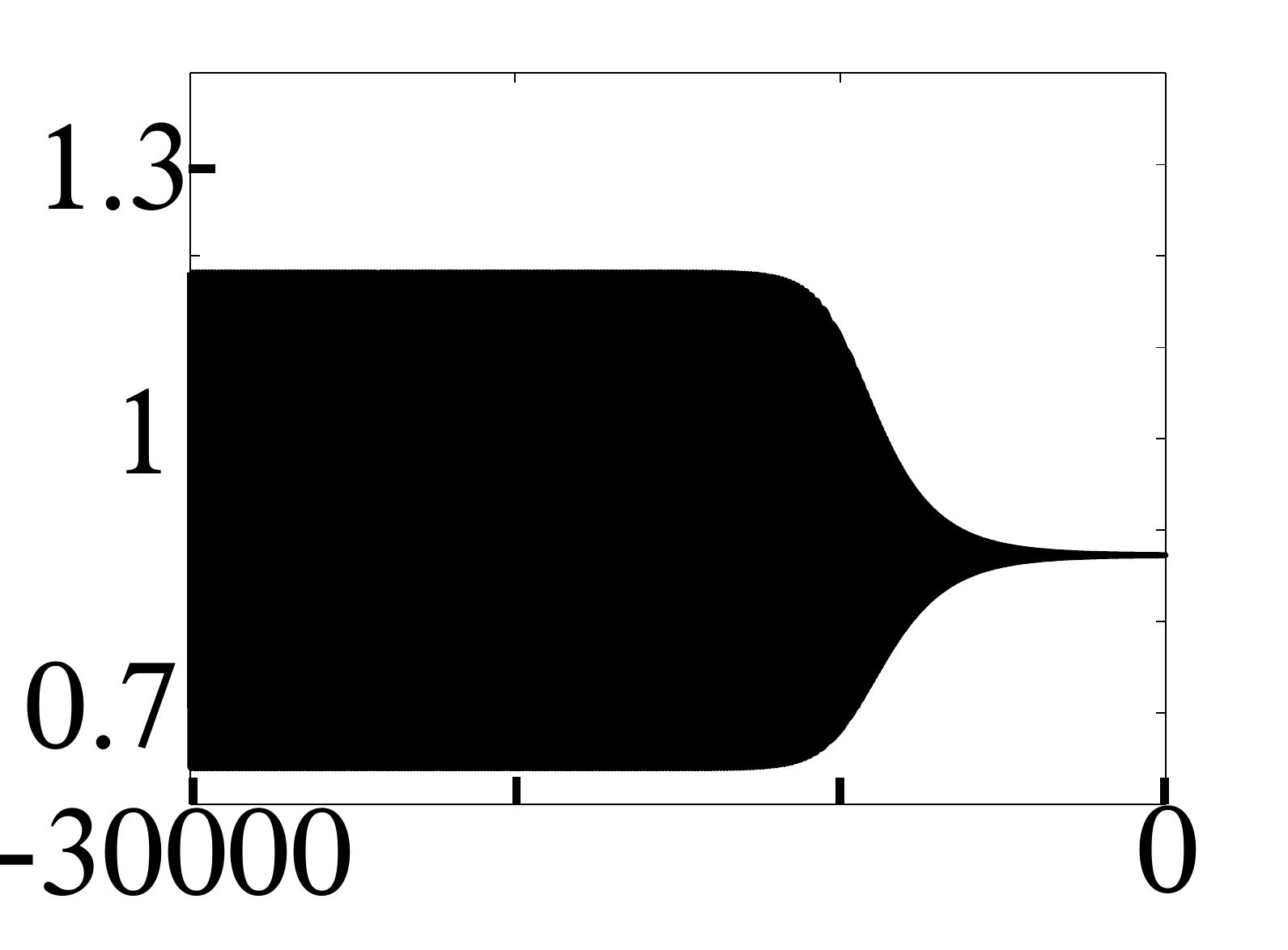}}
\put(38,33){\includegraphics[width=0.1\textwidth,height=0.07\textheight]
{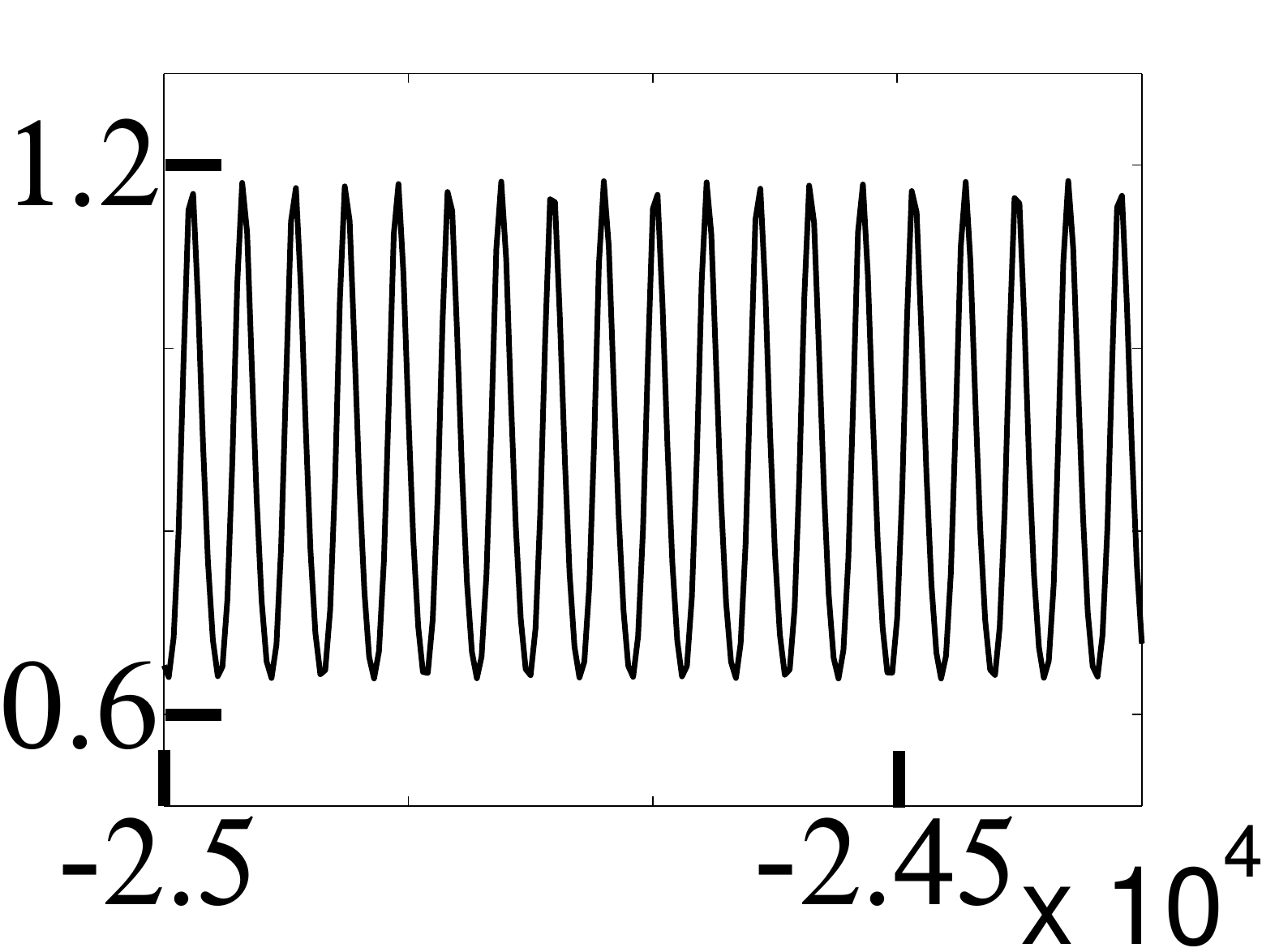}}
\put(75,60){\includegraphics[width=0.1\textwidth,height=0.07\textheight]
{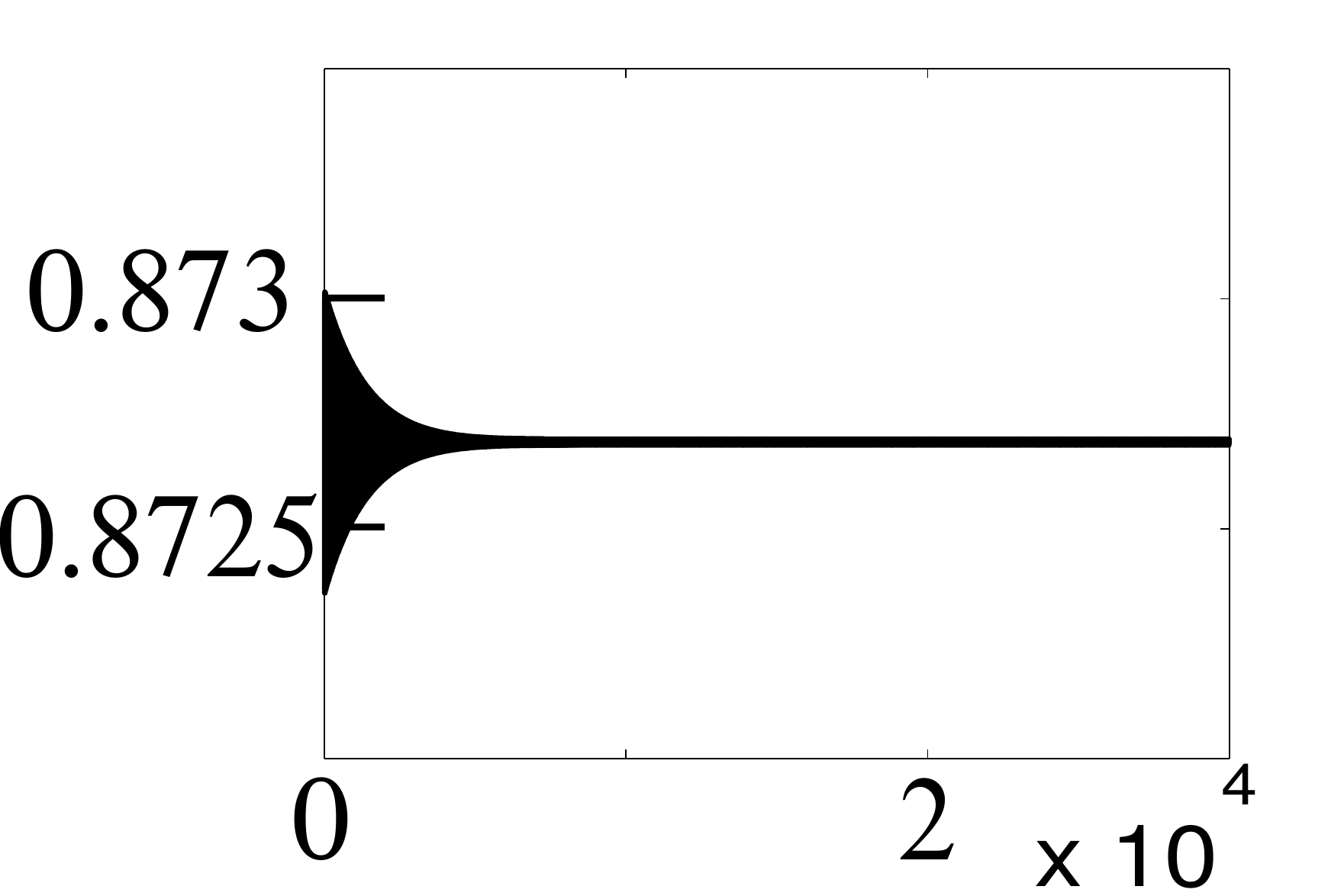}}
\put(48,82){(1a)}
\put(48,56){(1b)}
\put(86,82){(1c)}
\put(81,5){$\ast$}
\put(75,20){\includegraphics[width=0.1\textwidth,height=0.07\textheight]
{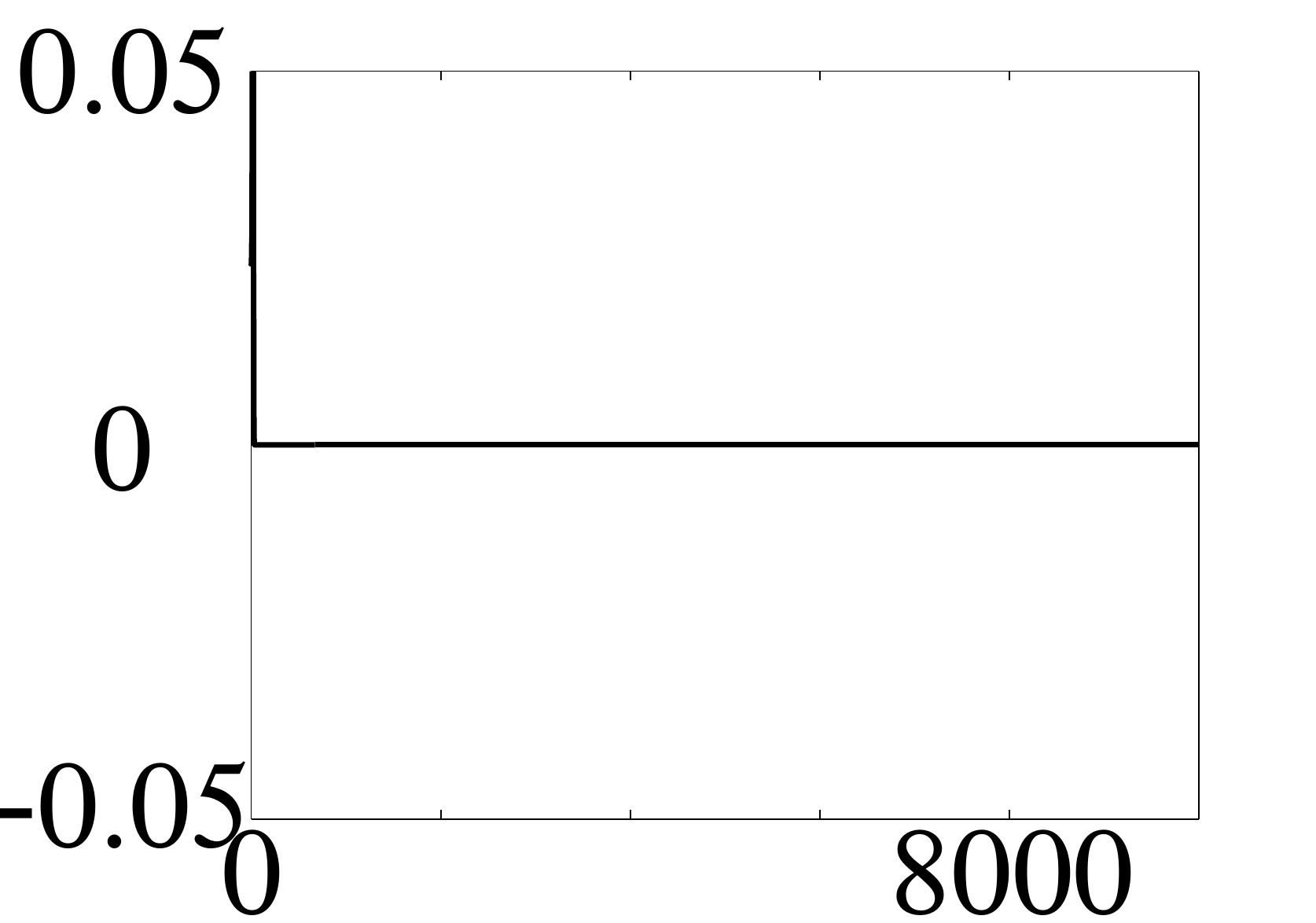}}
\put(86,42){(1d)}
\put(63,83){\vector(2,1){15}}
\put(81,83){\vector(0,1){6}}
\put(82,20){\vector(0,-1){8}}
\put(45,56){\vector(0,1){4}}
\end{overpic}
}
\hspace{0.5in}
\vspace{0.4in}
\subfloat{
\begin{overpic}[width=0.4\textwidth]{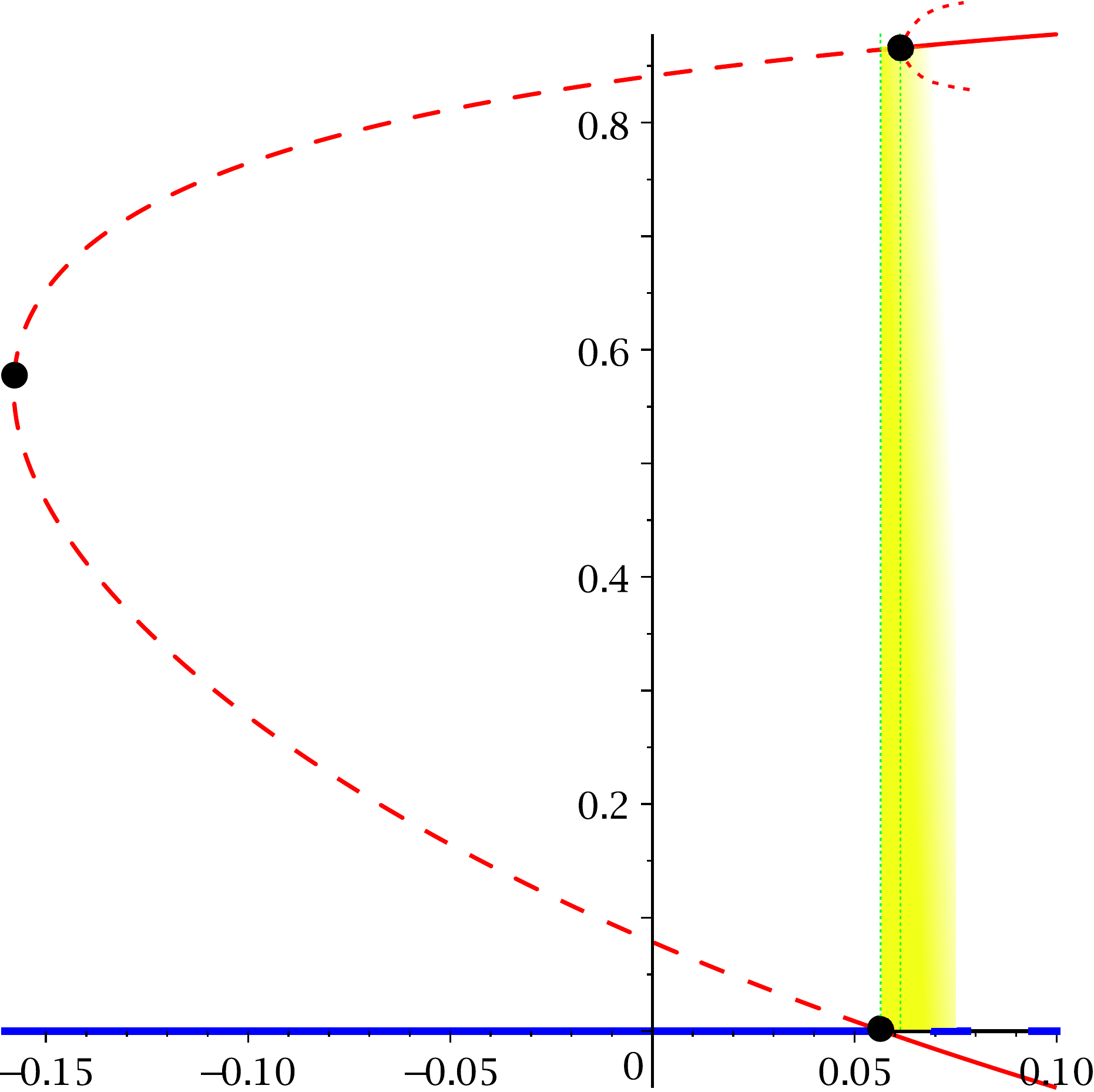}
\put(50,100){(2)}
\put(50,-2){$B$}
\put(50,55){$Y$}
\put(5,65){ Turning}
\put(80,89){ Hopf$_{\text{sub}}$}
\put(73,10){ Transcritical}
\put(85,67){$\ast$}
\put(85,47){$\ast$}
\put(81,57){$\ast$}
\put(91,80){\vector(-1,-2){5}}
\put(90,58){\vector(-1,0){6}}
\put(90,36){\vector(-1,2){5}}
\put(95,85){(2a)}
\put(95,65){(2b)}
\put(95,44){(2c)}
\put(90,70){\includegraphics[width=0.07\textwidth]
{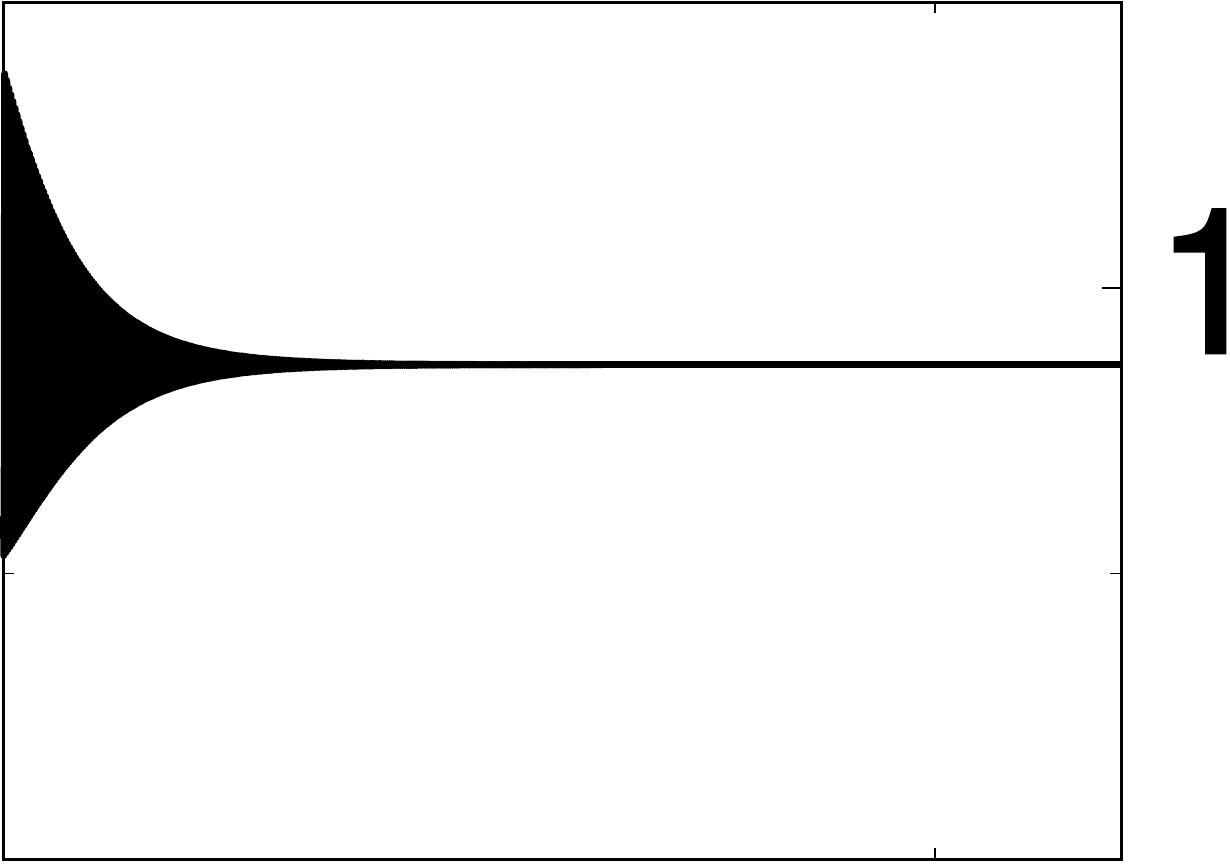}}
\put(90,50){\includegraphics[width=0.07\textwidth]
{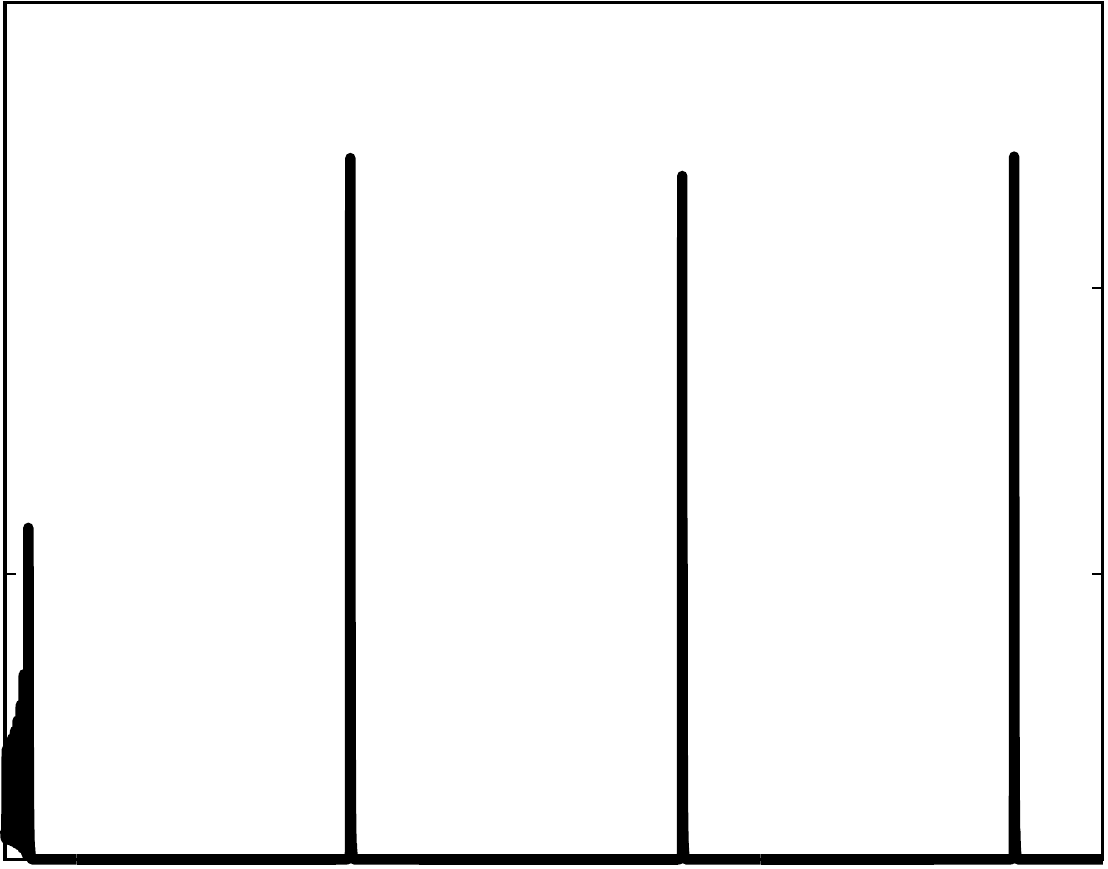}}
\put(90,30){\includegraphics[width=0.07\textwidth]
{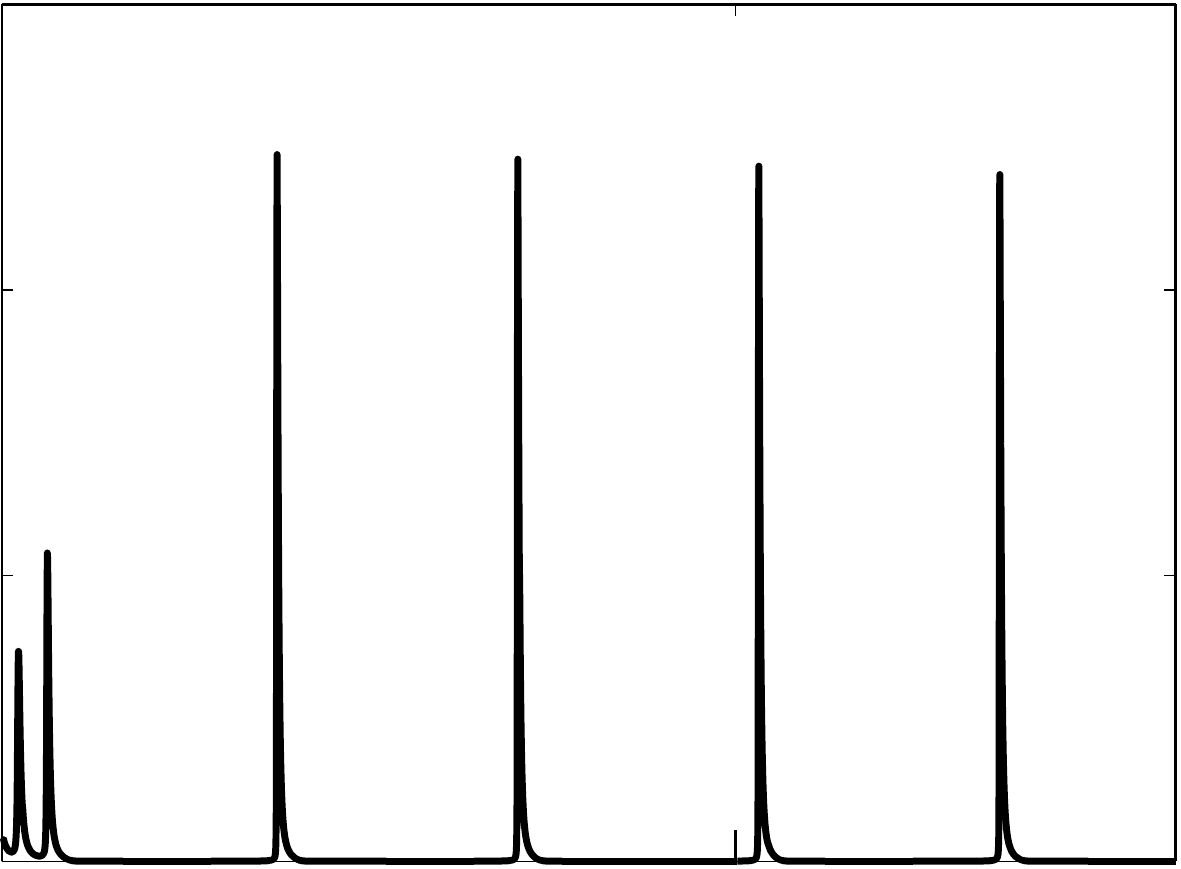}}
\end{overpic} 
}
\\ \vspace{-0.2cm}
\hspace{-.0in}
\subfloat{
\begin{overpic}[width=0.4\textwidth]{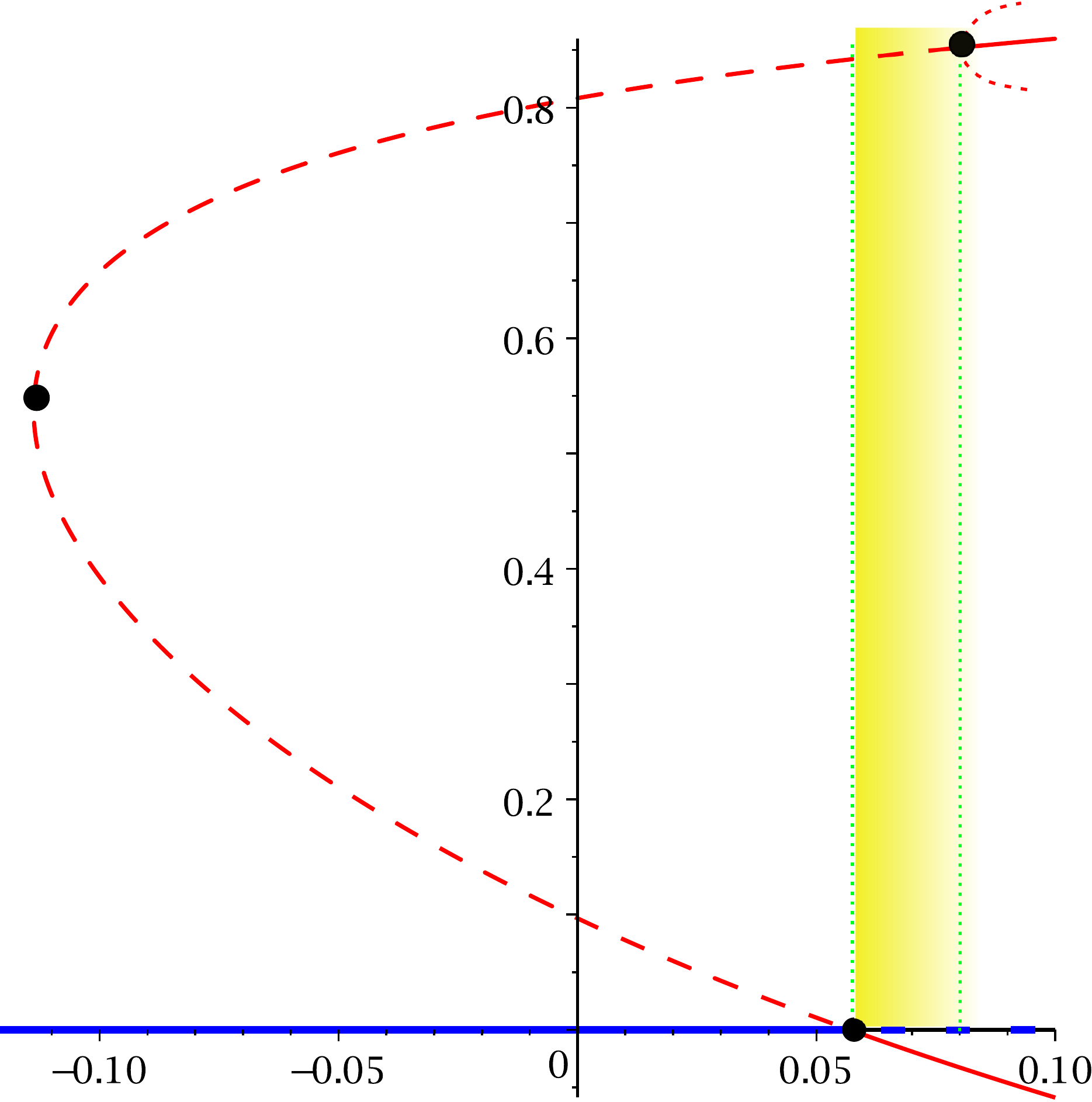}
\put(50,100){(3)}
\put(50,0){$B$}
\put(47,55){$Y$}
\put(-7,65){ Turning}
\put(70,97){ Hopf$_{\text{sub}}$}
\put(73,10){ Transcritical}
\put(83,93){$\ast$}
\put(65,85){\vector(2,1){15}}
\put(83,16){$\ast$}
\put(95,85){\vector(-1,1){5}}
\put(88,90){$\ast$}
\put(75,25){\vector(1,-1){7}}
\put(88,16){$\ast$}
\put(100,20){\vector(-2,-1){7}}
\put(55,70){\includegraphics[width=0.07\textwidth]
{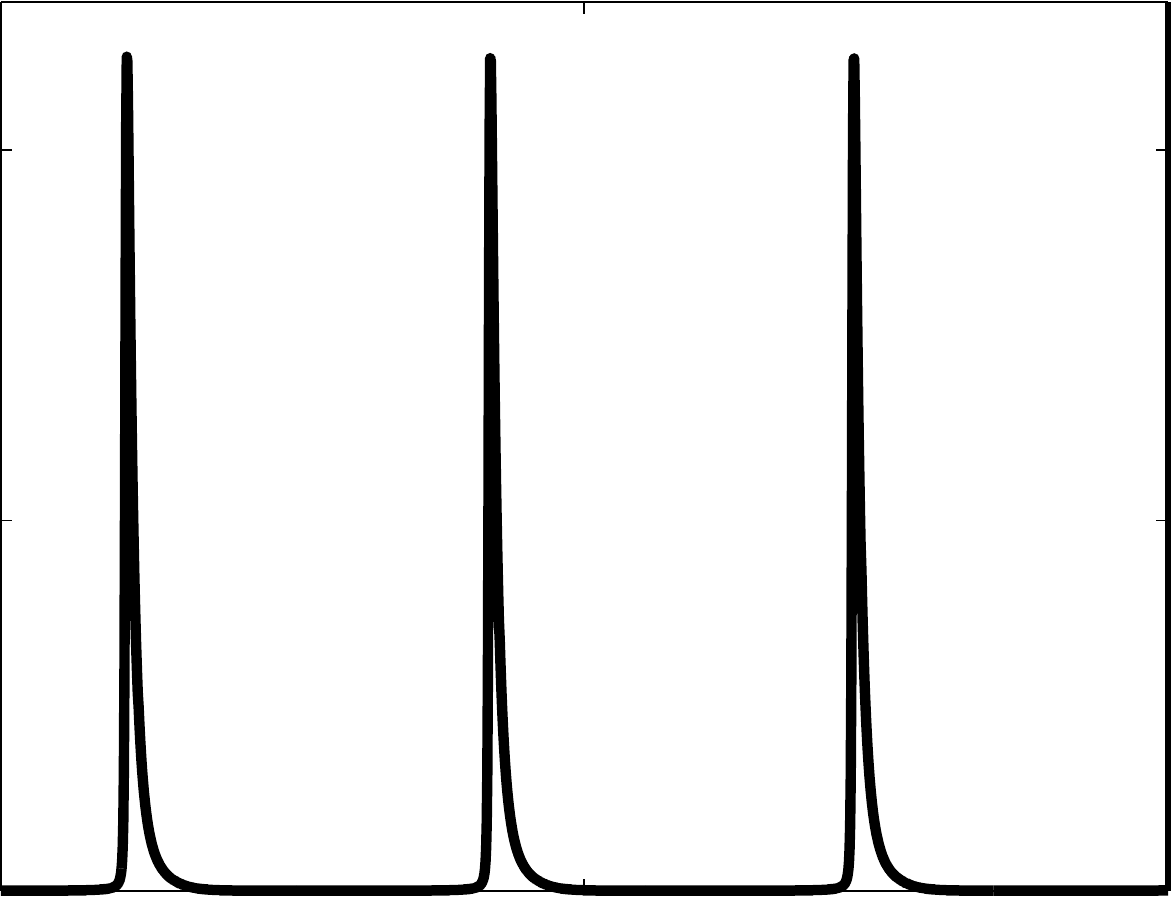}}
\put(60,85){(3a)}
\put(55,20){\includegraphics[width=0.07\textwidth]
{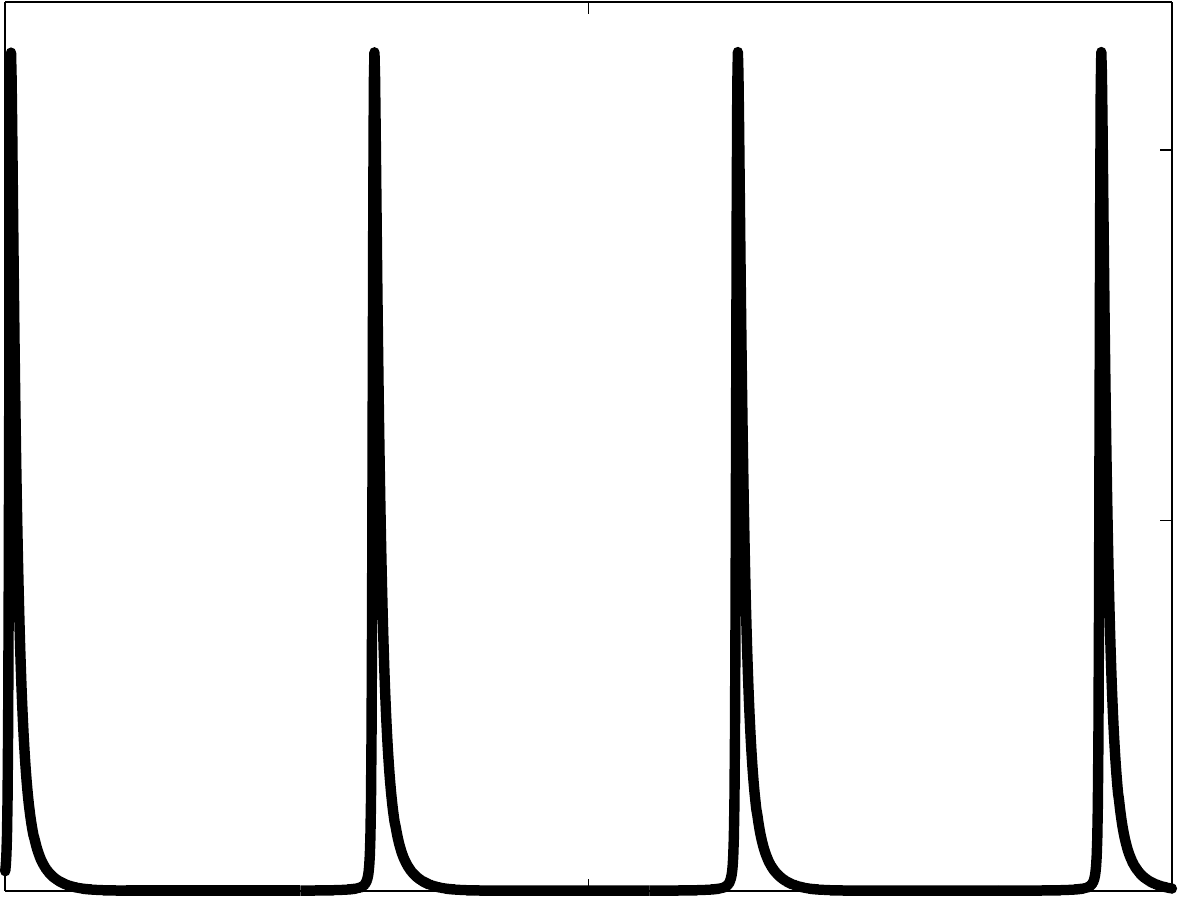}}
\put(60,35){(3b)}
\put(90,70){\includegraphics[width=0.08\textwidth]
{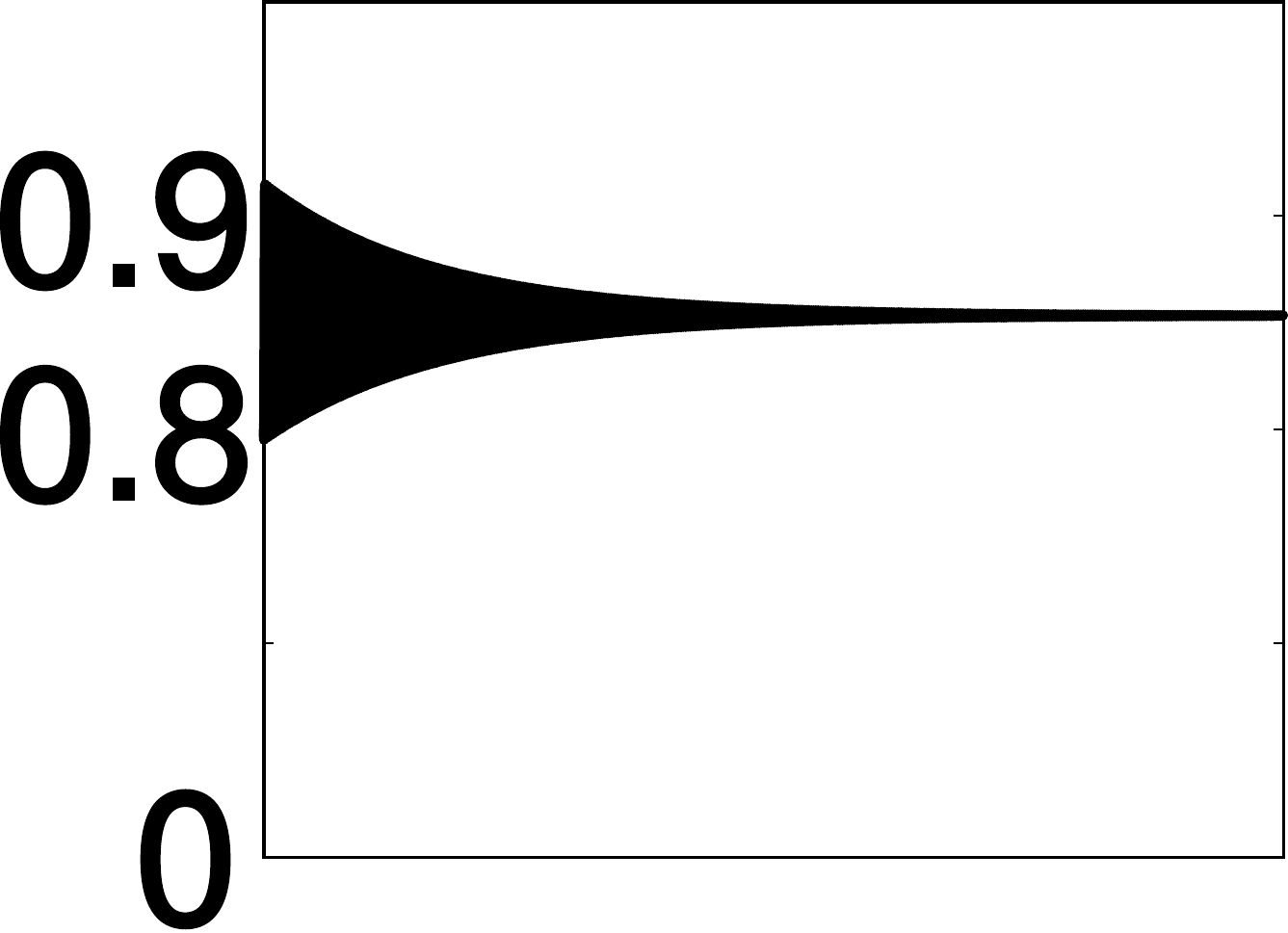}}
\put(95,86){(3c)}
\put(90,20){\includegraphics[width=0.07\textwidth]
{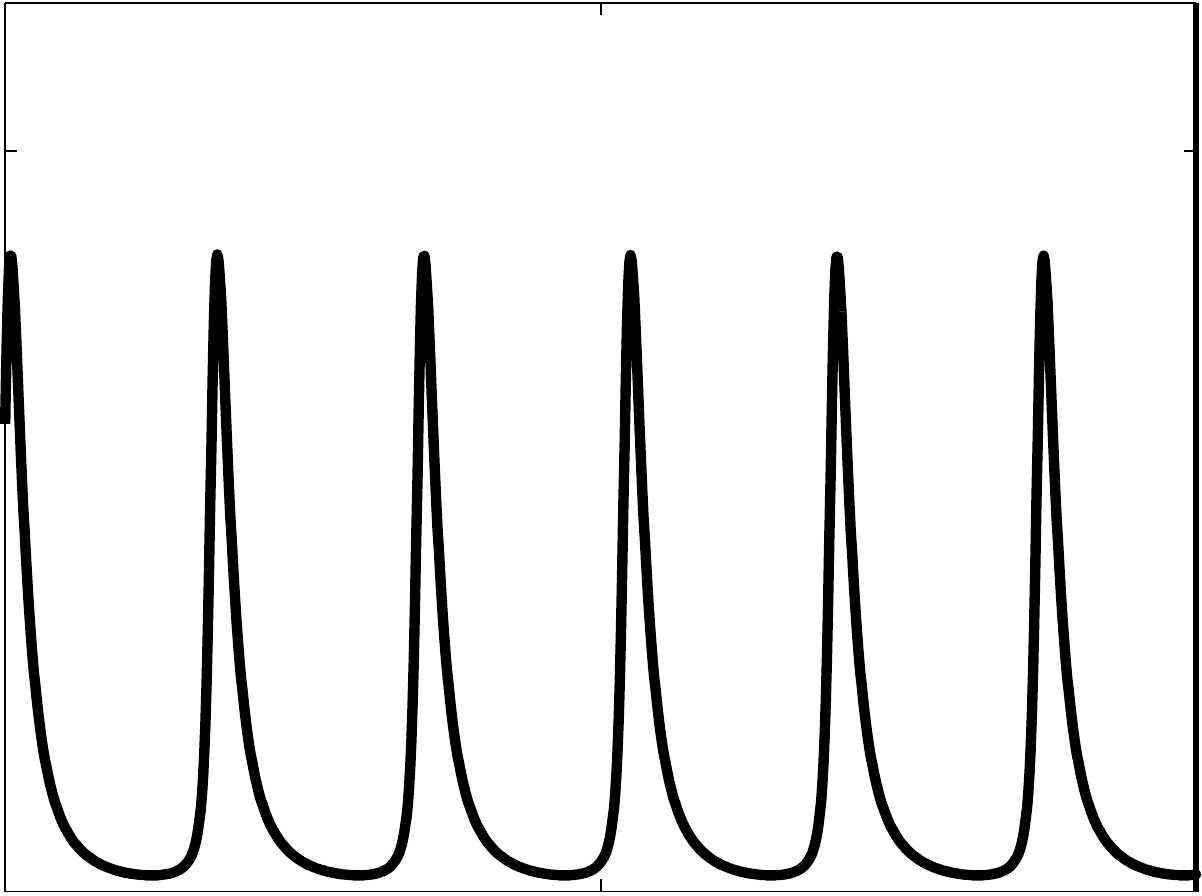}}
\put(95,35){(3d)}
\end{overpic}
}
\hspace{0.7in}
\subfloat{
\begin{overpic}[width=0.4\textwidth]{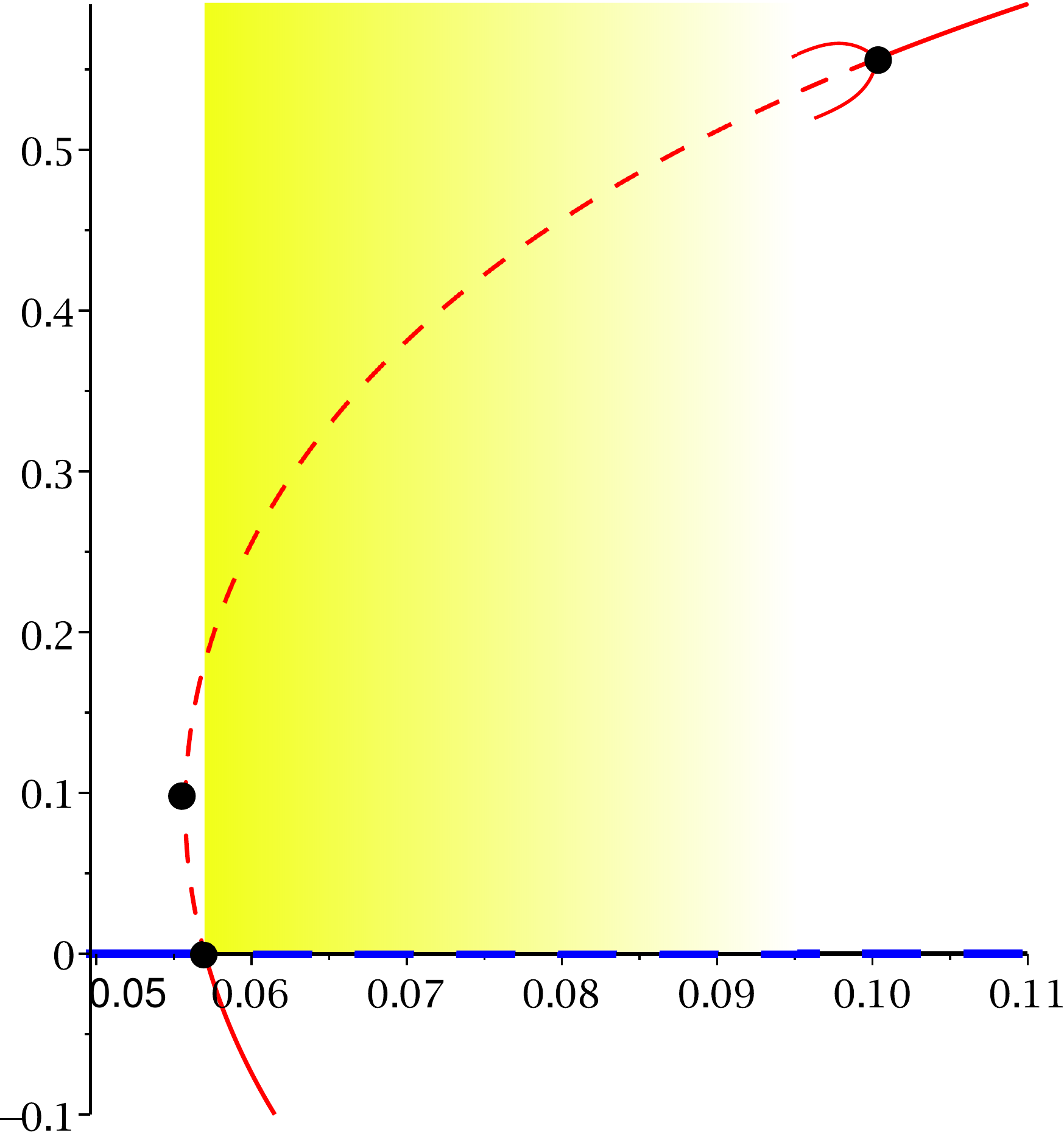}
\put(50,100){(4)}
\put(50,0){$B$}
\put(-2,50){$Y$}
\put(7,31){ Turning}
\put(75,91){ Hopf$_{\text{super}}$}
\put(7,18){ Transcritical}
\put(22,30){$\ast$}
\put(49,30){$\ast$}
\put(25,35){\includegraphics[width=0.07\textwidth]
{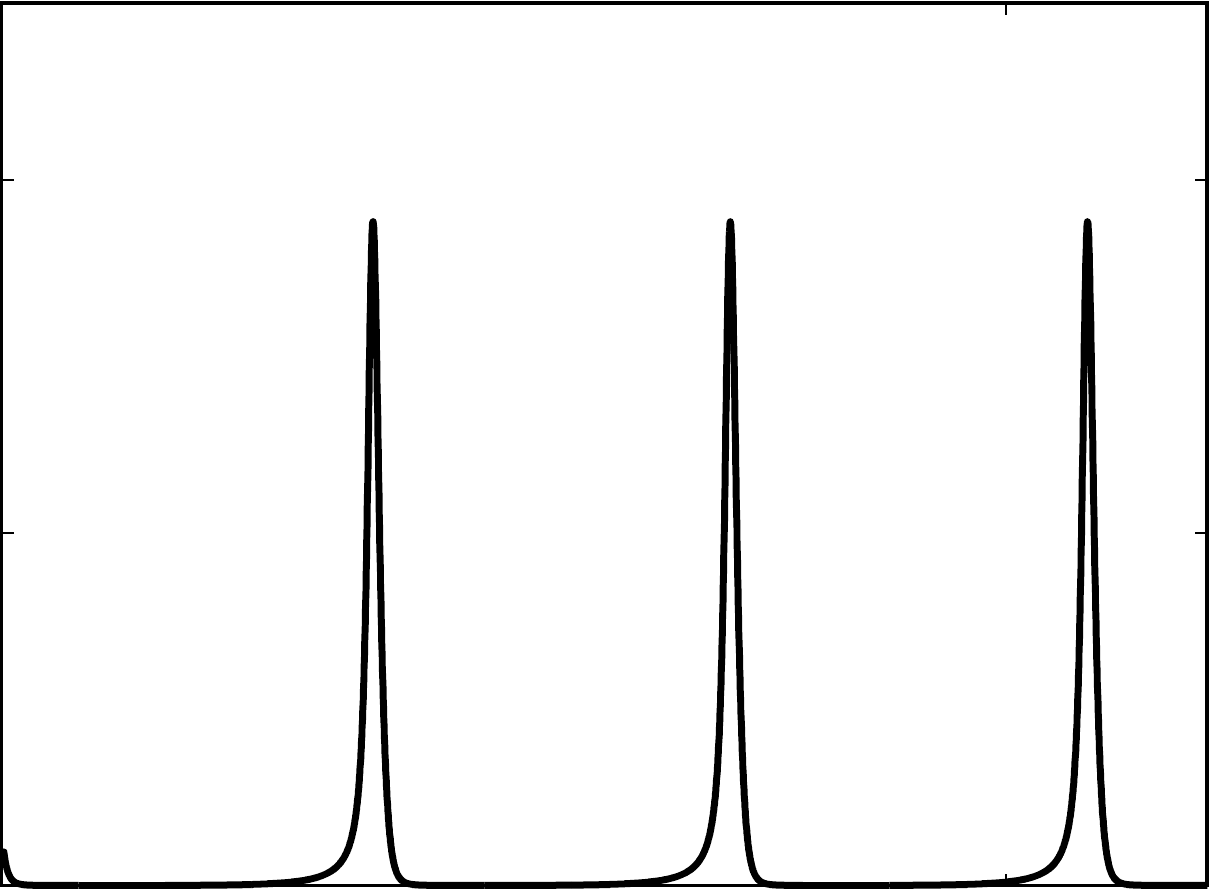}}
\put(25,35){\vector(-1,-1){3}}
\put(55,35){\includegraphics[width=0.07\textwidth]
{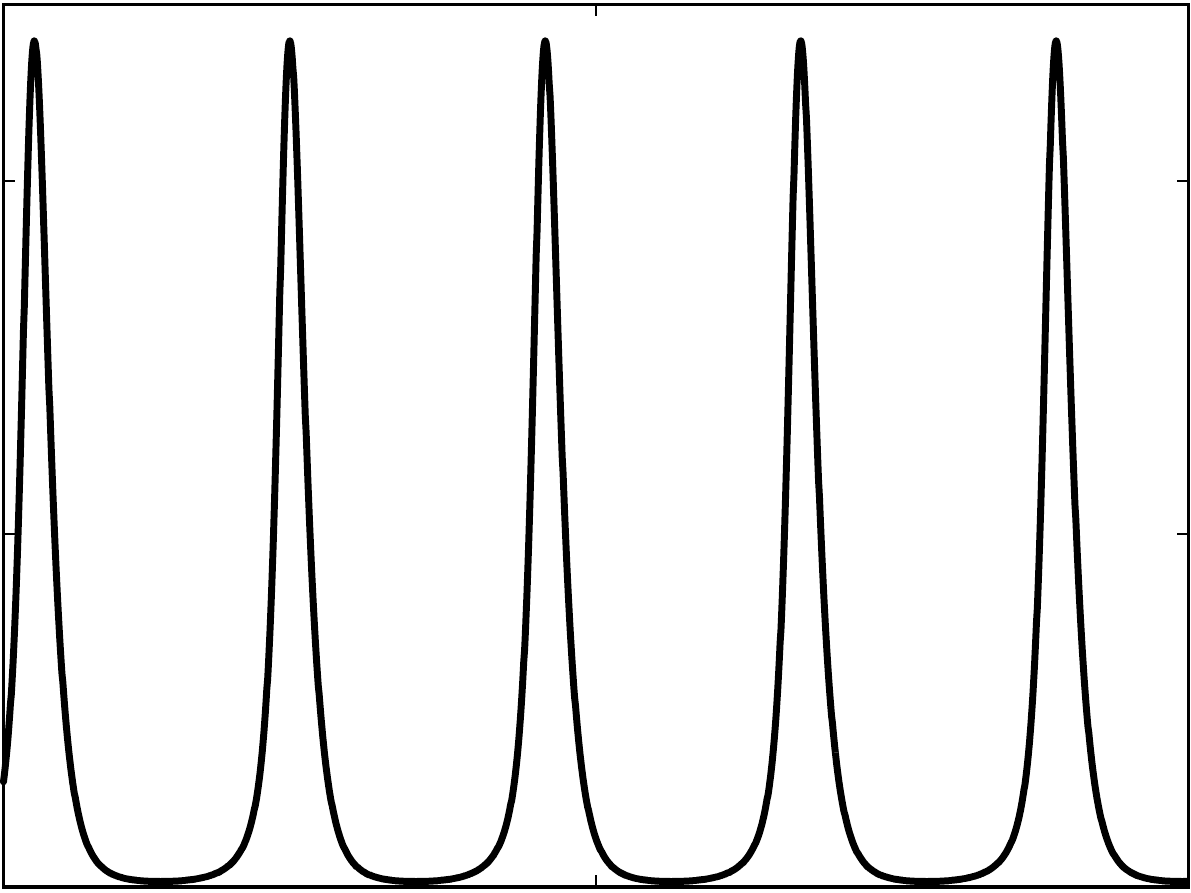}}
\put(55,35){\vector(-1,-1){3}}
\put(32,48){(4a)}
\put(60,48){(4b)}
\end{overpic}
}
}
\caption{Dynamical behaviors of system~(\ref{Paper3_Eq4}) corresponding to eight
cases listed in Table~\ref{Paper3_Table2} and \ref{Paper3_Table4}.
All insets are simulated time histories of $Y$ vs. $t$.
The yellow areas fading to white show regions in which recurrent
behavior occurs and fades to regular oscillations.}
\end{figure}

\begin{figure} 
\vspace{0.15in}  
\hspace{0.10in} 
\ContinuedFloat
{\tiny
\hspace{-0.1in} 
\subfloat{
\begin{overpic}[width=0.4\textwidth]{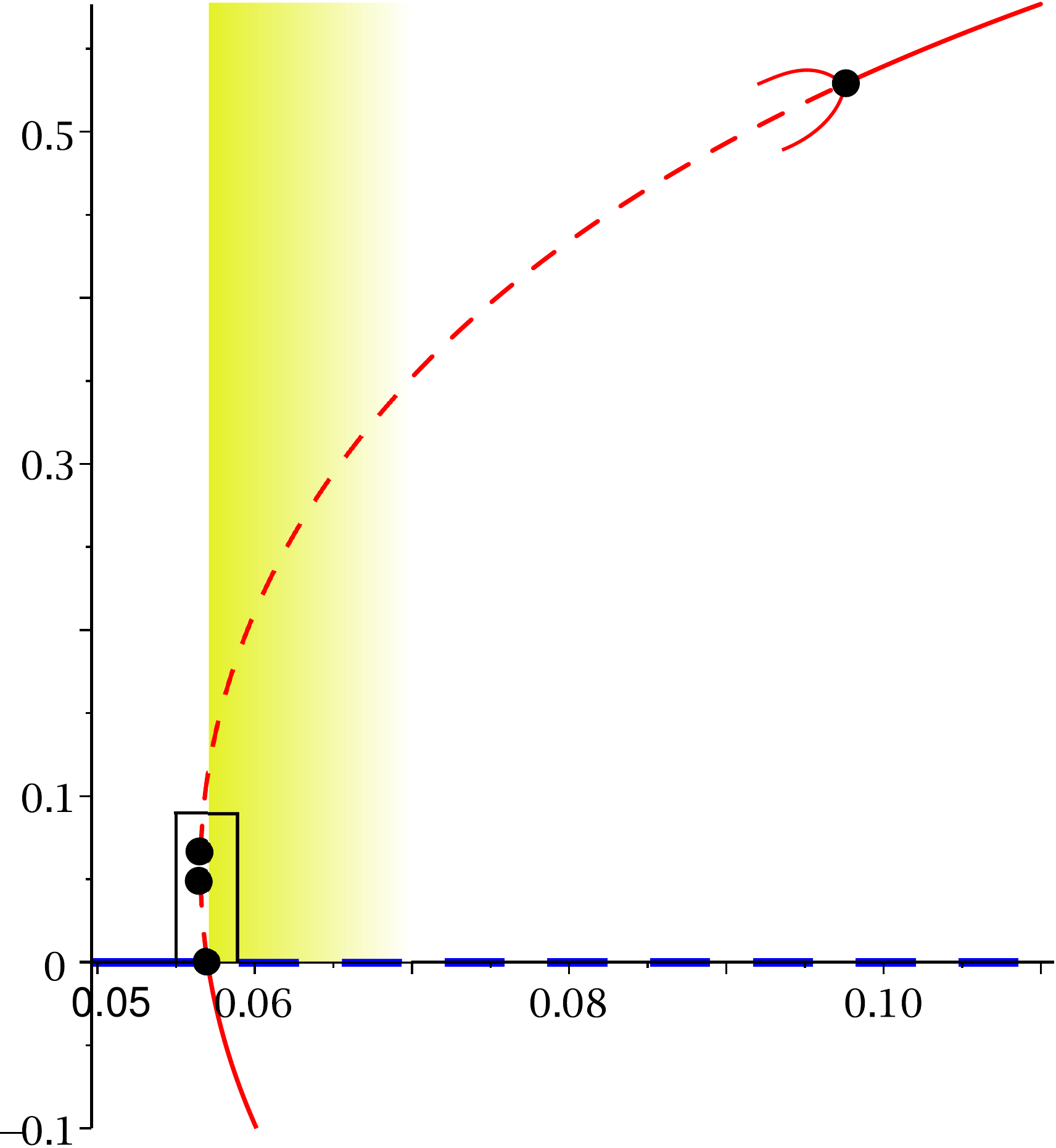}
\put(50,100){(5)}
\put(50,0){$B$}
\put(-2,50){$Y$}
\put(80,25){ Transcritical}
\put(60,60){ Turning}
\put(60,63){ Hopf$_{\text{super1}}$}
\put(76,91){ Hopf$_{\text{super2}}$}
\put(50,17){\includegraphics[width=0.2\textwidth,height=0.2\textheight]
{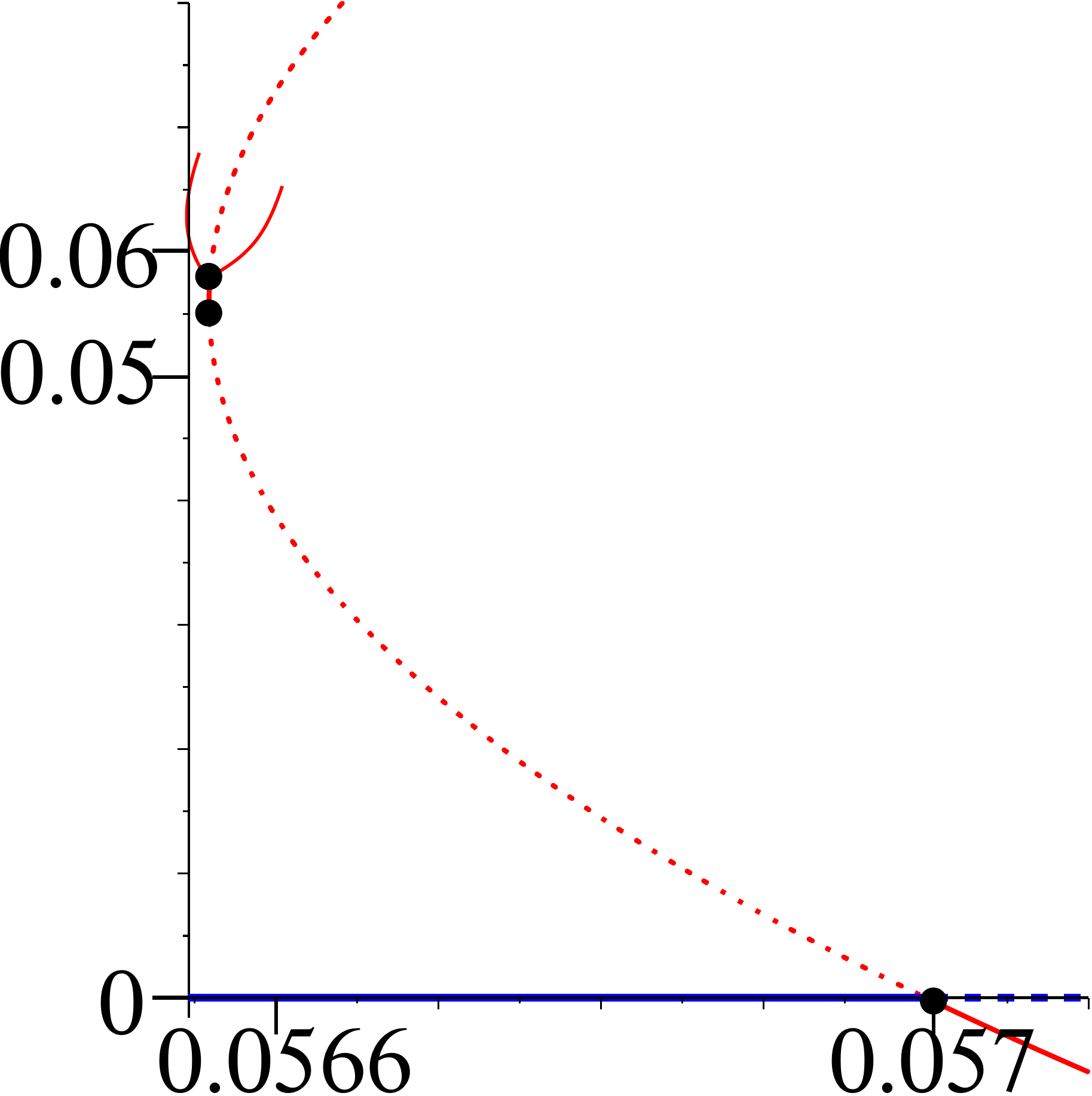}}
\put(21,27){\vector(1,0){30}}
\put(30,35){\includegraphics[width=0.1\textwidth,height=0.05\textheight]
{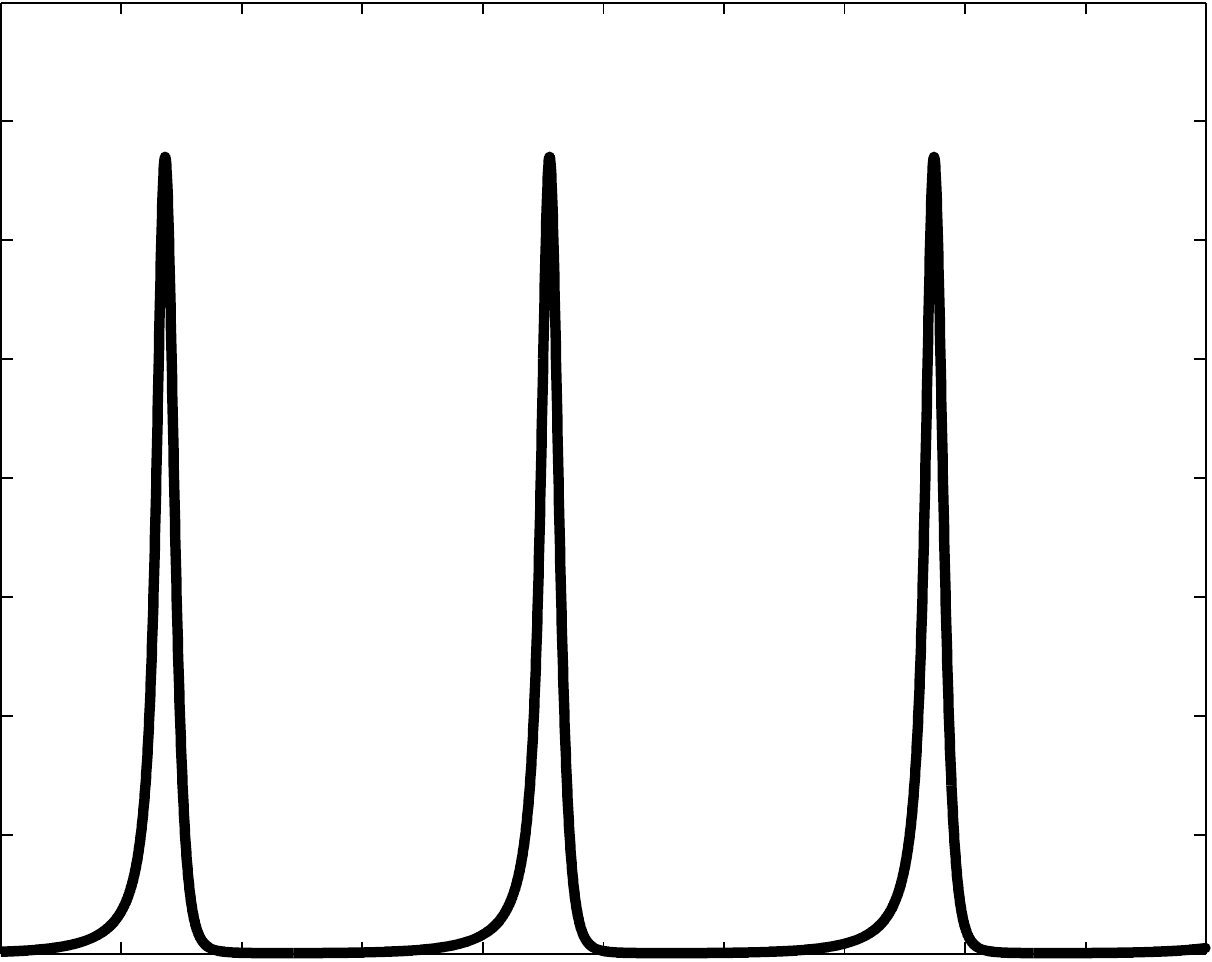}}
\put(30,35){\vector(-2,-1){7}}
\put(21,30){$\ast$}
\end{overpic} 
}
\hspace{0.6in}
\subfloat{
\begin{overpic}[width=0.4\textwidth]{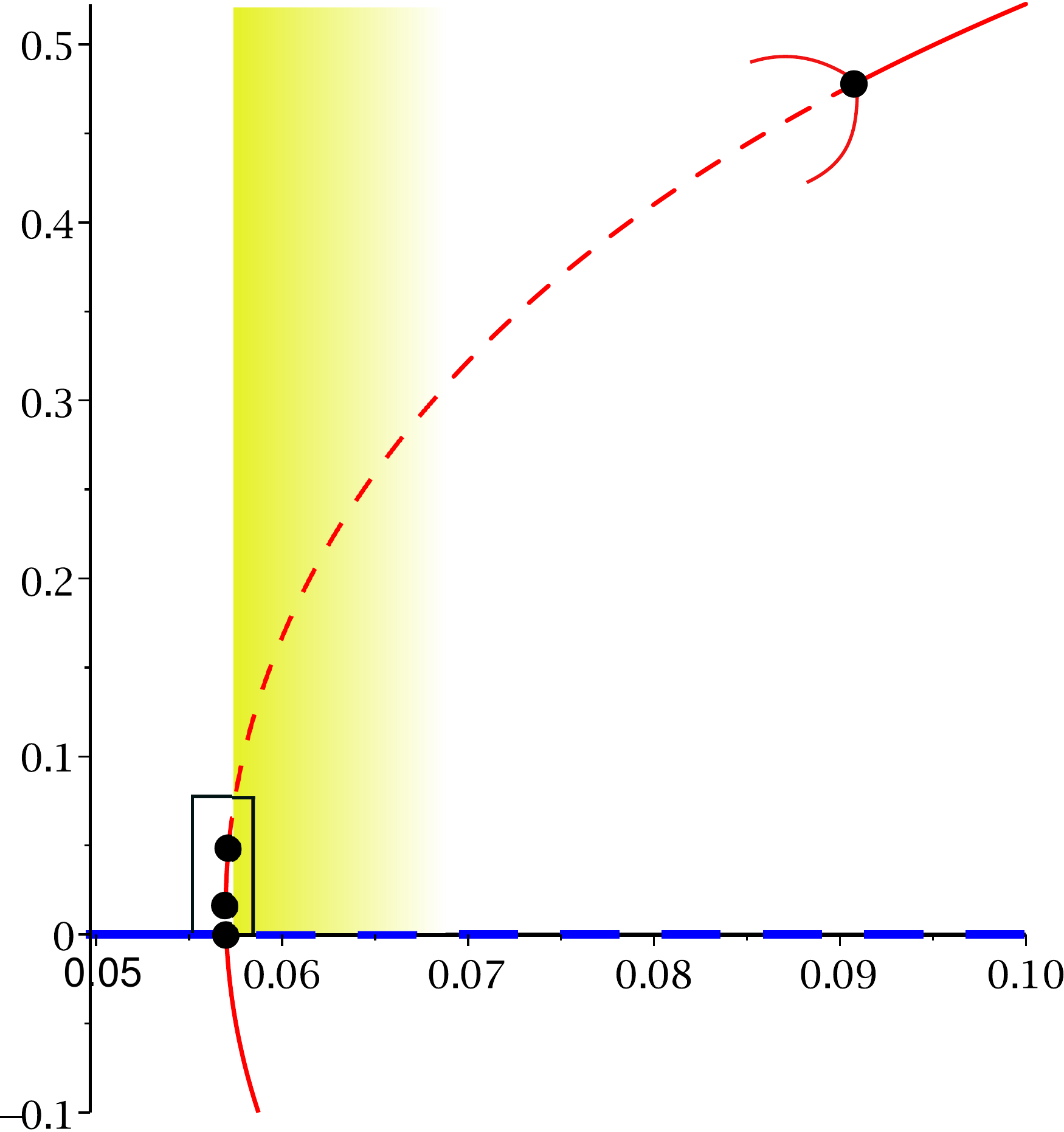}
\put(50,100){(6)}
\put(50,0){$B$}
\put(-2,50){$Y$}
\put(65,27){ Transcritical}
\put(90,65){ Hopf$_{\text{super1}}$}
\put(80,91){ Hopf$_{\text{super2}}$}
\put(65,35){ Turning}
\put(50,20){\includegraphics[width=0.2\textwidth,height=0.17\textheight]
{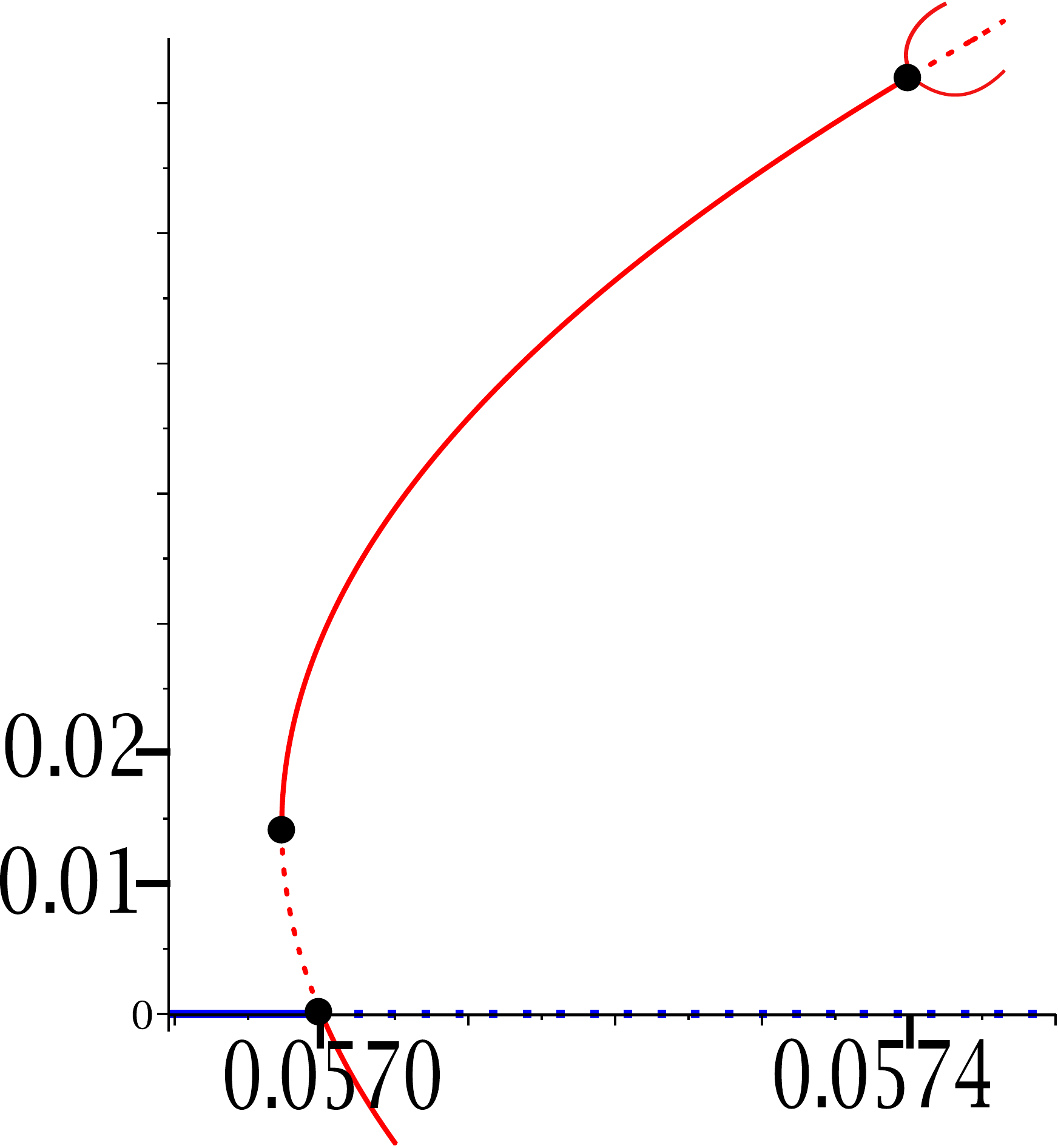}}
\put(10,60){\includegraphics[width=0.1\textwidth,height=0.05\textheight]
{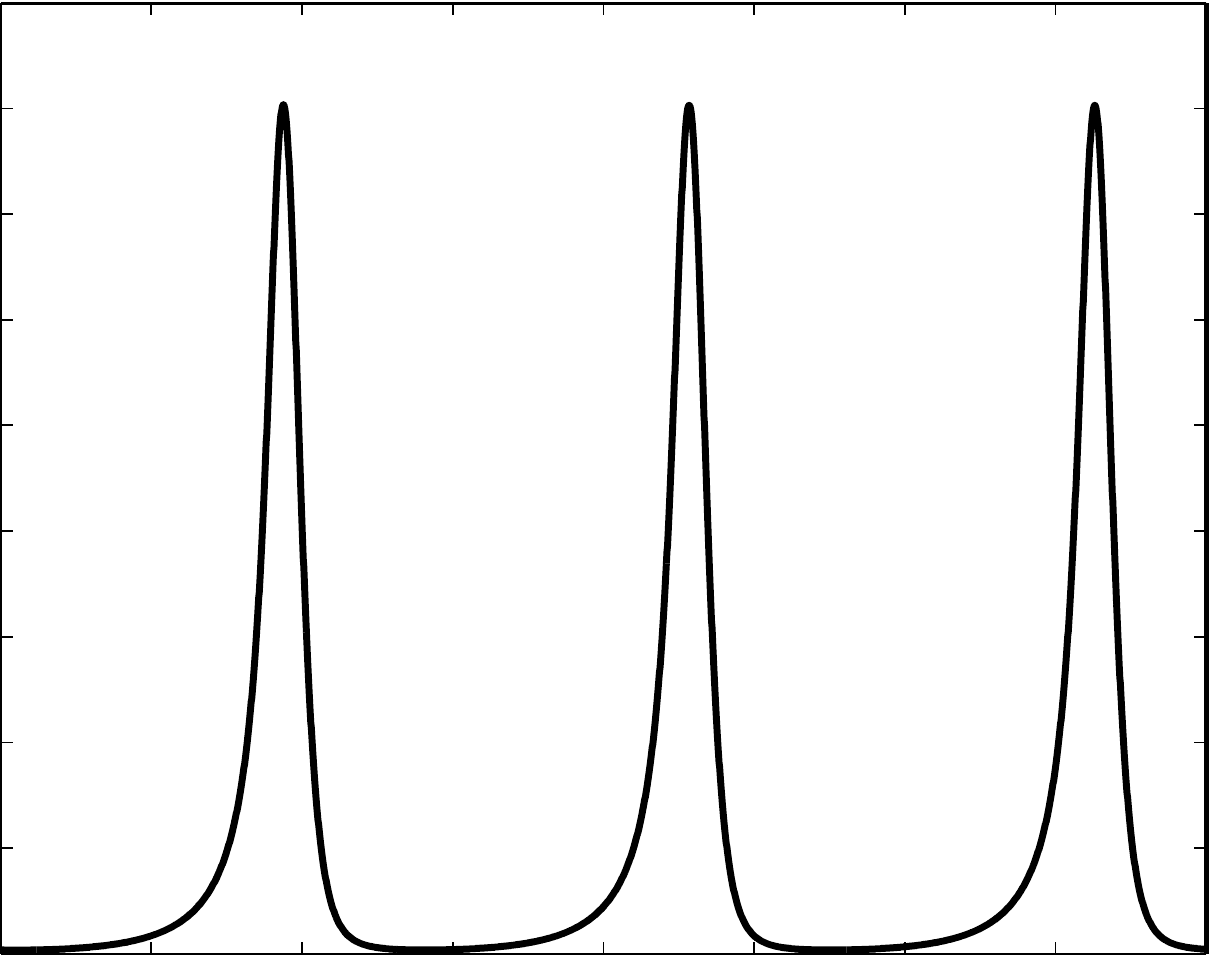}}
\put(24,33){$\ast$}
\put(25,60){\vector(0,-1){25}}
\end{overpic} 
}
 
\vspace{0.8cm}
\hspace{-0.0in}
\subfloat{
\begin{overpic}[width=0.4\textwidth]{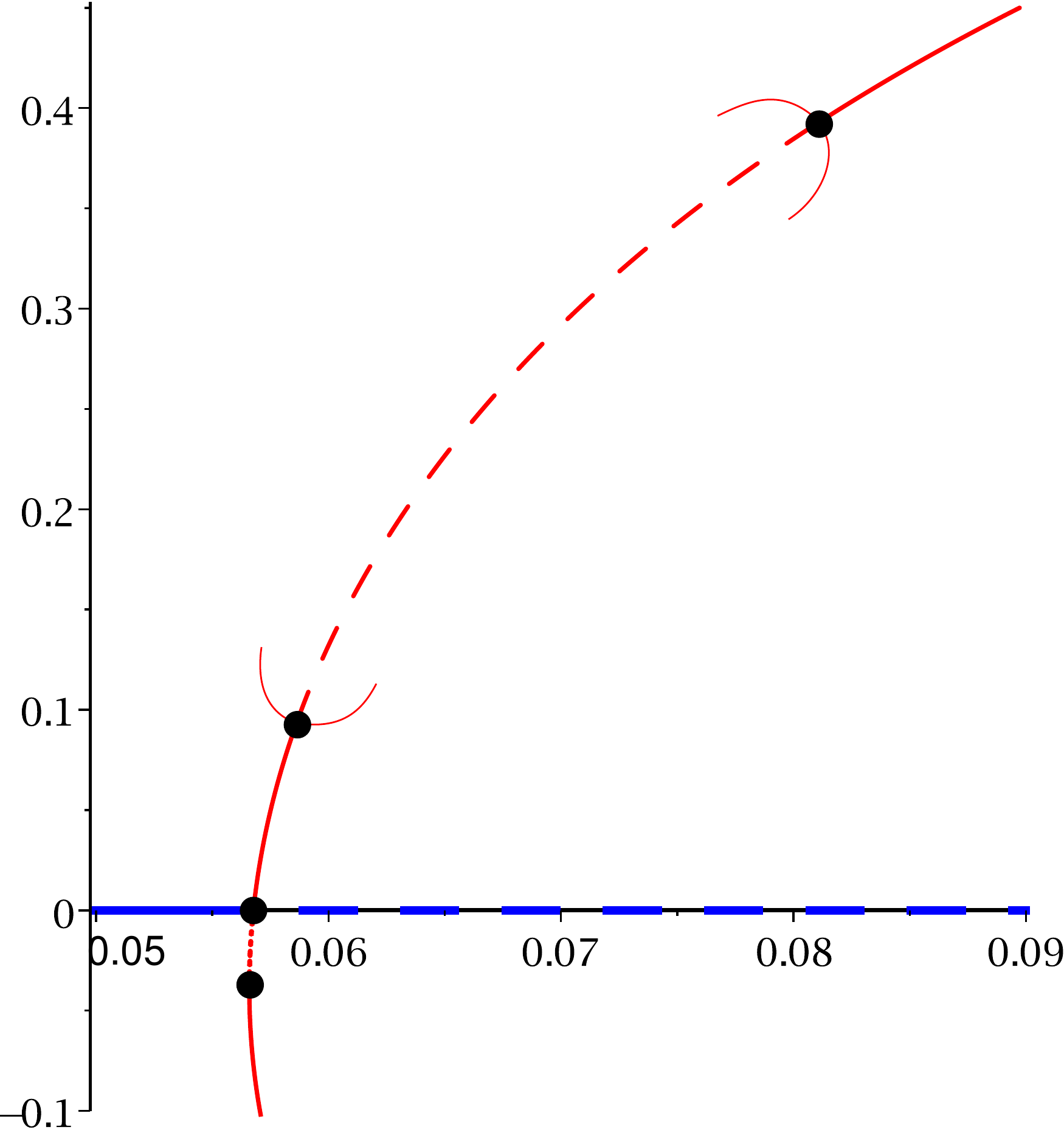}
\put(50,106){(7)}
\put(50,10){$B$}
\put(-2,50){$Y$}
\put(24,21){ Transcritical}
\put(7,35){ Hopf$_{\text{super1}}$}
\put(75,89){ Hopf$_{\text{super2}}$}
\put(22,10){ Turning}
\put(50,30){\includegraphics[width=0.15\textwidth,height=0.07\textheight]
{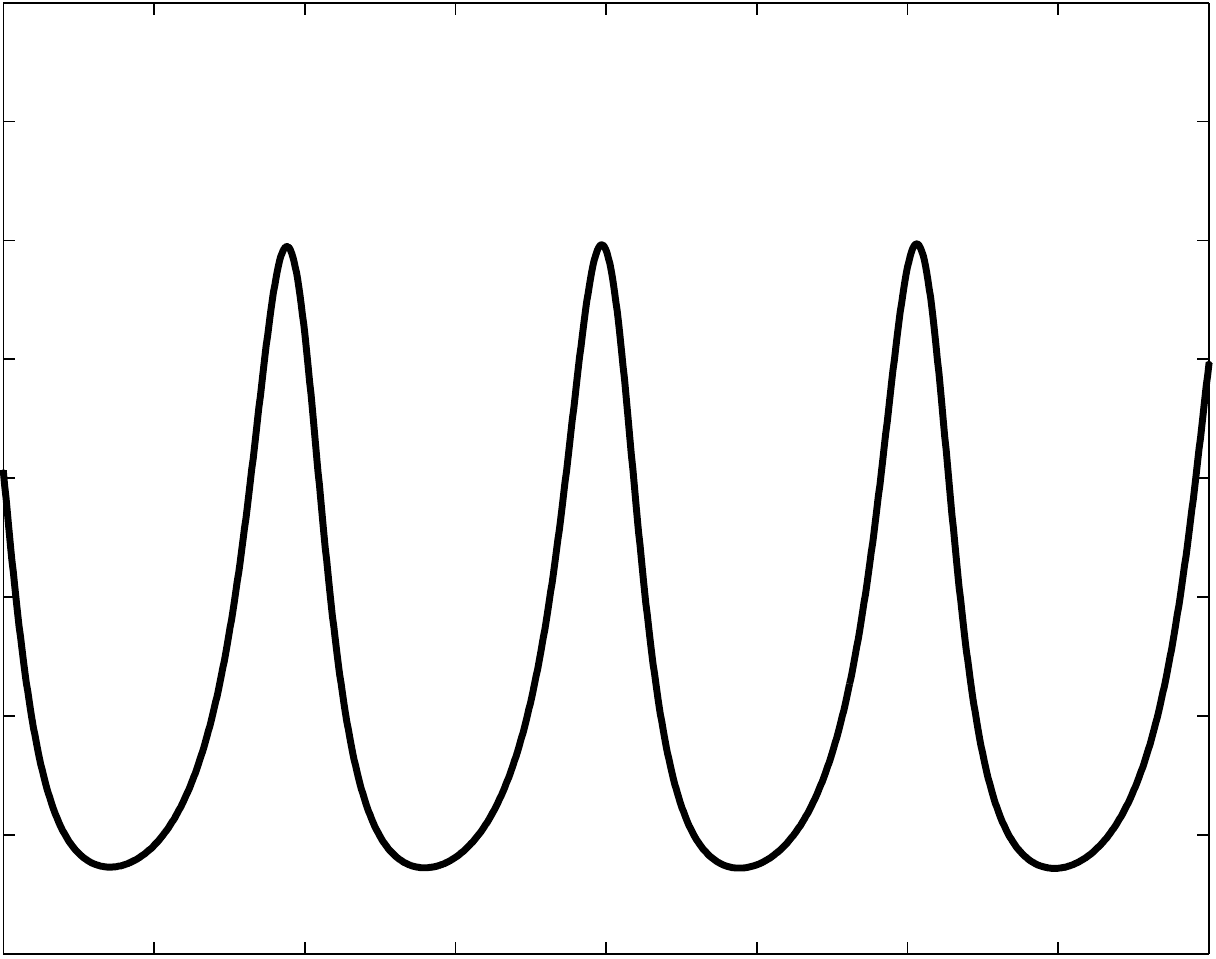}}
\put(30,37){$\ast$}
\put(50,38){\vector(-1,0){15}}
\end{overpic} 
}
\subfloat{
\hspace{0.6in}
\begin{overpic}[width=0.4\textwidth]{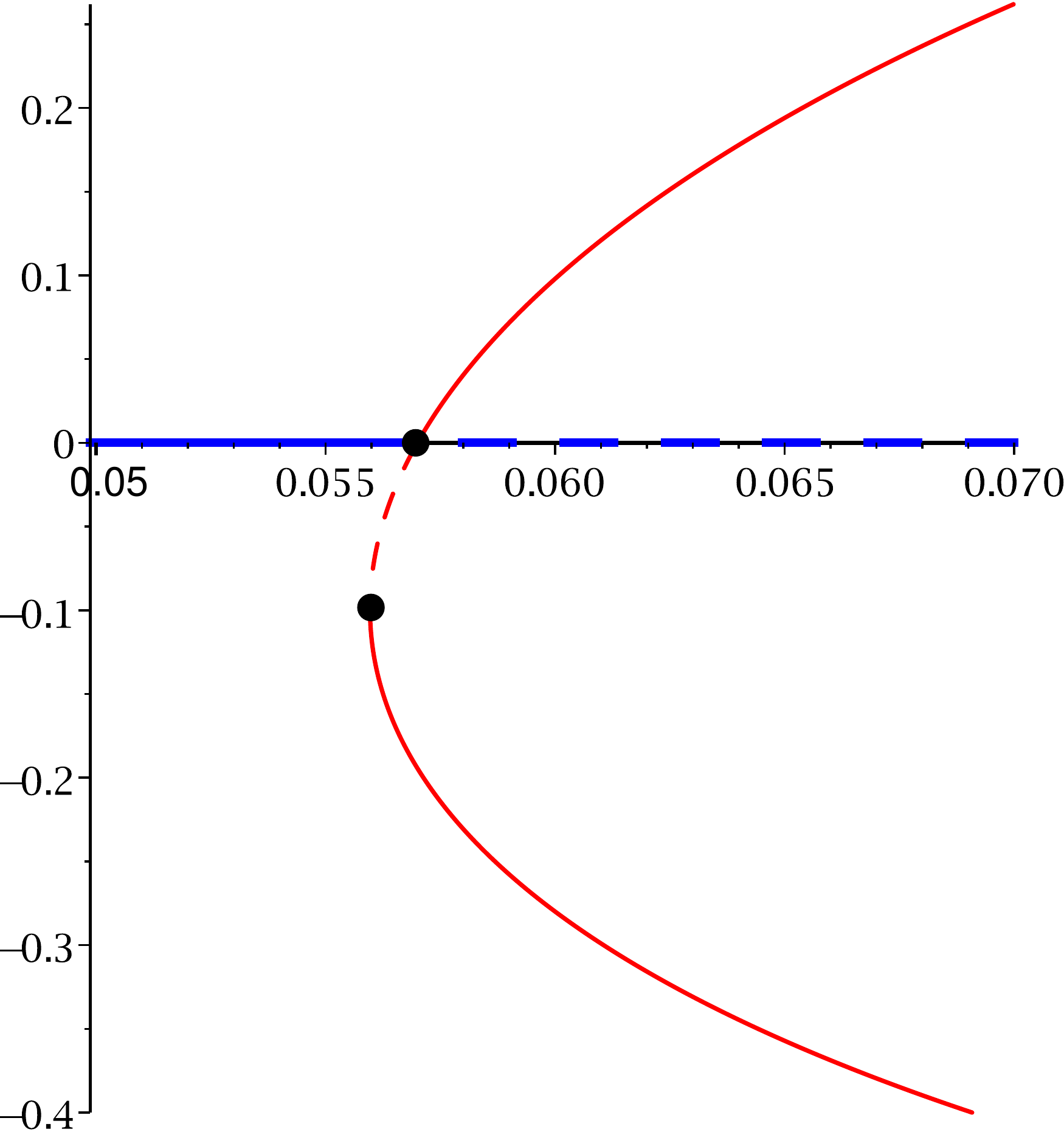}
\put(50,106){(8)}
\put(50,50){$B$}
\put(-2,50){$Y$}
\put(40,62){ Transcritical}
\put(35,45){ Turning}
\put(10,75){\includegraphics[width=0.15\textwidth,height=0.07\textheight]
{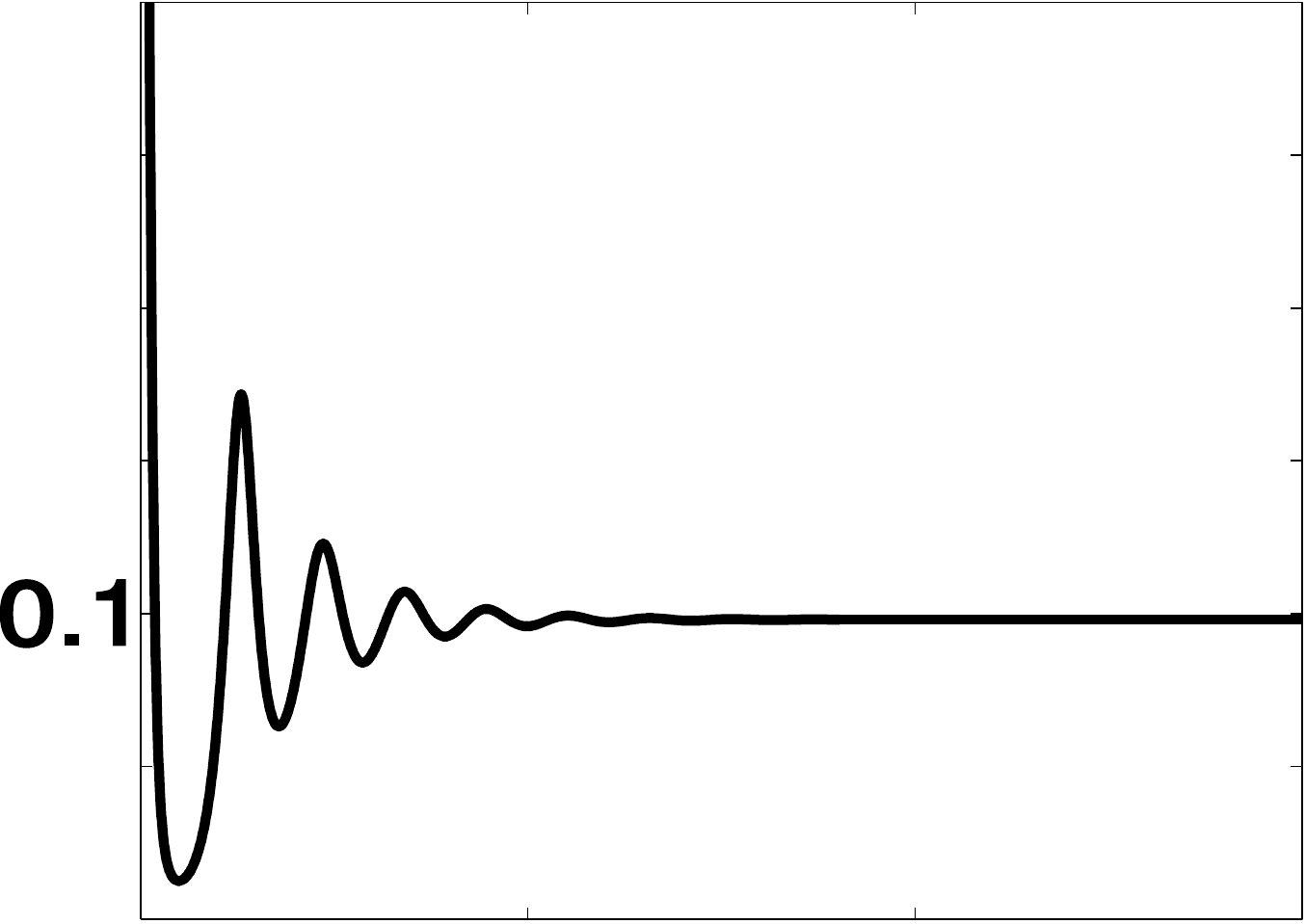}}
\put(48,68){$\ast$}
\put(40,74){\vector(2,-1){8}}
\end{overpic}
}
}
\caption{Dynamical behaviors of system~(\ref{Paper3_Eq4}) corresponding to eight
cases listed in Table~\ref{Paper3_Table2} and \ref{Paper3_Table4}.
All insets are simulated time histories of $Y$ vs. $t$.
The yellow areas fading to white show regions in which recurrent
behavior occurs and fades to regular oscillations.}
\label{Paper3_Fig4}
\end{figure}

\section{Negative backward bifurcation in an autoimmune disease model}

In the previous section, 
we examined three cases of negative backward bifurcation:
Table~\ref{Paper3_Table1} Case 4 for system~(\ref{Paper3_Eq7}) and
Table~\ref{Paper3_Table2} Case (7) and (8) for system~(\ref{Paper3_Eq4}).
The analytical and numerical results showed that solutions typically
converge to the infected equilibrium in these cases, 
and the parameter range for Hopf bifurcation is very limited.
As a result, negative backward bifurcation tends to give no interesting 
behavior.
In this section, however, we shall explore an established 
autoimmune model~\cite{AW2011} in which negative backward bifurcation occurs. 
We demonstrate that after modification, 
the autoimmune model can also exhibit recurrence.

The autoimmune model~\cite{AW2011} takes the form
\begin{equation}
\begin{array}{ll}
\dss\frac{{\mathrm d} A}{{\mathrm d} t} = f\tilde{v} G 
- (\sigma_1 R_n + b_1) A - \mu_A A \\[1.0ex]
\dss\frac{{\mathrm d} R_n}{{\mathrm d} t}= 
(\pi_1 E + \beta)A - \mu_n R_n \\[1.0ex]
\dss\frac{{\mathrm d} E}{{\mathrm d} t} = \lambda_E A - \mu_E E \\[1.0ex]
\dss\frac{{\mathrm d} G}{{\mathrm d} t} = \gamma E - \tilde{v} G - \mu_G G,
\end{array}
\label{Paper3_Eq16}
\end{equation}
where mature pAPCs ($A$) undergo maturation by intaking self-antigen ($G$),
at rate $f\tilde{v}$, and are suppressed by specific regulatory T cells,
T$_\text{Reg}$ cells ($R_n$), at rate $\sigma_1$;
$b_1$ represents additional non-specific background suppression.
The T$_\text{Reg}$ cells are activated by mature pAPCs 
at a rate proportional to the number of auto-reactive effector T cells ($E$) at
rate $\pi_1$, and by other sources at rate $\beta$.
Active auto-reactive effector T cells ($E$) come from the activation process
initiated by mature pAPCs, at rate $\lambda_E$, then attack
healthy body tissue and release free self-antigen ($G$) at rate $\gamma$,
which is ready for mature pAPCs to engulf; the antigen engulfing rate is
$\tilde{v}$. The death rates of the populations $A$, $R_n$, $E$, and $G$
are denoted by $\mu_A$, $\mu_n$, $\mu_E$, and $\mu_G$, respectively.

Following the steps described by Zhang et al. (submitted for 
publication), 
system~(\ref{Paper3_Eq16}), can be reduced via quasi-steady
state analysis to a 2-dimensional system: 
\begin{equation}
\begin{array}{ll}
\dss\frac{{\mathrm d} A}{{\mathrm d} t} = 
[\frac{f\tilde{v}\gamma\lambda_E}{\mu_E(\tilde{v}+\mu_G)}-b_1-\mu_A]A -
\sigma_1 R_n A, \\[1.5ex] 
\dss\frac{{\mathrm d} R_n}{{\mathrm d} t} = 
(\frac{\pi_1 \lambda_E}{\mu_E}A + \beta) A - \mu_n R_n.
\end{array}
\label{Paper3_Eq17}
\end{equation}
For simplicity, we set 
$a = \frac{f\tilde{v}\gamma\lambda_E}{\mu_E(\tilde{v}+\mu_G)}-b_1-\mu_A$
and $b = \frac{\pi_1 \lambda_E}{\mu_E}$.
For the stability and bifurcation analysis, 
we choose $\lambda_E$ as the bifurcation parameter.
System~(\ref{Paper3_Eq17}) has a disease-free equilibrium $\bar{E}_0=(0,\,0)$,
which is stable if $a>0$  or 
$\lambda_E>\frac{(b_1+\mu_A)(\tilde{v}+\mu_G)\mu_E}{f\tilde{v}\gamma}$;
and unstable if $a<0$  or 
$\lambda_E<\frac{(b_1+\mu_A)(\tilde{v}+\mu_G)\mu_E}{f\tilde{v}\gamma}$.
Thus a static bifurcation occurs on $\bar{E}_0$ when $a=0$ or 
$\lambda_E=\frac{(b_1+\mu_A)(\tilde{v}+\mu_G)\mu_E}{f\tilde{v}\gamma}$.
The disease equilibrium is given by $\bar{E}_1=(\bar{A},\,\bar{R}_n)$, in which
$\bar{R}_n=\frac{(b\bar{A}+\beta)\bar{A}}{\mu_n}$, and $\bar{A}$ is
given by the roots of
the following equation, 
\begin{equation}
f_8(A) = b\sigma_1 A^2 + \beta \sigma_1 A - \mu_n a.
\label{Paper3_Eq18}
\end{equation}
Equation~(\ref{Paper3_Eq18}) has two roots with negative signs if $a<0$, 
with opposite signs if $a>0$, and only one zero root if $a=0$.
This means that a negative backward bifurcation is possible in 
system~(\ref{Paper3_Eq17}) with proper parameter values.
We further examine the characteristic equation at $\bar{E}_1$, which shares
the same form as equation~(\ref{Paper3_Eq13}), with
$\mathrm{Tr}(J|_{\bar{E}_1})=\frac{1}{\mu_n}
(b\sigma_1 A^2 + \beta \sigma_1 A +\mu_n^2 -a \mu_n) := a_{11}$ and 
$\mathrm{Det}(J|_{\bar{E}_1})=3b \sigma_1 A^2+2\beta \sigma_1 A-a\mu_n:=a_{12}$.
Solving $f_8(A)=0$ and $a_{12}=\mathrm{Det}(J|_{\bar{E}_1})=0$,
gives the static bifurcation point of $\bar{E}_1$ at
$(\bar{A},\,a)=(0,\,0)$ or $(\bar{A},\,\lambda_E)=(0,\,
\frac{(b_1+\mu_A)(\tilde{v}+\mu_G)\mu_E}{f\tilde{v}\gamma})$,
which is a transcritical bifurcation point between $\bar{E}_0$ 
and $\bar{E}_1$.
Moreover, Hopf bifurcation can happen if and only if $f_8(A)=0$ and 
$a_{11}=\mathrm{Tr}(J|_{\bar{E}})=0$, which can 
be satisfied only if $\mu_n=0$. This implies that the positive branch 
of $\bar{E}_1$ is stable for any positive values of $\mu_n$. 
Thus, this model cannot exhibit recurrence, bistability, 
or even regular oscillation. 
The same conclusion was obtained in Zhang et al., (submitted 
for publication) 
for the original 4-dimensional model~(\ref{Paper3_Eq16}). 

However,
a recent experimental discovery~\cite{Baecher2006} 
has revealed a new class of terminally differentiated T$_{\text{Reg}}$ cells.
As described in detail in Zhang et al., (submitted 
for publication), introducing
this cell population, denoted $R_d$, into the model yields the full system
\begin{equation*}
\begin{array}{ll}
\frac{{\mathrm d} A}{{\mathrm d} t} = f\tilde{v}G - \sigma_1(R_n+dR_d)A - (b_1+\mu_A)A \\[0.5ex]
\frac{{\mathrm d} R_n}{{\mathrm d} t}= (\pi_1 E+\beta)A - \mu_n R_n - \xi R_n \\[0.5ex]
\frac{{\mathrm d} R_d}{{\mathrm d} t} = c \xi R_n - \mu_d R_d \\[0.5ex]
\frac{{\mathrm d} E}{{\mathrm d} t} = \lambda_E A - \mu_E E \\[0.5ex]
\frac{{\mathrm d} G}{{\mathrm d} t} = \gamma E - \tilde{v}G - \mu_G G
\end{array}
\end{equation*} 
and quasi-steady state analysis then yields a 
reduced 3-dimensional model in the form
\begin{equation*}
\begin{array}{ll}
\frac{{\mathrm d} A}{{\mathrm d} t} = 
          [\frac{f\tilde{v}\gamma\lambda_E}{(\tilde{v}+\mu_G)\mu_E}
           -(b_1+\mu_A)]A - \sigma_1(R_n+dR_d)A, \\[0.3ex]
\frac{{\mathrm d} R_n}{{\mathrm d} t}= 
          (\frac{\pi_1 \lambda_E}{\mu_E}A+\beta)A 
           - \mu_n R_n - \xi R_n, \\[0.3ex]
\frac{{\mathrm d} R_d}{{\mathrm d} t} = c \xi R_n - \mu_d R_d.
\end{array}
\label{Paper3_Eq19}
\end{equation*}
Again, here $\lambda_E$ is chosen as the bifurcation parameter
for stability and bifurcation analysis. It is easy to show that  
system~(\ref{Paper3_Eq19}) still has 
a disease-free equilibrium $\bar{E}_0$ as $(A,\,R_n,\,R_d)=(0,\,0,\,0)$,
and a disease equilibrium $\bar{E}_1$ as $(\bar{A},\,\bar{R}_n,\,\bar{R}_d)$,
where $\bar{R}_d=\frac{c\xi\bar{R}_n}{\mu_d}$, 
$\bar{R}_n=\frac{\beta\mu_E+\pi_1\lambda_E \bar{A}}{\mu_E(\mu_n+\xi)}\bar{A}$,
and $\bar{A}$ is determined from the following quadratic equation:
\begin{equation*}
f_9(A) = \pi_1\lambda_E A^2 + \beta\mu_E A +
\frac{\mu_d(\mu_n+\xi)}{(\tilde{v}+\mu_G)(cd\xi+\mu_d)\sigma_1}
[-f\gamma\tilde{v}\lambda_E+(b_1+\mu_A)(\mu_G+\tilde{v})\mu_E],
\end{equation*}
which gives two negative roots if 
$\lambda_E<\lambda_{ES}=
\frac{(b_1+\mu_A)(\mu_G+\tilde{v})\mu_E}{f\gamma\tilde{v}}$,
and two roots with opposite signs when
$\lambda_E>\lambda_{ES}$. The critical point is determined 
by $\lambda_E = \lambda_{ES}$, which is actually the intersection 
point of $\bar{E}_0$ and $\bar{E}_1$.
The two equilibrium solutions exchange their stability at 
$\lambda_{ES}$, leading to a transcritical bifurcation at
$(\bar{A},\,\lambda_E)=(0,\,\lambda_{ES})$. 
Note that the negative backward bifurcation still happens 
in system~\ref{Paper3_Eq19}.
Moreover, a Hopf bifurcation occurs from 
the upper branch of $\bar{E}_1$, giving rise to oscillation and recurrence.

\begin{figure}
\vspace{-0.00in} 
\begin{center}  
\hspace{0.30in} 
\begin{overpic}[width=0.35\textwidth
]{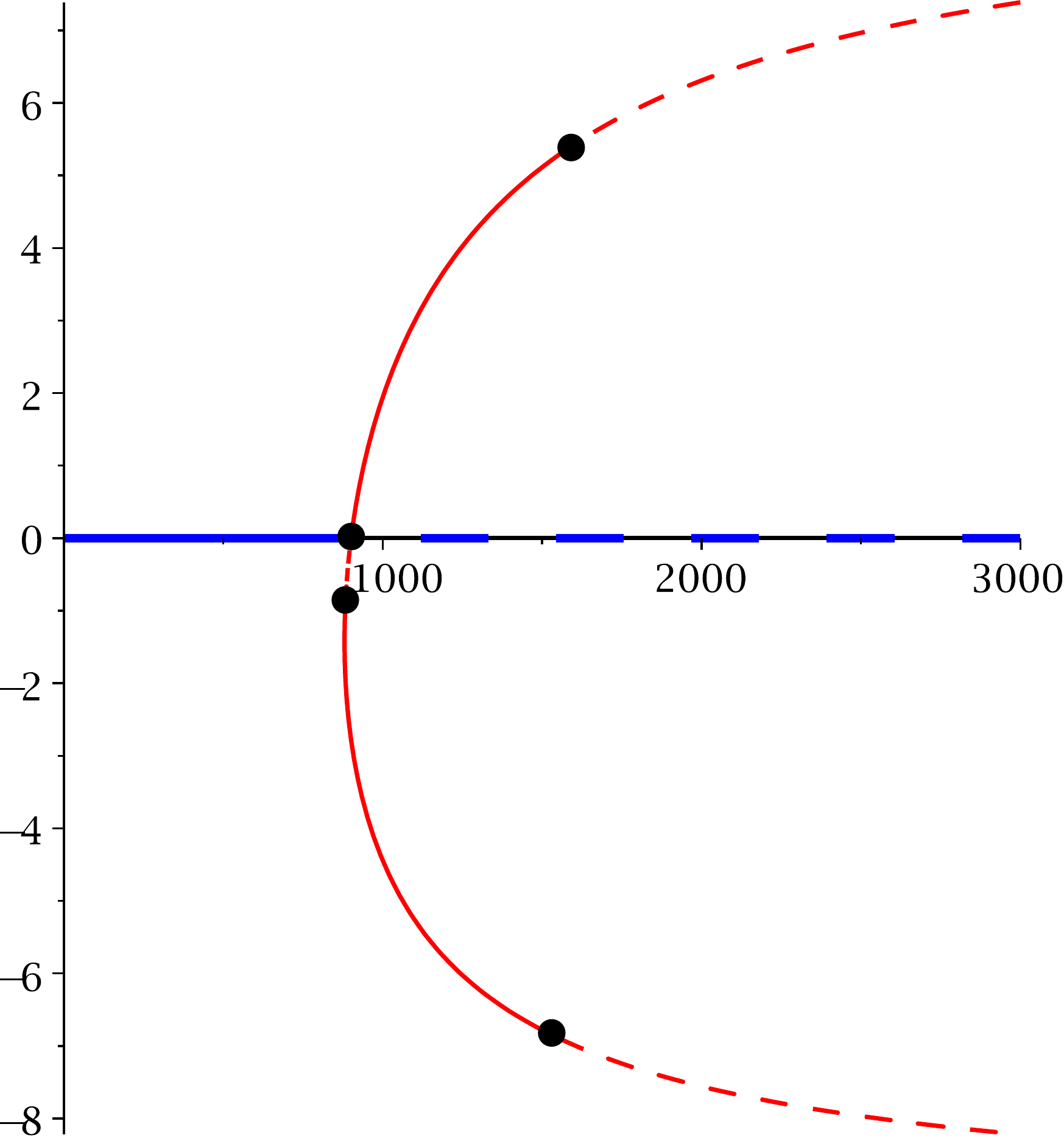}
\put(-4,60){$A$}
\put(75,45){$\lambda_E$}
\put(55,85){Hopf}
\put(60,5){Hopf}
\put(35,55){Transcritical}
\put(7,46){Turning}
\put(50,110){($a$)}
\end{overpic}
\hspace{0.30in}
\begin{overpic}[width=0.50\textwidth,height=0.28\textheight
]{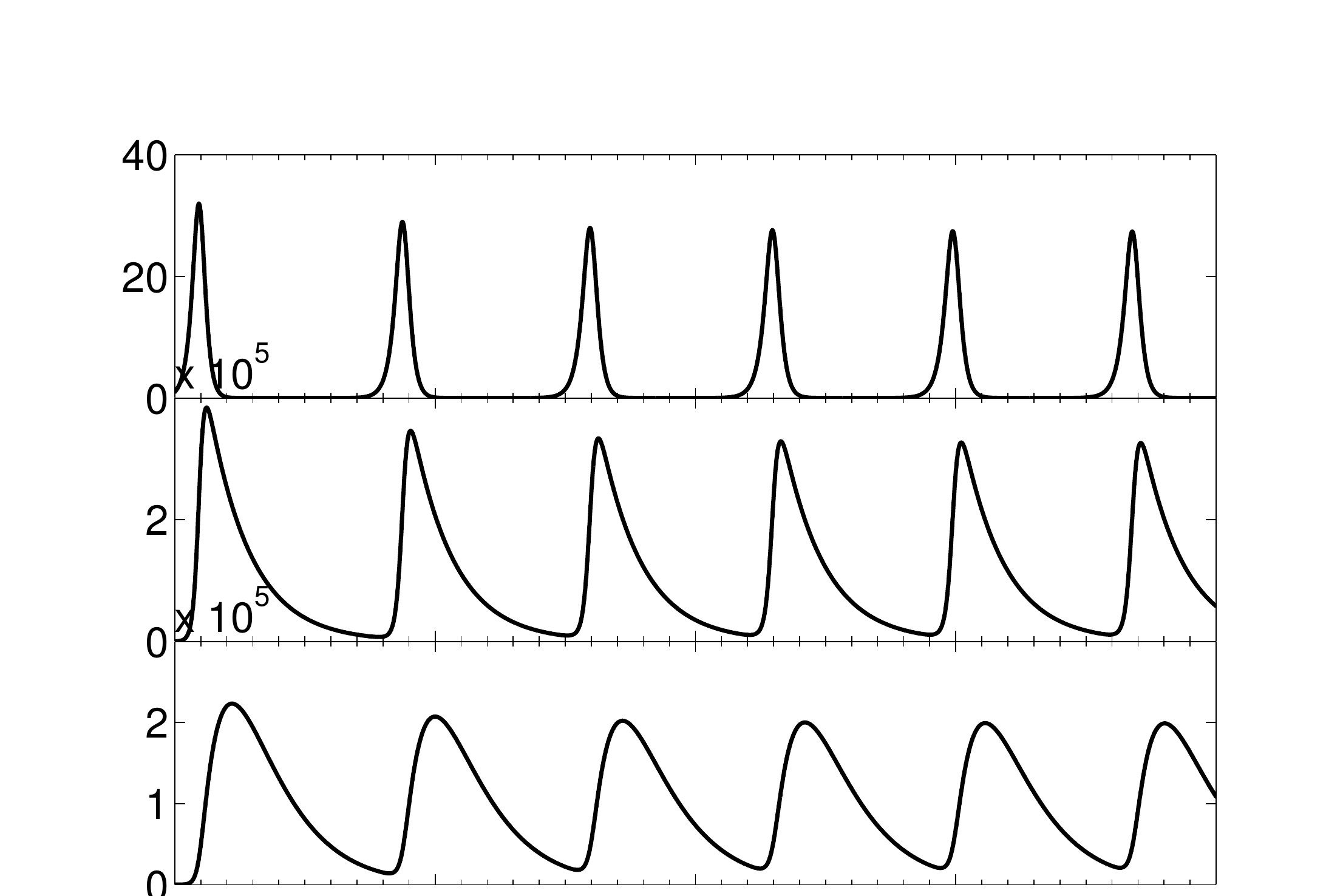}
\put(0,60){$A$}
\put(0,35){$R_n$}
\put(0,10){$R_d$}
\put(50,-8){$t$}
\put(50,82){($b$)}
\end{overpic} \vspace{0.20in} 
\caption{Dynamics of system (\ref{Paper3_Eq17}):
(a) bifurcation diagram; and 
(b) simulated time history for 
$\lambda_E=\lambda_{EH}+1000$.}
\label{Paper3_Fig5}
\end{center} 
\end{figure}

Realistic parameter values have been obtained 
in Zhang et al., (submitted for publication), 
and are given as follows:
\begin{equation*}
\begin{array}{ll} 
f = 1\times 10^{-4},\ \ \tilde{v} = 0.25\times 10^{-2},\ \
\sigma_1 = 3\times 10^{-6},\ \ b_1 = 0.25,\  \mu_A = 0.2,\ \ 
\pi_1 = 0.016,\\[0.5ex]  
\beta = 200,\ \ \mu_n = 0.1,\ \ \mu_E = 0.2,\ \
\gamma = 2000,\ \ \mu_G = 5,\ \ \mu_d =0.2,\ \ c=8,\ \ d =2,\ \ \xi=0.025.
\end{array}
\end{equation*}
For the above parameter values, the Hopf critical point is 
obtained at $(A_H,\,\lambda_{EH})\!=\!(5.6739,$ $1691.6414)$,
while the turning point is at 
$(A_T,\,\lambda_{ET})\!=\!(-1.4205,\,879.9848)$,
and the transcritical bifurcation point is at
$(A_S,\,\lambda_{ES})\!=\!(0,\,900.45)$.
These three bifurcation points and the stability of equilibrium solutions 
are shown in the bifurcation diagram given in Figure~\ref{Paper3_Fig5}(a), 
and the simulated recurrent time history is plotted 
in Figure~\ref{Paper3_Fig5}(b) for $\lambda_E=\lambda_{EH}+1000$.

In summary, when negative backward bifurcation occurs,
that is, the turning point is located in the negative state variable
space, less complex dynamical behavior will be present.
Hopf bifurcation in a biologically feasible area does not happen in 
the reduced 2-dimensional system~(\ref{Paper3_Eq17}), nor in the original 
system~(\ref{Paper3_Eq16})~(Zhang et al., submitted 
for publication). 
However, if we increase the dimension of the system, Hopf bifurcation 
and complex dynamical phenomena can emerge, as 
shown in our results for system~(\ref{Paper3_Eq19}).

\section{Conclusion}
In this paper, we first review previous work on
a reduced 2-dimensional infection model with a concave incidence rate 
\cite{Korobeinikov2005}.
The authors proved that
the disease equilibrium will emerge and be globally stable 
when the basic reproduction number $R_0$ is greater than $1$.
This means that no complex dynamical phenomenon can occur in such 
models.
However, by adding an extra saturating treatment term to this
simple 2-dimensional infection model, the resulting system~(\ref{Paper3_Eq6})
considered in \cite{Zhou2012} can exhibit backward bifurcation,
which increases the parameter range for  
Hopf bifurcation, which in turn leads to recurrent, bistable and 
regular oscillating behaviors. 

Instead of adding an extra term, a 2-dimensional infection model
with a convex incidence function can likewise show rich dynamics 
due to the occurrence of backward bifurcation,
giving rise to two types of Hopf bifurcation. 
Biologically, a convex incidence rate implies that  
existing infection makes the host more vulnerable to further infection, 
showing a cooperative effect in disease progression. 
From the view point of mathematics, the convex incidence function 
enables backward bifurcation to occur on the positive branch of the disease 
equilibrium solution, which further generates Hopf bifurcation.
The location and direction of Hopf bifurcation(s), determined by parameter
values, can further give rise to bistable, recurrent, and regular
oscillating behaviors. 

Cooperative effects also occur during the progression of autoimmune disease.
However, for an autoimmune model with negative backward bifurcation,
in which the turning point is located on the negative state variable
space, the biologically feasible parameter range in which
Hopf bifurcation may occur is limited. 
By introducing an additional state variable 
to the autoimmune model, recurrent phenomenon are once again observed.

\pagestyle{plain}

\end{document}